\documentstyle[12pt]{article}
\font\teneufm=eufm10
\font\seveneufm=eufm7
\font\fiveeufm=eufm5
\newfam\eufmfam
\textfont\eufmfam=\teneufm
\scriptfont\eufmfam=\seveneufm
\scriptscriptfont\eufmfam=\fiveeufm
\def\frak#1{{\fam\eufmfam\relax#1}}

\newfam\msbfam
\font\tenmsb=msbm10 scaled \magstep1  \textfont\msbfam=\tenmsb
\font\sevenmsb=msbm7 scaled \magstep1 \scriptfont\msbfam=\sevenmsb
\font\fivemsb=msbm5 scaled \magstep1  \scriptscriptfont\msbfam=\fivemsb
\def\Bbb{\fam\msbfam \tenmsb}

\def\RR{{\Bbb R}}
\def\CC{{\Bbb C}}

\def\PP{{\Bbb P}}

\def\Im{\hbox{Im}\,}
\def\ra{\rightarrow}

\def\a{\alpha}

\def\HollowBoxx #1#2#3{{\dimen0=#1 \advance\dimen0 by -#2       
       \dimen1=#1 \advance\dimen1 by #3                       
        \vrule height 0pt depth #3 width #2                   
       \hskip -#3
       \vrule height #1 depth #3 width #3}}                   
 \def\LeftContraction{\mathord{\kern1.45pt \HollowBoxx{6pt}{3.5pt}{.4pt}}\,}

 \def\HollowBox #1#2#3{{\dimen0=#1 \advance\dimen0 by -#3       
       \dimen1=#1 \advance\dimen1 by #3                       
        \vrule height #1 depth #3 width #3                    
        \vrule height 0pt depth #3 width #2                   
        \hskip -#3}}                                             
 \def\RightContraction{\mathord{\, \HollowBox{6pt}{3.1pt}{.4pt}} \kern1.6pt}

\newtheorem{theorem}{THEOREM}[section]

\newtheorem{lemma}[theorem]{Lemma}

\newtheorem{remark}[theorem]{REMARK}
\newtheorem{proposition}[theorem]{PROPOSITION}

\begin{document}
\begin{center}
{\Large \bf Invariants of Elliptic and Hyperbolic} 
\medskip \\
{\Large \bf $CR$-Structures of Codimension 2}  
\medskip \\
{\normalsize \rm V. V. Ezhov\footnote{Research supported by the
Australian Research Council.}, \ \ \ A. V. Isaev, \ \ \ G. Schmalz}
\end{center} 

\begin{quotation} 
\small \sl
We reduce $CR$-structures on smooth elliptic and hyperbolic
manifolds of $CR$-codimension 2 to parallelisms thus solving the
problem of global equivalence for such manifolds. The parallelism that we
construct is defined on a sequence of two principal bundles over the manifold, takes values in the Lie algebra of infinitesimal
automorphisms of the quadric corresponding to the Levi form of the manifold, and behaves
``almost'' like a Cartan connection. The construction is explicit and
allows us to study the properties of the parallelism as well as those of its curvature
form. It also leads to a 
natural class of ``semi-flat'' manifolds for which the two bundles
reduce to a single one and the parallelism
turns into a true Cartan connection.

In addition, for real-analytic manifolds we describe certain
local normal forms that do not require passing to bundles, but in many ways agree with the structure of the parallelism.     
\end{quotation}

\markboth{V. V. Ezhov, A. V. Isaev, G. Schmalz}{CR-structures of
codimension 2}

\setcounter{section}{-1}

\section{Introduction}

We start with a brief overview of necessary definitions and facts from
$CR$-geometry (see e.g. \cite{Tu1} for a more detailed exposition). 

A $CR$-{\it structure} on a smooth  real connected  manifold $M$ of dimension $m$ is a smooth distribution of
subspaces in the tangent spaces $T_p^c(M)\subset T_p(M)$, $p\in M$,
with operators of complex structure $J_p:T_p^c(M)\ra T_p^c(M)$,
$J_p^2\equiv -\hbox{id}$, that depend smoothly on $p$. A manifold $M$
equipped with a $CR$-structure is called a $CR$-{\it manifold}. It
follows that the number
$CR\hbox{dim}M:=\hbox{dim}_{\CC}T_p^c(M)$ does not depend on
$p$; it is called the $CR$-{\it dimension} of  $M$. The number
$CR\hbox{codim}M:=m-2CR\hbox{dim}M$ is called the $CR$-{\it
codimension} of $M$. $CR$-structures naturally arise on real submanifolds in complex manifolds. Indeed, if, for example,
$M$ is a real submanifold of $\CC^K$, then one can define the distribution
$T_p^c(M)$ as follows:
$$
T_p^c(M):=T_p(M)\cap iT_p(M).
$$
On each $T_p^c(M)$ the operator $J_p$ is then defined as the
operator of multiplication by $i$. Then $\{T_p^c(M), J_p\}_{p\in M}$
form a $CR$-structure on $M$, if $\dim_{\CC}T_p^c(M)$ is constant. This is always the case, for example, if $M$ is a real
hypersurface in $\CC^K$ (in which case $CR\hbox{codim}M=1$). We say that such a $CR$-structure is {\it
induced by $\CC^K$}.

A mapping between two $CR$-manifolds $f:M_1\ra M_2$ is called a
$CR$-{\it mapping}, if for every $p\in M_1$: {\bf (i)} $df(p)$ maps
$T_p^c(M_1)$ to $T_{f(p)}^c(M_2)$, and {\bf (ii)} $df(p)$ is complex linear
on $T_p^c(M_1)$. Two $CR$-manifolds $M_1$, $M_2$ are called $CR$-{\it
equivalent}, if there is a $CR$-diffeomorphism from $M_1$ onto
$M_2$. Such a $CR$-diffeomorphism $f$ is called a $CR$-{\it isomorphism}.

In this paper we are interested in the equivalence problem for $CR$-manifolds. This problem can be viewed as a special case of the
equivalence problem for $G$-{\it structures}.
Let $G\subset GL(m,\RR)$ be a subgroup. A $G$-structure on an
$m$-dimensional manifold $M$ is a subbundle $Q$ of the frame bundle $F(M)$ over
$M$ which is a principal $G$-bundle over $M$. Two $G$-structures
$Q_1$, $Q_2$ on manifolds $M_1$, $M_2$ respectively are called {\it
equivalent}, if there is a diffeomorphism $f$ from $M_1$ onto $M_2$
such that the induced mapping $f_{*}:F(M_1)\ra F(M_2)$ maps $Q_1$ onto
$Q_2$. Such a diffeomorphism $f$ is called an {\it isomorphism of the
$G$-structures}.
A $CR$-structure on a manifold $M$ of $CR$-dimension $n$ and
$CR$-codimension $k$ (so that $m=2n+k$) is a $G$-structure, where $G$
is the group of linear transformations of $\CC^n\oplus\RR^k$ that
preserve the first component and are complex linear on it. The notion of
equivalence of such $G$-structures is then exactly that of
$CR$-structures.

\'E. Cartan developed a general approach to the equivalence problem for
$G$-structures (see \cite{St}) that works for many important examples
of $G$-structures. We will be looking for a solution to
the equivalence problem in the spirit of Cartan's work. Namely, we
will be trying to reduce the $CR$-structure in consideration to an
$\{e\}$-structure, or {\it parallelism}, where $\{e\}$ is the one-element group. An
$\{e\}$-structure on an $N$-dimensional manifold $P$ is given by a
1-form $\omega$ on $P$ with values in $\RR^N$ such that for every
$x\in P$, $\omega(x)$ is an isomorphism of $T_x(P)$ onto
$\RR^N$. The problem of local equivalence for generic
parallelisms is well-understood (see \cite{St}).  

Let ${\cal C}$ be a class of manifolds equipped with
$G$-structures. We say that $G$-structures on manifolds from ${\cal C}$ are
{\it $s$-reducible to
parallelisms} if for any $M\in{\cal C}$  there is a sequence of principle bundles
$$
P^s\stackrel{\pi^s}{\rightarrow}\dots\stackrel{\pi^3}{\rightarrow}P^2\stackrel{\pi^2}{\rightarrow}P^1\stackrel{\pi^1}{\rightarrow}M
$$
and a parallelism $\omega$ on $P^s$ such that:
\smallskip\\

\noindent {\bf (i)} Any isomorphism of $G$-structures $f:M_1\ra M_2$
for $M_1,M_2\in{\cal C}$ can be lifted to a diffeomorphism $F: P^s_1\ra
P^s_2$ such that $F^{*}\omega_2=\omega_1$;
\smallskip\\

\noindent {\bf (ii)} Any diffeomorphism $F: P^s_1\ra P^s_2$ such that
$F^{*}\omega_2=\omega_1$, is a lift of an isomorphism of the
$G$-structures $f:M_1\ra M_2$, for $M_1,M_2\in{\cal C}$.
\smallskip\\

\noindent In the above definition we say that $F$ is a {\it lift} of $f$ if
$$
\pi_2^1\circ\dots\circ\pi_2^s\circ
F=f\circ\pi_1^1\circ\dots\circ\pi_1^s.
$$

From now on we will concentrate on solving the equivalence problem, in
the sense of reducing to parallelisms, for $CR$-structures of certain classes ${\cal C}$
that we are now going to introduce. Let $M$ be a $CR$-manifold. For
every $p\in M$ consider the complexification
$T_p^c(M)\otimes_{\RR}\CC$. Clearly, this complexification can be
represented as the direct sum
$$
T_p^c(M)\otimes_{\RR}\CC=T_p^{(1,0)}(M)\oplus T_p^{(0,1)}(M),
$$
where
\begin{eqnarray*}
T_p^{(1,0)}(M)&:=&\{X-iJ_pX:X\in T_p^c(M)\},\\
T_p^{(0,1)}(M)&:=&\{X+iJ_pX:X\in T_p^c(M)\}.
\end{eqnarray*}
The $CR$-structure on $M$ is called {\it integrable} if for any local
sections $Z,Z'$ of the bundle $T^{(1,0)}(M)$, the vector field $[Z,Z']$
is also a section of $T^{(1,0)}(M)$. It is not difficult to see that
if $M\subset\CC^K$ and the $CR$-structure on $M$ is induced by
$\CC^K$, then it is integrable. In this paper all $CR$-structures are
assumed to be integrable.

An important characteristic of a $CR$-structure called the {\it Levi
form} comes from taking commutators of local sections of
$T^{(1,0)}(M)$ and $T^{(0,1)}(M)$. Let $p\in M$, $z, z'\in
T_p^{(1,0)}(M)$, and $Z$, $Z'$ be local sections of $T^{(1,0)}(M)$ near
$p$ such that $Z(p)=z$, $Z'(p)=z'$. The Levi form of $M$ at
$p$ is the Hermitian form on $T_p^{(1,0)}(M)$ with values in $(T_p(M)/T_p^c(M))\otimes_{\RR}\CC$ given by
$$
{\cal L}_M(p)(z,z'):=i[Z,\overline{Z'}](p)(\hbox{mod}\,T_p^c(M)\otimes_{\RR}\CC).
$$
The Levi form is defined uniquely up to the choice of coordinates in
$(T_p(M)/T_p^c(M))\otimes_{\RR}\CC$, and, for fixed $z$ and $z'$,
its value  does not depend on the choice of $Z$ and $Z'$.

Let $H=(H_1,\dots,H_k):\CC^n\times\CC^n\ra
\RR^k$ be a $\RR^k$-valued Hermitian form on $\CC^n$. We say that $H$ is
{\it non-degenerate} if:
\smallskip\\

\noindent{\bf (i)} The scalar Hermitian forms $H_1,\dots,H_k$ are
linearly independent over $\RR$;
\smallskip\\

\noindent{\bf (ii)} $H(z,z')=0$ for all $z'\in\CC^n$ implies $z=0$.
\smallskip\\

A $CR$-structure on $M$ is
called {\it Levi non-degenerate}, if its
Levi form at any $p\in M$ is non-degenerate. In this paper we consider
only Levi-nondegenerate manifolds.

For any Hermitian form $H$ there is a corresponding standard $CR$-manifold
$Q_H\subset\CC^{n+k}$ of $CR$-dimension $n$ and $CR$-codimension
$k$ defined as follows:
$$
Q_H:=\{(z,w):\Im w=H(z,z)\},
$$
where $z:=(z_1,\dots,z_n)$, $w:=(w_1,\dots,w_k)$ are coordinates in
$\CC^{n+k}$.
The manifold $Q_H$ is often called the {\it quadric associated with
the form $H$}. Clearly, the Levi form of $Q_H$ at any point is given by
$H$. 

An important tool in all the considerations below is
the {\it automorphism group of $Q_H$}. Let $\hbox{Aut}(Q_H)$ denote
the collection of all local $CR$-isomorphisms of $Q_H$ to itself that
we call {\it local $CR$-automorphisms}. It
turns out that, if $H$ is non-degenerate, then any
$f\in\hbox{Aut}(Q_H)$ extends to a rational (more precisely, a matrix
fractional quadratic) map of $\CC^{n+k}$ \cite{KT},
\cite{F}, \cite{Tu2}, \cite{ES5}. Thus, for a non-degenerate $H$, $\hbox{Aut}(Q_H)$
is a finite-dimensional Lie group, and we denote by
$\hbox{Aut}_{e}(Q_H)$ its identity component. Note that $Q_H$ is a
homogeneous manifold since the global $CR$-automorphisms
\begin{eqnarray*}
z&\mapsto& z+z^0,\\
w&\mapsto&w+w^0+2iH(z,z^0), \qquad \qquad \qquad \qquad (0.1)
\end{eqnarray*}
for $(z^0,w^0)\in Q_H$, act transitively on $Q_H$. Therefore, it is
important to consider the isotropy group of a point of $Q_H$, say the
origin, i.e. the group of local $CR$-automorphisms of $Q_H$ preserving
the origin. We denote this subgroup of
$\hbox{Aut}(Q_H)$ by $\hbox{Aut}_0(Q_H)$ and its identity
component by $\hbox{Aut}_{0,e}(Q_H)$. We also introduce the group
$\hbox{Aut}_{lin}(Q_H)\subset\hbox{Aut}_0(Q_H)$ of linear
automorphisms of $Q_H$ and its identity
component $\hbox{Aut}_{lin,e}(Q_H)$. All these groups are
finite-dimensional Lie groups. We call a Levi
non-degenerate $CR$-manifold $M$ {\it weakly uniform}, if for any pair
of points $p,q\in M$,
the groups $\hbox{Aut}_{lin,e}(Q_{{\cal L}_M(p)})$, $\hbox{Aut}_{lin,e}(Q_{{\cal L}_M(q)})$  are isomorphic,
and the isomorphism extends to an isomorphism between
$\hbox{Aut}_{0,e}(Q_{{\cal L}_M(p)})$ and $\hbox{Aut}_{0,e}(Q_{{\cal L}_M(q)})$.

Let $H^1,H^2$ be two $\RR^k$-valued Hermitian forms on $\CC^n$. We say
that $H^1$ and $H^2$ are {\it equivalent}, if there exist linear
transformations $A$ of $\CC^n$ and $B$ of $\RR^k$ such that
$$
H^2(z,z)=BH^1(Az,Az)
$$
for all $z\in\CC^n$. We call a $CR$-manifold $M$ {\it strongly uniform},
if the forms ${\cal L}_M(p)$ are equivalent for all $p\in M$. If, for
example, $M$ is Levi non-degenerate and $CR\hbox{codim}M=1$ then $M$
is strongly uniform. Clearly,
for a Levi non-degenerate $M$, strong uniformity implies weak
uniformity.

Existing results on the equivalence problem for
$CR$-structures treat two large classes of Levi-nondegenerate
manifolds: {\bf (i)} strongly uniform manifolds
and {\bf (ii)} weakly uniform
manifolds, for which, in addition, the groups
$\hbox{Aut}_0(Q_{{\cal L}_M(p)})$ are ``sufficiently small''; in
particular, $\hbox{Aut}_{0}(Q_{{\cal L}_M(p)})=\hbox{Aut}_{lin}(Q_{{\cal L}_M(p)})$.

\'E. Cartan \cite{Ca} solved the problem for 3-dimensional
Levi-nondegenerate $CR$-manifolds of $CR$-dimension 1 (such manifolds
are, of course, strongly uniform). Tanaka in 1967 obtained a solution
for all Levi-nondegenerate strongly uniform manifolds \cite{Ta1}, but his result
became widely known only after the Chern-Moser work in 1974 \cite{CM} where the
problem was solved independently for Levi-nondegenerate manifolds of
$CR$-codimension 1 (see also \cite{Ta2}, \cite{J}). We note that although Tanaka's pioneering construction is very
important and applies
to very general situations (that include also geometric structures
other than $CR$-structures), his treatment of the case of
$CR$-codimension 1 is far
less detailed and clear than that due to Chern and Moser (see \cite{K} for a
discussion of this matter). For example, Tanaka's
construction gives 3-reducibility to parallelisms, whereas Chern's
original construction gives 2-reducibility and
even 1-reducibility \cite{BS}. The structure group of the single
bundle $P\rightarrow M$ that arises in Chern's construction, is
$\hbox{Aut}_{0,e}(Q_H)$ (or, alternatively, $\hbox{Aut}_0(Q_H)$), where $H$ is a Hermitian
form equivalent to any of ${\cal L}_M(p)$, $p\in M$, and the
parallelism $\omega$ takes values in the Lie algebra of infinitesimal
automorphisms of $Q_H$ (we denote it by ${\frak g}_H$). This algebra is the Lie algebra of the group
$\hbox{Aut}(Q_H)$ and is well-understood (see \cite{Sa}, \cite{B1}); in
particular, ${\frak g}_H$ is a graded Lie algebra: ${\frak
g}_H=\oplus_{k=-2}^2{\frak g}^k_H$, where the components ${\frak
g}^1_H$, ${\frak g}^2_H$ are responsible for non-linear automorphisms.
In Tanaka's
construction, however, the parallelism takes values in a certain maximal 
prolongation $\tilde{{\frak g}}_H$ of $\oplus_{k=-2}^{0}{\frak
g}^{k}_H$; the
coincidence of $\tilde{{\frak g}}_H$ and ${\frak g}_H$ is not in
general 
established (see Section 5 for a discussion). Further,
it is shown in \cite{CM} (see also \cite{BS}) that the parallelism
$\omega$ from Chern's construction is in fact a {\it Cartan
connection}, i.e. changes in a regular way under the action of the
structure group of the bundle. Namely, if for $\eta\in\hbox{Aut}_{0,e}(Q_H)$, $L_{\eta}$
denotes the (left) action by $\eta\in\hbox{Aut}_{0,e}(Q_H)$ on $P$,
then $L_{\eta}^{*}\omega=\hbox{Ad}(\eta)\omega$, where $\hbox{Ad}$ is the
adjoint representation of $\hbox{Aut}_{e}(Q_H)$ on ${\frak g}_H$. It is
not clear from \cite{Ta1}  (even in the case of $CR$-codimension 1) whether the sequence of bundles $P^3\ra
P^2\ra P^1\ra M$ there can be
reduced to a single bundle and whether the
parallelism defined on $P^3$ behaves in any sense like a Cartan
connection (see, however, later work in \cite{Ta2}, \cite{Ta3}). Being more detailed, Chern's construction also
allows one to investigate the important {\it curvature form of $\omega$}, i.e. the
2-form $\Omega:=d\omega-\frac{1}{2}[\omega,\omega]$ and find its
precise expansion. It also can be used
to introduce special invariant curves on the manifold called {\it chains} that
turned out to be very important in the study of real hypersurfaces in
$\CC^K$. These and other differences between Tanaka's and Chern's
construction motivated our work. 

The results in \cite{Ca}, \cite{Ta1}, \cite{CM}, in particular, solve
the equivalence problem for Levi-nondegenerate $CR$-manifolds of
$CR$-codimension 1, thus we concentrate on the case of higher
$CR$-codimensions. Certain Levi-nondegenerate weakly uniform $CR$-structures of
codimension 2 were treated in \cite{La}, \cite{M}. The conditions
imposed on the Levi form in these papers are stronger than non-degeneracy and
force the groups $\hbox{Aut}_0(Q_{{\cal L}_M(p)})$, $p\in
M$, to be 
minimal possible; in
particular, they contain only linear transformations of a special
form (this of course implies that ${\frak g}^k_{{\cal L}_M(p)}=0$ for $k=1,2$). A similar situation for Levi non-degenerate manifolds $M$ with
$CR\hbox{codim}M>2$ and the additional condition $CR\hbox{dim}M>(CR\hbox{codim}M)^2$ was
treated in \cite{GM}. A motivation to consider manifolds with Levi
form satisfying conditions as in \cite{M} (for $CR\hbox{dim}M\ge 7$), \cite{La}, \cite{GM} is
that, in the considered situations, these conditions are {\it stable}, i.e., if they are
satisfied at a single point $p$ of a manifold $M$, they are also satisfied
on an neighbourhood of $p$. Moreover, quadrics associated with Levi
forms as in \cite{M} (for $CR\hbox{dim}M\ge 7$), \cite{La} are dense in
the space of all Levi non-degenerate quadrics.

In this paper we consider the case
$CR\hbox{dim}M=CR\hbox{codim}M=2$. This is one of two exceptional
cases among all
$CR$-structures with 
$CR\hbox{codim}M>1$ in the following sense: typically (in fact, always
except for
the cases $CR\hbox{dim}M=CR\hbox{codim}M=2$ and
$(CR\hbox{dim}M)^2=CR\hbox{codim}M$) for generic non-degenerate Hermitian forms, the corresponding quadrics
have only linear automorphisms \cite{M}, \cite{B2}, \cite{ES6}. In the case that we consider, however, generic
quadrics always have plenty of non-linear automorphisms. 
Any non-degenerate Hermitian
form $H=(H_1,H_2)$ on $\CC^2$ is equivalent to one of the following:
\begin{eqnarray*}
H^1(z,z)&:=&(|z_1|^2+|z^2|^2, z_1\overline{z_2}+z_2\overline{z_1}),\\
H^{-1}(z,z)&:=&(|z_1|^2-|z_2|^2, z_1\overline{z_2}+z_2\overline{z_1}),\\
H^0(z,z)&:=&(|z_1|^2, z_1\overline{z_2}+z_2\overline{z_1}).
\end{eqnarray*}
These forms are called respectively {\it hyperbolic}, {\it elliptic}
and {\it parabolic}. We will be interested in the case of strongly
uniform $CR$-manifolds whose Levi form is either at every point
equivalent to the hyperbolic form, or it is at every point
equivalent to the elliptic form. We will call such manifolds {\it
hyperbolic} and {\it elliptic} respectively. Clearly, the conditions of
hyperbolicity and ellipticity are stable: if the Levi form of a
$CR$-manifold $M$ at $p\in M$ is equivalent to the hyperbolic or
elliptic form, then there is a neighbourhood of $p$ which is respectively a
hyperbolic or elliptic manifold. Moreover, the collection of all
hyperbolic and elliptic quadrics is an open dense subset in the space of
all Levi non-degenerate quadrics of $CR$-codimension and $CR$-dimension
2. We denote the sets of all hyperbolic and elliptic manifolds by
${\cal C}^{1}$ and ${\cal C}^{-1}$ respectively. 

The equivalence problem for hyperbolic and elliptic $CR$-manifolds is,
of course, covered by Tanaka's construction in
\cite{Ta1}. Therefore, our main goal is to produce, in this particular
case, a construction different from that in \cite{Ta1}, such that it would be
more detailed and easier to use
in applications. We achieve this by
following the main steps of Chern's reduction
in \cite{CM}, although a great many things will
have to be treated differently. Although we study just manifolds with
$CR\hbox{dim}M=CR\hbox{codim}M=2$, the considered case possesses a
rich structure: the groups
$\hbox{Aut}_0(Q_{H^{\delta}})$ are large in the sense
that they contain substantial non-linear
part (here $\hbox{dim}({\frak g}^1_{H^{\delta}}\oplus{\frak
g}^2_{H^{\delta}})=6$ \cite{ES1}). 

We will formulate our main
result in Section 1 and discuss it in Section 3; here we list just a
few features of our
construction and its applications:
\smallskip\\

\noindent {\bf (i)} We obtain {\it 2-reducibility} to parallelisms, i.e.
sequences of two principal bundles $P^{2,\delta}\ra P^{1,\delta}\ra M$
for $M\in{\cal C}^{\delta}$.
\smallskip\\

\noindent {\bf (ii)} The structure groups of $P^{1,\delta}$, $P^{2,\delta}$ are simply
described groups $G^{1,\delta}$, $G^{2,\delta}$, where $G^{2,\delta}$ is a subgroup of
$\hbox{Aut}_{0,e}(Q_{H^{\delta}})$, 
whereas in Tanaka's constructions the
structure groups on each step are found as certain
very special factor-groups of subgroups of $\hbox{Aut}_{0,e}(Q_{H^{\delta}})$. 
\smallskip\\

\noindent {\bf (iii)} The parallelism $\omega^{\delta}$ defined on $P^{2,\delta}$
takes values in ${\frak g}_{H^{\delta}}$ rather than in the abstractly defined
Lie algebra $\tilde{{\frak g}}_{H^{\delta}}$ as in \cite{Ta1}.
\smallskip\\

\noindent {\bf (iv)} There is an explicit transformation formula for
$\omega^{\delta}$ under the action of $G^{2,\delta}$ on $P^{2,\delta}$
that 
shows that $\omega^{\delta}$ is
not ``too far'' from being a Cartan connection. We also explicitly
calculate the obstructions for $\omega^{\delta}$ to be a Cartan connection. The
obstructions are given by two scalar $CR$-invariants
(i.e. $CR$-invariant functions) on $P^{2,\delta}$, and
we study manifolds for which these invariants vanish; we term such
manifolds {\it semi-flat}. It turns out that the invariant theory on semi-flat
manifolds is completely analogous to that in the case of
$CR$-dimension and 
$CR$-codimension 1, if one substitutes in all the formulas scalars by
matrices from a certain commutative algebra.
\smallskip\\

\noindent {\bf (v)} We calculate precisely the obstruction to
1-reducibility, that is, we can say when the sequence of two bundles
$P^{2,\delta}\ra P^{1,\delta}\ra M$ can be reduced to a single bundle with structure
group $\hbox{Aut}_{0,e}(Q_{H^{\delta}})$. The obstruction is given by a single scalar $CR$-invariant
on $P^{2,\delta}$.
\smallskip\\

\noindent {\bf (vi)} We obtain exact expansions for the curvature form of
$\omega^{\delta}$ in terms of the components of $\omega^{\delta}$. This allows us, for example,
to describe all $CR$-flat manifolds in much the same way as in the case of
$CR$-codimension 1: all such manifolds must be locally $CR$-equivalent
to $Q_{H^{\delta}}$.
\smallskip\\

There is one more issue that does not seem to be tractable from
Tanaka's construction and that in fact was a starting point for our
work. Namely, we are interested in finding analogues of chains for
$CR$-manifolds of $CR$-codimension $>1$. In the case of
$CR$-codimension 1, chains arise in \cite{CM}
independently in the geometric construction as well as in the
construction of the analytic normal form
for a defining function of a real-analytic hypersurface in
$\CC^K$. In the geometric approach chains are the projections to $M$ of the
integral manifolds of a certain distribution on $P$ that consists of parallel
subspaces with respect to the parallelism. In the analytic approach chains
are certain curves that locally become straight lines of a special form in
normal coordinates. For some classes of real-analytic $CR$-submanifolds of
$\CC^K$ of $CR$-codimension $>1$  analogues of the Chern-Moser normal
forms have been found in \cite{Lo1}, \cite{ES2}--\cite{ES4}.
These
normal forms have led to certain analogues of chains that are submanifolds of
$M$ of dimension equal to $CR\hbox{codim}M$. However, they do not
inherit all the properties of chains in the case
$CR\hbox{codim}M=1$. In particular, translations along such chains do
not preserve all conditions of the normal forms; in other words, such chains
can be regarded as proper chains only at a single point 
(the center of normalization).  
The first motivation for the present work was the fact that we did not
have any reasonable explanation to this phenomenon. Our initial idea
was to construct proper analogues of chains (or to understand why
such construction is impossible) by using a reduction of the $CR$-structure to
parallelisms rather than normal forms. Our approach to some extend clarifies the
matter. Namely:
\smallskip\\

\noindent {\bf (vii)} The construction leads to a certain distribution
on $P^{2,\delta}$ (that we call the {\it chain distribution}) which is analogous to
Chern's distribution. However, this distribution is {\it not} in general involutive, and
thus does not in general have integral manifolds. It is worth noting that the
obstructions for the distribution to be integrable exactly coincide
with those for the parallelism to be a Cartan connection. In
particular, the distribution gives
proper chains (that we call $G$-{\it chains}) on semi-flat manifolds.
\smallskip\\

Thus, the parallelism in general
does not generate proper chains. However, there are many
submanifolds of $P^{2,\delta}$ whose tangent space
at a given point is an element of the chain distribution. Most likely,
the
projections of a family of such submanifolds to $M$ are the chains arising in
\cite{Lo1}, \cite{ES2} and thus are ``chains at a single
point''.
It may be that in applications one should be content with considering
the chain distribution
itself without trying to integrate it, that is, with considering only
``infinitesimal chains''.

The paper is organized as follows. In Section 1 we collect all
necessary facts concerning the groups $\hbox{Aut}_{e}(Q_{H^{\delta}})$,
$\hbox{Aut}_0(Q_{H^{\delta}})$,
$\hbox{Aut}_{0,e}(Q_{H^{\delta}})$ and the algebra ${\frak g}_{H^{\delta}}$, $\delta=\pm 1$, and formulate
our main result (Theorem 1.1). We prove Theorem 1.1 in Section 2. In
Section 3 we discuss some corollaries of Theorem 1.1 and
applications of the
construction used in its proof; in particular, we introduce semi-flat
manifolds as manifolds for which the curvature form of the
parallelism behaves in some sense like that in the case of
$CR$-codimension 1. In the real-analytic case, we also introduce
so-called {\it matrix surfaces} as submanifolds of $\CC^4$ whose
defining functions are given by
power series of a
special form. Matrix surfaces are examples of semi-flat
manifolds, and it is very likely that semi-flat manifolds in the 
real-analytic case locally coincide with matrix surfaces. We conclude
Section 
3 by proving this statement for hyperbolic manifolds.
In
Section 4 we reintroduce local normal forms for defining functions of
real-analytic hyperbolic and elliptic $CR$-manifolds in $\CC^4$ that
are certain interpretations of the normal forms constructed in \cite{Lo1}, \cite{ES2}. These normal
forms in many ways agree with our reduction process of $CR$-structures to
parallelisms in the proof of Theorem 1.1. In particular, such normal
forms independently define proper
chains on matrix surfaces, and it turns out that these chains
coincide with $G$-chains. We conclude the paper with Section 5 where
we discuss some questions that arose during our work and
that we consider important for a better understanding of
high-codimensional $CR$-structures.

We would like to thank N. Kruzhilin and A. Tumanov for stimulating
discussions and for showing us useful references. A significant part of this work was done while the second author was
an Alexander von Humboldt fellow at the University of Wuppertal and
the third author was visiting the University of Adelaide. The
work was completed while the first author was visiting the Centre for
Mathematics and Its Applications, The Australian National University.

\section{Preliminaries and Formulation of the Main Result}

Before we formulate our main result, we will discuss the structure of the groups $\hbox{Aut}_{e}(Q_{H^{\delta}})$,
$\hbox{Aut}_0(Q_{H^{\delta}})$,
$\hbox{Aut}_{0,e}(Q_{H^{\delta}})$ and the algebra ${\frak
g}_{H^{\delta}}$, $\delta=\pm 1$.

The groups $\hbox{Aut}_{0,e}(Q_{H^{\delta}})$ were explicitly found in \cite{ES1} (see
also \cite{B2} for the case $\delta=1$). One of the possible
interpretations of the descriptions in \cite{ES1}, \cite{B2} is as
follows. Denote by ${\frak A}^{\delta}$ the commutative algebra of
matrices of the form
$
\left(
\begin{array}{ccc}
   \displaystyle a&\delta b\\
   \displaystyle b&a\\
\end{array}
\right),
$
where $a,b\in\CC$. Let $SU^{\delta}(2,1)$ be the group of
$3\times3$ matrices $U$ with elements from ${\frak A}^{\delta}$ such
that
$$
U\left(
\begin{array}{ccc}
\displaystyle 0&0&-\frac{i}{2}E\\
\displaystyle 0&E&0\\
\displaystyle \frac{i}{2}E&0&0\\
\end{array}
\right)\overline{U}^T=
\left(
\begin{array}{ccc}
\displaystyle 0&0&-\frac{i}{2}E\\
\displaystyle 0&E&0\\
\displaystyle \frac{i}{2}E&0&0\\
\end{array}
\right),
$$
where $E$ is the $2\times2$ identity matrix, and such that
$\hbox{det}\,U=E$. Let ${\frak A}^{\delta *}$ denote the set of
invertible elements in ${\frak A}^{\delta}$, $\hbox{Re}\,{\frak A}^{\delta}$
the set of elements of ${\frak A}^{\delta}$ with real entries, and by 
$\hbox{Re}\,{\frak A}^{\delta *}$ the set of invertible elements in
$\hbox{Re}\,{\frak A}^{\delta}$. It is shown in \cite{ES1} that any
element of $\hbox{Aut}_{0,e}(Q_{H^{\delta}})$ can be viewed as
a transformation of the form
$$
\left(
\begin{array}{ccc}
   \displaystyle E&0&0\\
   \displaystyle -2i\overline{A}&C&0\\
\displaystyle -(R+iA\overline{A})&CA&C\overline{C}\\
\end{array}
\right), \eqno{(1.1)}
$$
of the
``projective space'' ${\frak A}^{\delta}\PP^3:=({\frak
A}^{\delta}\oplus
{\frak A}^{\delta}\oplus{\frak A}^{\delta})/{\frak A}^{\delta *}$,
where $A\in{\frak A}^{\delta}$, $C\in{\frak
A}^{\delta *}$, $R\in\hbox{Re}{\frak
A}^{\delta}$. To make a matrix of the form (1.1) belong to
$SU^{\delta}(2,1)$ we need to multiply it by $\sigma\in{\frak
A}^{\delta *}$ such that $\sigma\overline{\sigma}C\overline{C}=E$ and
$\sigma^3C^2\overline{C}=E$. Note that this does not
change mapping (1.1) as a transformation of ${\frak A}^{\delta}\PP^3$.  
It is not difficult to show that such a
$\sigma$ always exists and is unique up to multiplication by
$\lambda\in{\frak A}^{\delta *}$ with $\lambda\overline{\lambda}=E$,
$\lambda^3=E$. We denote the set of such $\lambda$'s by 
${\hat{\frak A}}^{\delta}$. A straightforward computation gives that
$$
{\hat{\frak A}}^{-1}=\left\{aE:a^3=1\right\},
$$
$$
{\hat{\frak A}}^1=\left\{aE:a^3=1\right\}\bigcup
\left\{
a\left(
\begin{array}{ccc}
\displaystyle 1&\pm i\sqrt{3}\\
\displaystyle \pm i\sqrt{3}&1\\
\end{array}
\right): a^3=-\frac{1}{8}\right\}.
$$
Therefore, $\hbox{Aut}_{0,e}(Q_{H^{\delta}})$ is isomorphic to the subgroup of
$SU^{\delta}(2,1)$ of matrices of the form
$$
\left(
\begin{array}{ccc}
   \displaystyle \sigma&0&0\\
   \displaystyle -2i\sigma\overline{A}&\sigma C&0\\
\displaystyle -\sigma(R+iA\overline{A})&\sigma CA&\sigma C\overline{C}\\
\end{array}
\right), \eqno{(1.2)}
$$
with $A\in{\frak A}^{\delta}$, $C,\sigma\in{\frak A}^{\delta *}$, $R\in\hbox{Re}\,{\frak
A}^{\delta}$, $\sigma\overline{\sigma}C\overline{C}=E$, $\sigma^3C^2\overline{C}=E$,
factorized by the subgroup $Z^{\delta}$ of matrices
$$
\lambda\left(
\begin{array}{ccc}
   \displaystyle E&0&0\\
   \displaystyle 0&E&0\\
\displaystyle 0&0&E\\
\end{array}
\right),
$$
with $\lambda\in\hat{\frak A}^{\delta}$.
Note that $Z^{\delta}$ is a discrete subgroup of $SU^{\delta}(2,1)$.

Analogously, the group of
transformations of the form (0.1) is isomorphic to the subgroup of
$SU^{\delta}(2,1)$ of matrices
$$
\left(
\begin{array}{ccc}
   \displaystyle E&P&Q+iP\overline{P}\\
   \displaystyle 0&E&2i\overline{P}\\
\displaystyle 0&0&E\\
\end{array}
\right),
$$
with $P\in{\frak A}^{\delta}$, $Q\in\hbox{Re}\,{\frak A}^{\delta}$.
Since any element of $\hbox{Aut}_{e}(Q_{H^{\delta}})$ is the composition
of an automorphism from $\hbox{Aut}_{0,e}(Q_{H^{\delta}})$ and an automorphism of the form
(0.1), it follows that  $\hbox{Aut}_{e}(Q_{H^{\delta}})$ is isomorphic
to $SU^{\delta}(2,1)/Z^{\delta}$. Therefore, ${\frak g}_{H^{\delta}}$
is isomorphic to ${\frak {su}}^{\delta}(2,1)$, the Lie algebra of
$SU^{\delta}(2,1)$, which is clearly the algebra of matrices
$$
\left(
\begin{array}{ccc}
   \displaystyle X&Y&Z\\
   \displaystyle W&-2i\hbox{Im}X&2i\overline{Y}\\
\displaystyle V&-\frac{i}{2}\overline{W}&-\overline{X}\\
\end{array}
\right),\eqno{(1.3)}
$$
where $X,Y,W\in{\frak A}^{\delta}$, $Z,V\in\hbox{Re}\,{\frak
A}^{\delta}$.

Further, the group $\hbox{Aut}_{0,e}(Q_{H^{\delta}})$ turns out to be
isomorphic to the group of matrices of the form
$$
\left(
\begin{array}{cccc}
   \displaystyle C^{-1}\overline{C}^{-1}&0&0&0\\
   \displaystyle T&\overline{C}^{-1}&0&0\\
\displaystyle \overline{T}&0&C^{-1}&0\\
\displaystyle S&iC\overline{T}&-i\overline{C}T&E
\end{array}
\right),\eqno{(1.4)}
$$
where $T\in{\frak A}^{\delta}$, $C\in{\frak
A}^{\delta *}$, $S\in\hbox{Re}\,{\frak A}^{\delta}$. The isomorphism
that we denote by $\chi^{\delta}$ is
given explicitly as follows: let an element
$\eta\in\hbox{Aut}_{0,e}(Q_{H^{\delta}})$ be represented by matrix
(1.2);
then we set
$$
\chi^{\delta}(\eta)=\left(
\begin{array}{cccc}
   \displaystyle C^{-1}\overline{C}^{-1}&0&0&0\\
   \displaystyle -2AC^{-1}\overline{C}^{-1}&C^{-1}&0&0\\
\displaystyle -2\overline{A}C^{-1}\overline{C}^{-1}&0&\overline{C}^{-1}&0\\
\displaystyle -4RC^{-1}\overline{C}^{-1}&-2i\overline{A}C^{-1}&2iA\overline{C}^{-1}&E
\end{array}
\right).\eqno{(1.5)}
$$
It is straightforward to check that $\chi^{\delta}$ is a group
isomorphism.

The primary description of $\hbox{Aut}_{0,e}(Q_{H^{\delta}})$ in
\cite{ES1} was in fact given in terms of rational mappings of
$\CC^4$. In particular, it was shown that all automorphisms
from $\hbox{Aut}_{lin,e}(Q_{H^{\delta}})$ have the form
\begin{eqnarray*}
z^{*}&=&Cz,\\
w^{*}&=&C\overline{C}w,\qquad\qquad\qquad\qquad (1.6)\\
\end{eqnarray*}
with $C\in{\frak A}^{\delta *}$. Any element of
$\hbox{Aut}_{0,e}(Q_{H^{\delta}})$ is then a composition of a
rational mapping $z^{*}=z^{*}(z,w), w^{*}=w^{*}(z,w)$, such that
$\partial z^{*}/\partial z(0)=E$, $\partial w^{*}/\partial z(0)=0$,
$\partial w^{*}/\partial w(0)=E$, and an automorphism of the form
(1.6). It is also easy to see that the full group
$\hbox{Aut}_0(Q_{H^{\delta}})$ has exactly two connected
components, and that the second component is obtained by taking the
compositions of mappings from
$\hbox{Aut}_{0,e}(Q_{H^{\delta}})$ and the linear automorphism
$$
z^{*}=
\left(
\begin{array}{ccc}
\displaystyle 1&0\\
\displaystyle 0&-1\\
\end{array}
\right)z,\qquad
w^{*}=
\left(\begin{array}{ccc}
\displaystyle 1&0\\
\displaystyle 0&-1\\
\end{array}
\right)w.
$$
Thus, automorphisms from $\hbox{Aut}_{0,e}(Q_{H^{\delta}})$
are characterized among all elements of
$\hbox{Aut}_0(Q_{H^{\delta}})$ by the condition
$\hbox{det}\left(\partial w^{*}/\partial w(0)\right)>0$.

We are now ready to formulate the main result of this paper.
Let $G^{1,\delta}$ be the group of elements of $\hbox{Re}\,{\frak
A}^{\delta *}$ of the form $C\overline{C}$, $C\in{\frak A}^{\delta
*}$, and $G^{2,\delta}$ be the subgroup of
$\hbox{Aut}_{0,e}(Q_{H^{\delta}})$ defined by the
condition $C\overline{C}=E$.

\begin{theorem} \sl Let $M\in{\cal C}^{\delta}$ be an oriented
manifold. Then there are principal bundles $P^{1,\delta}$, $P^{2,\delta}$
$$
P^{2,\delta}\stackrel{\pi^{2,\delta}}{\rightarrow} P^{1,\delta}\stackrel{\pi^{1,\delta}}{\rightarrow}M
$$
with structure groups $G^{1,\delta}$, $G^{2,\delta}$ respectively and a 1-form $\omega^{\delta}$
on $P^{2,\delta}$ such that at any point $x\in P^{2,\delta}$, $\omega^{\delta}(x)$ is an
isomorphism between $T_x(P^{2,\delta})$ and ${\frak
{su}}^{\delta}(2,1)$, and the following holds:
\smallskip\\

\noindent {\bf (i)} Any $CR$-isomorphism $f:M_1\ra M_2$
between oriented manifolds $M_1,M_2\in{\cal C}^{\delta}$ that
preserves orientation, can be lifted to a diffeomorphism $F: P^{2,\delta}_1\ra
P^{2,\delta}_2$ such that $F^{*}\omega_2^{\delta}=\omega_1^{\delta}$;
\smallskip\\

\noindent {\bf (ii)} Any diffeomorphism $F: P^{2,\delta}_1\ra P^{2,\delta}_2$ such that
$F^{*}\omega_2^{\delta}=\omega_1^{\delta}$, is a lift of a $CR$-isomorphism
$f:M_1\ra M_2$ that preserves orientation, for $M_1,M_2\in{\cal C}^{\delta}$.
\smallskip\\

Moreover, there exists an explicit transformation formula for
$\omega^{\delta}$ under the action on $G^{2,\delta}$ on $P^{2,\delta}$: if for $\eta\in\hbox{Aut}_{0,e}(Q_H)$, $L_{\eta}$
denotes the left action of $G^{2,\delta}$ on $P^{2,\delta}$ by 
$\eta$, then $L_{\eta}^{*}\omega^{\delta}=\hbox{Ad}(\eta)\omega^{\delta}+\dots$, where $\hbox{Ad}$ is the
adjoint representation of $\hbox{Aut}_{e}(Q_{H^{\delta}})$ on ${\frak
{su}}^{\delta}(2,1)$, and $\dots$ denotes  an error term (see formula
(2.59) below).
\end{theorem}

\begin{remark} \rm The condition for the manifolds to be oriented is not
important. Theorem 1.1 could be formulated for any manifold from ${\cal
C}^{\delta}$, but then the group $G^{1,\delta}$ would have to be replaced by
$$
\tilde{G}^{1,\delta}:=\left\{C\overline{C}: C\in{\frak A}^{\delta *}\right\}
\bigcup
\left\{\left(
\begin{array}{ccc}
\displaystyle 1&0\\
\displaystyle 0&-1\\
\end{array}
\right)C\overline{C}:C\in{\frak A}^{\delta *}\right\}.
$$
The group $\tilde{G}^{1,\delta}$ is disconnected. Thus, the bundle
$P^{1,\delta}$ would have a disconnected fibre, and, for an oriented
$M$, the bundle $P^{1,\delta}$ and therefore the bundle $P^{2,\delta}$
would be
disconnected. To avoid these kinds of disconnectedness, we require
the manifolds to be oriented.
\end{remark}

\begin{remark} \rm Everywhere in this paper we assume manifolds to be
$C^{\infty}$-smooth. However, an inspection of the proof of Theorem
1.1 below shows that it is enough to require only some
sufficiently high, but finite, degree of smoothness.
\end{remark}

\section{Proof of Theorem 1.1}

Let $M$ be an oriented connected manifold from ${\cal
C}^{\delta}$. We now fix $\delta$ and drop it in all superscripts.
For any $p\in M$ denote by ${\frak M}_p$ the set of all pairs
$(\theta^1,\theta^2)$ of real linearly independent 1-forms defined in a
neighbourhood of $p$ such that:
\smallskip\\

\noindent {\bf (i)} $T_q^c(M)=\left\{X\in T_q(M):
\theta^1(q)(X)=\theta^2(q)(X)=0\right\}$ for $q$ close to $p$,
\smallskip\\

\noindent {\bf (ii)} There exist complex 1-forms $\omega^1, \omega^2$
near $p$ such that: for all $q$ close to $p$ they are complex linear on $T^c_q(M)$;
$\left(\theta^{\alpha}(q),\hbox{Re}\,\omega^{\alpha}(q),
\hbox{Im}\,\omega^{\alpha}(q)\right)$ is a coframe, and near $p$ the
following holds
\begin{eqnarray*}
d\theta^1&=&i\left(\omega^1\wedge\overline{\omega^1}+
\delta\omega^2\wedge\overline{\omega^2}\right)\qquad\left(\hbox{mod}\,\theta^{\alpha}\right),\\
d\theta^2&=&i\left(\omega^1\wedge\overline{\omega^2}+\omega^2\wedge\overline{\omega^1}\right)
\qquad\left(\hbox{mod}\,\theta^{\alpha}\right).\qquad\qquad\qquad\qquad
(2.1)
\end{eqnarray*}

The integrability of the $CR$-structure and the type of the Levi form
imply that ${\frak M}_p\ne\emptyset$ for any $p\in M$.

We define a smooth bundle $P^1\ra M$ as
$$
P^1:=\left\{{\frak M}_p^{+}\right\}_{p\in M},
$$
where ${\frak M}_p^{+}$ is the set of pairs
$y:=(\theta^1(p),\theta^2(p))$ with $(\theta^1,\theta^2)\in{\frak M}_p$ such
that the orientation that
they define on the cotangent space $T^{*}_p(M)$ coincides with that induced on
$T^{*}_p(M)$ by the original orientation of $M$. Clearly, $P^1$ so
defined is a principal $G^1$-bundle over $M$. We introduce fibre
coordinates $(a,b)$ on $P^1$ via the entries of $C\overline{C}$:
$$
\left(
\begin{array}{ccc}
\displaystyle a&\delta b\\
\displaystyle b&a\\
\end{array}
\right)=C\overline{C}.
$$

To construct the bundle $P^2\ra P^1$ we need the following technical
lemma.

\begin{lemma} \sl Let $(\theta^1,\theta^2)\in {\frak M}_p$ be such
that $(\theta^1(p),\theta^2(p))\in{\frak M}_p^{+}$. Then
$\omega^1,\omega^2$ in (2.1) can be chosen so that the following
holds:
\begin{eqnarray*}
d\theta^1&=&i\left(\omega^1\wedge\overline{\omega^1}+
\delta\omega^2\wedge\overline{\omega^2}\right)+\theta^1\wedge\phi^1+\delta\theta^2\wedge\phi^2,\\
d\theta^2&=&i\left(\omega^1\wedge\overline{\omega^2}+\omega^2\wedge\overline{\omega^1}\right)+\theta^1\wedge\phi^2+\theta^2\wedge\phi^1+\\
&{}&2\theta^1\wedge\hbox{Re}\,\left(\delta
r_1\omega^1+r_2\omega^2\right)+
2\theta^2\wedge\hbox{Re}\,\left(
r_2\omega^1+r_1\omega^2\right),\qquad\qquad\qquad\qquad
(2.2)
\end{eqnarray*}
where $\phi^1,\phi^2$ are real 1-forms and $r_1,r_2$ are smooth
complex-valued functions near $p$.
\end{lemma}

{\bf Proof.} By Proposition 3.2 of \cite{M} we can assume that $\omega^1,\omega^2$ are
chosen in such a way that
\begin{eqnarray*}
d\theta^1&=&i\left(\omega^1\wedge\overline{\omega^1}+
\delta\omega^2\wedge\overline{\omega^2}\right)+\theta^1\wedge\phi^{1'},\\
d\theta^2&=&i\left(\omega^1\wedge\overline{\omega^2}+\omega^2\wedge\overline{\omega^1}\right)+\theta^2\wedge\phi^{2'},\qquad\qquad\qquad\qquad
(2.3)
\end{eqnarray*}
where $\phi^{\alpha'}$ are real 1-forms near
$p$. Since
$(\theta^{\alpha},\hbox{Re}\,\omega^{\alpha},\hbox{Im}\,\omega^{\alpha})$
gives a coframe at every point near $p$, we have
$$
\phi^{\alpha'}=a^{\alpha}_{\gamma}\omega^{\gamma}+\overline{a^{\alpha}_{\gamma}\omega^{\gamma}}+
b^{\alpha}_{\gamma}\theta^{\gamma},\qquad \alpha=1,2,
$$
where $\a^{\alpha}_{\gamma}$ are complex-valued and
$b^{\alpha}_{\gamma}$ are real-valued functions near $p$. 

We choose the new forms $\omega^{\alpha *}$ as follows:
\begin{eqnarray*}
\omega^{1 *}&:=&\omega^1,\\
\omega^{2 *}&:=&\omega^2+\frac{i\delta}{2}(\overline{a_2^2}-\overline{a_2^1})\theta^1+
\frac{i}{2}(\overline{a_1^1}-\overline{a_1^2})\theta^2.\\
\end{eqnarray*}
It is now straightforward to check that under this transformation
equations (2.3) take the form (2.2).

The lemma is proved.\hfill $\Box$

Let $\tilde{\theta^1},\tilde{\theta^2}$ be the following globally
defined tautological 1-forms on $P^1$. For $y:=(\theta^1(p),\theta^2(p))\in P^1$ set
$$
\tilde{\theta^{\alpha}}(y)=(\pi^{1*}\theta^{\alpha})(y),\qquad \alpha=1,2,
$$
where $\pi^1:P^1\ra M$ is the natural projection: $\pi^1(y)=p$.
We now define the bundle $P^2$ over $P^1$ as follows: the fibre over
$y\in P^1$ is the collection of all coframes at $y$ of the form
$(\tilde{\theta^{\alpha}}(y),\hbox{Re}\,\tilde{\omega^{\alpha}},
\hbox{Im}\,\tilde{\omega^{\alpha}},\tilde{\phi^{\alpha}})$
such that:
\smallskip

\noindent {\bf (i)} $\tilde{\omega^{\alpha}}=\pi^{1
*}\omega^{\alpha}(y)$, for some complex covectors $\omega^{\alpha}$ at
$p$ that are complex-linear on $T_p^c(M)$;
\smallskip

\noindent {\bf (ii)} $\tilde{\phi^{\alpha}}$ are real covectors at
$y$;
\smallskip\\

\noindent {\bf (iii)} For some $\tilde{r_\alpha}\in\CC$
the following holds:
\begin{eqnarray*}
d{\tilde\theta^1}(y)&=&i\left(\tilde{\omega^1}\wedge\overline{\tilde{\omega^1}}+
\delta\tilde{\omega^2}\wedge\overline{\tilde{\omega^2}}\right)+\tilde{\theta^1}(y)\wedge\tilde{\phi^1}+\delta\tilde{\theta^2}(y)\wedge\tilde{\phi^2},\\
d\tilde{\theta^2}(y)&=&i\left(\tilde{\omega^1}\wedge\overline{\tilde{\omega^2}}+\tilde{\omega^2}\wedge\overline{\tilde{\omega^1}}\right)+\tilde{\theta^1}(y)\wedge\tilde{\phi^2}+\tilde{\theta^2}(y)\wedge\tilde{\phi^1}+\\
&{}&2\tilde{\theta^1}(y)\wedge\hbox{Re}\,\left(\delta
\tilde{r_1}\tilde{\omega^1}+\tilde{r_2}\tilde{\omega^2}\right)+
2\tilde{\theta^2}(y)\wedge\hbox{Re}\,\left(
\tilde{r_2}\tilde{\omega^1}+\tilde{r_1}\tilde{\omega^2}\right).
\end{eqnarray*}
The existence of such coframes follows from Lemma 2.1.

From now on we will write an element of $x\in P^2$ in the form:
$x:=(\tilde{\theta^{\alpha}}(y),
\tilde{\omega^{\alpha}},\overline{\tilde{\omega^{\alpha}}},
\tilde{\phi^{\alpha}})$. It is
a routine calculation to verify that the most general linear
transformation that, being applied to $x$, gives an element from
the fibre of $P^2$ over $y$, is of the form (1.4):
$$
\left(
\begin{array}{cccc}
   \displaystyle E&0&0&0\\
   \displaystyle T&C&0&0\\
\displaystyle \overline{T}&0&\overline{C}&0\\
\displaystyle S&iC\overline{T}&-i\overline{C}T&E
\end{array}
\right),\eqno{(2.4)}
$$
where $T\in{\frak A}^{\delta}$, $C\in{\frak
A}^{\delta *}$, $C\overline{C}=E$, $S\in\hbox{Re}\,{\frak
A}^{\delta}$. The group of matrices (2.4) is isomorphic to $G^2$ 
by the isomorphism $\chi$ defined in (1.5). Therefore, $P^2$ is a
principle bundle with structure group $G^2$. We will treat the independent
entries of $T,C,S$ as fibre coordinates. 

We now introduce globally defined tautological 1-forms on $P^2$. Let 
$x=(\tilde{\theta^{\alpha}}(y),\tilde{\omega^{\alpha}},
\overline{\tilde{\omega^{\alpha}}},\tilde{\phi^{\alpha}})\in
P^2$. Then we set:
\begin{eqnarray*}
\hat{\theta^{\alpha}}(x)&:=&(\pi^{2*}\tilde{\theta^{\alpha}})(x),\\
\hat{\omega^{\alpha}}(x)&:=&(\pi^{2*}\tilde{\omega^{\alpha}})(x),\\
\hat{\phi^{\alpha}}(x)&:=&(\pi^{2*}\tilde{\phi^{\alpha}})(x),\qquad\qquad\qquad\qquad\qquad
(2.5)
\end{eqnarray*}
for $\alpha=1,2$, where $\pi^2:P^2\rightarrow P^1$ is the projection: $\pi^2(x)=y$.
To simplify notation, until the end of this section we drop hats in the
forms defined in (2.5). These forms satisfy the equations
\begin{eqnarray*}
d\theta^1&=&i\left(\omega^1\wedge\overline{\omega^1}+
\delta\omega^2\wedge\overline{\omega^2}\right)+\theta^1\wedge\phi^1+
\delta\theta^2\wedge\phi^2,\\
d\theta^2&=&i\left(\omega^1\wedge\overline{\omega^2}+
\omega^2\wedge\overline{\omega^1}\right)+\theta^1
\wedge\phi^2+\theta^2\wedge\phi^1+\\
&{}&2\theta^1\wedge\hbox{Re}\,\left(\delta
r_1\omega^1+r_2\omega^2\right)+
2\theta^2\wedge\hbox{Re}\,\left(
r_2\omega^1+r_1\omega^2\right),\qquad\qquad
(2.6)
\end{eqnarray*}
for some uniquely determined smooth complex-valued functions $r_{\alpha}$ on
$P^2$.

Next, it follows from the integrability of the $CR$-structure that
$$
d\omega^{\alpha}=\omega^{\beta}\wedge\phi^{\alpha}_{\beta}+\theta^{\beta}\wedge
\mu^{\alpha}_{\beta},\eqno{(2.7)}
$$
for some locally defined 1-forms
$\phi^{\alpha}_{\beta}$, $\mu^{\alpha}_{\beta}$.
Differentiating (2.6) and plugging (2.6), (2.7) in the resulting expressions, we get
\begin{eqnarray*}
&{}&i\omega^1\wedge\overline{\omega^1}\wedge\left(\phi^1-2\hbox{Re}\,\phi_1^1\right)+i\omega^1\wedge\overline{\omega^2}\wedge\left(
\delta\phi^2-\delta\phi_1^2-\overline{\phi_2^1}\right)+\\
&{}&i\omega^2\wedge\overline{\omega^1}\wedge\left(\delta\phi^2-\phi_2^1-\delta\overline{\phi_1^2}\right)+i\delta
\omega^2\wedge\overline{\omega^2}\left(
\phi^1-2\hbox{Re}\,\phi_2^2\right)+\\
&{}&\theta^1\wedge\Biggl(-d\phi^1+2\hbox{Re}\,\left(i\mu_1^1\wedge\overline{\omega^1}+i\delta
\mu_1^2\wedge\overline{\omega^2}+
(r_1\omega^1+\delta
r_2\omega^2)\wedge\phi^2\right)\Biggr)+\\
&{}&\theta^2\wedge\Biggl(-\delta
d\phi^2+2\hbox{Re}\,\left(i\mu_2^1\wedge\overline{\omega^1}+i\delta
\mu_2^2\wedge\overline{\omega^2}+
\delta(r_2\omega^1+r_1\omega^2)\wedge\phi^2\right)\Biggr)=0,\qquad\qquad\qquad
(2.8.a)\\
&{}&i\omega^1\wedge\overline{\omega^1}\wedge\left(\phi^2-2\hbox{Re}\,\phi_1^2\right)+i\omega^1\wedge\overline{\omega^2}\wedge\left(
\phi^1-\phi_1^1-\overline{\phi_2^2}\right)+\\
&{}&i\omega^2\wedge\overline{\omega^1}\wedge\left(\phi^1-\phi_2^2-\overline{\phi_1^1}\right)+
i\omega^2\wedge\overline{\omega^2}\left(
\delta\phi^2-2\hbox{Re}\,\phi_2^1\right)+\\
&{}&\theta^1\wedge\Biggl(-d\phi^2+2\hbox{Re}\,\biggl(i\mu_1^2\wedge\overline{\omega^1}+
i\mu_1^1\wedge\overline{\omega^2}-
(r_2\omega^1+
r_1\omega^2)\wedge\phi^2-\\
&{}&\delta
dr_1\wedge\omega^1-dr_2\wedge\omega^2-
\delta
r_1(\omega^1\wedge\phi_1^1+\omega^2\wedge\phi_2^1+\theta^2\wedge\mu_2^1)-\\
&{}&r_2(\omega^1\wedge\phi_1^2+\omega^2\wedge\phi_2^2+\theta^2\wedge\mu_2^2)\biggr)+4\hbox{Re}\,(\delta
r_1\omega^1+r_2\omega^2)\wedge\hbox{Re}\,(r_2\omega^1+r_1\omega^2)
\Biggr)+\\
&{}&\theta^2\wedge\Biggl(
-d\phi^1+2\hbox{Re}\,\biggl(i\mu_2^2\wedge\overline{\omega^1}+
i\mu_2^1\wedge\overline{\omega^2}-
(r_1\omega^1+\delta
r_2\omega^2)\wedge\phi^2-\\
&{}&dr_2\wedge\omega^1-dr_1\omega^2-
r_2(\omega^1\wedge\phi_1^1+\omega^2\wedge\phi_2^1+\theta^1\wedge\mu_1^1)-\\
&{}&r_1(\omega^1\wedge\phi_1^2+\omega^2\wedge\phi_2^2+\theta^1\wedge\mu_1^2)\biggr)
\Biggr)=0.
\qquad\qquad\qquad (2.8.b)
\end{eqnarray*}

\begin{lemma} \sl There exist $\phi^{\alpha}_{\beta}$ satisfying (2.7)
such that the following holds
\begin{eqnarray*}
\hbox{Re}\,\phi_1^1&=&\frac{1}{2}\phi^1,\\
\hbox{Re}\,\phi_1^2&=&\frac{1}{2}\phi^2,\\
\hbox{Re}\,\phi^2_2&=&\frac{1}{2}\phi^1-\hbox{Re}\,\left(\rho_1\omega^1+\rho_2\omega^2\right),\\
\hbox{Re}\,\phi_2^1&=&\frac{1}{2}\delta\phi^2-\hbox{Re}\,\left(\rho_2\omega^1+\delta\rho_1\omega^2\right),\\
\phi_1^1+\overline{\phi_2^2}&=&\phi^1-\rho_1\omega^1-\rho_2\omega^2,\\
\delta\phi_1^2+\overline{\phi_2^1}&=&\delta\phi^2-\rho_2\omega^1-\delta\rho_1\omega^2,\qquad\qquad\qquad\qquad
(2.9)
\end{eqnarray*}
where $\rho_{\alpha}$ are locally defined smooth complex-valued functions on $P^2$.
\end{lemma}

{\bf Proof.} It follows from (2.8) that
\begin{eqnarray*}
\phi^1-2\hbox{Re}\,\phi_1^1&=&2\hbox{Re}\,\left(a_1^1\omega^1+b_1^1\omega^2\right)+c_1^1\theta^1+d_1^1\theta^2,\\
\phi^2-2\hbox{Re}\,\phi_1^2&=&2\hbox{Re}\,\left(a_1^2\omega^1+b_1^2\omega^2\right)+c_1^2\theta^1+d_1^2\theta^2,\\
\phi^1-2\hbox{Re}\,\phi_2^2&=&2\hbox{Re}\,\left(a_2^2\omega^1+b_2^2\omega^2\right)+c_2^2\theta^1+d_2^2\theta^2,\\
\delta\phi^2-2\hbox{Re}\,\phi_2^1&=&2\hbox{Re}\,\left(a_2^1\omega^1+b_2^1\omega^2\right)+c_2^1\theta^1+d_2^1\theta^2,\\
\phi^1-\phi_2^2-\overline{\phi_1^1}&=&a_1^{1'}\omega^1+b_1^{1'}\omega^2+
a_1^{{1'}'}\overline{\omega^1}+b_1^{{1'}'}\overline{\omega^2}+c_1^{1'}\theta^1+d_1^{1'}\theta^2,\\
\delta\phi^2-\delta\phi_1^2-\overline{\phi_2^1}&=&a_2^{1'}\omega^1+b_2^{1'}\omega^2+a_2^{{1'}'}\overline{\omega^1}+b_2^{{1'}'}\overline{\omega^2}+c_2^{1'}\theta^1+d_2^{1'}\theta^2,\qquad
(2.10)
\end{eqnarray*}
where $a^{\alpha}_{\beta},b^{\alpha}_{\beta},a^{\alpha
'}_{\beta},b^{\alpha '}_{\beta},a^{{\alpha'}'}_{\beta}, b^{{\alpha'}'
}_{\beta}, c^{\alpha '}_{\beta},d^{\alpha
'}_{\beta}$ are complex-valued and
$c^{\alpha}_{\beta},d^{\alpha}_{\beta}$ are real-valued functions satisfying
the following relations
\begin{eqnarray*}
a_1^1+a_2^2&=&a_1^{1 '}+\overline{a_1^{{1'}'}},\\
b_1^1+b_2^2&=&b_1^{1 '}+\overline{b_1^{{1'}'}},\\
a_2^1+\delta a_1^2&=&a_2^{1'}+\overline{a_2^{{1'}'}},\\
b_2^1+\delta b_1^2&=&b_2^{1'}+\overline{b_2^{{1'}'}},\\
b_1^1&=&\overline{a_2^{{1'}'}},\\
a_2^2&=&\delta b_2^{1'},\\
a_2^1&=&\overline{b_1^{{1'}'}},\\
b_1^2&=&a_1^{1'},\\
c_1^1+c_2^2&=&2\hbox{Re}\,c_1^{1'},\\
d_1^1+d_2^2&=&2\hbox{Re}\,d_1^{1'},\\
c_2^1+\delta c_1^2&=&2\hbox{Re}\,c_2^{1'},\\
d_2^1+\delta d_1^2&=&2\hbox{Re}\,d_2^{1'}.\qquad\qquad\qquad\qquad (2.11)
\end{eqnarray*}

We set
\begin{eqnarray*}
\phi_1^{1*}&:=&\phi_1^1+a_1^1\omega^1+b_1^1\omega^2+\left(\frac{1}{2}c_1^1+i\hbox{Im}\,c_1^{1'}\right)\theta^1+\left(\frac{1}{2}d_1^1+i\hbox{Im}\,d_1^{1'}\right)\theta^2 ,\\
\phi_1^{2*}&:=&\phi_1^2+a_1^2\omega^1+b_1^2\omega^2+\frac{1}{2}c_1^2\theta^1+\frac{1}{2}d_1^2\theta^2,\\
\phi_2^{1*}&:=&\phi_2^1+b_1^1\omega^1+(b_2^1+\delta b_1^2-\delta
a_2^2)\omega^2+\left(\frac{1}{2}c_2^1+i\hbox{Im}\,c_2^{1'}\right)\theta^1+\left(\frac{1}{2}d_2^1+i\hbox{Im}\,d_2^{1'}\right)\theta^2 ,\\
\phi_2^{2*}&:=&\phi_2^2+b_1^2\omega^1+(b_1^1+b_2^2-a_2^1)\omega^2+\frac{1}{2}c_2^2\theta^1+\frac{1}{2}d_2^2\theta^2  .
\end{eqnarray*}
It now follows from (2.10), (2.11) that these forms satisfy (2.7) and (2.9).

The lemma is proved.\hfill $\Box$

From now on we will assume that $\phi_{\alpha}^{\beta}$ in (2.7)
satisfy conditions (2.9). Then (2.8.a) implies:
\begin{eqnarray*}
d\phi^1&=&2\hbox{Re}\,\left(i\mu_1^1\wedge\overline{\omega^1}+i\delta\mu_1^2\wedge\overline{\omega^2}+\left(r_1\omega^1+\delta
r_2\omega^2\right)\wedge\phi^2\right)+\\
&{}&\theta^1\wedge\psi^1+\theta^2\wedge\psi^2,\qquad\qquad\qquad\qquad
(2.12.a)\\
d\phi^2&=&2\hbox{Re}\,\left(i\delta\mu_2^1\wedge\overline{\omega^1}+i\mu_2^2\wedge\overline{\omega^2}+\left(r_2\omega^1+
r_1\omega^2\right)\wedge\phi^2\right)+\\
&{}&\theta^1\wedge\psi^3+\theta^2\wedge\psi^4.\qquad\qquad\qquad\qquad
(2.12.b)
\end{eqnarray*}
Analogously, (2.8.b) implies:
\begin{eqnarray*}
d\phi^1&=&2\hbox{Re}\,\Biggl(i\mu_2^2\wedge\overline{\omega^1}+i\mu_2^1\wedge\overline{\omega^2}-\left(r_1\omega^1+\delta
r_2\omega^2\right)\wedge\phi^2-\\
&{}&dr_2\wedge\omega^1-dr_1\wedge\omega^2-
r_2\left(\omega^1\wedge\phi_1^1+\omega^2\wedge\phi_2^1+\theta^1\wedge\mu_1^1\right)-
r_1\biggl(\omega^1\wedge\phi_1^2+\\
&{}&\omega^2\wedge\phi_2^2+
\theta^1\wedge\mu_1^2\biggr)\Biggr)+
\theta^1\wedge\psi^5+\theta^2\wedge\psi^6,\qquad\qquad
(2.12.c)\\
d\phi^2&=&2\hbox{Re}\,\Biggl(i\mu_1^2\wedge\overline{\omega^1}+i\mu_1^1\wedge\overline{\omega^2}-
\left(r_2\omega^1+
r_1\omega^2\right)\wedge\phi^2-\\
&{}&\delta dr_1\wedge\omega^1-dr_2\wedge\omega^2-
\delta r_1\left(\omega^1\wedge\phi_1^1+\omega^2\wedge\phi_2^1+\theta^2\mu_2^1\right)-\\
&{}&r_2\left(\omega^1\wedge\phi_1^2+\omega^2\wedge\phi_2^2+\theta^2\wedge\mu_2^2\right)\Biggr)+4\hbox{Re}\,\left(\delta
r_1\omega^1+r_2\omega^2\right)\wedge\\
&{}&\left(r_2\omega^1+r_1\omega^2\right)+
\theta^1\wedge\psi^7+\theta^2\wedge\psi^8.\qquad\qquad\qquad\qquad
(2.12.d)
\end{eqnarray*}
In formulas (2.12) $\psi^{\alpha}$ are real locally defined 1-forms such that
\begin{eqnarray*}
\psi^2&=&\delta\psi^3+s_1\theta^1+s_2\theta^2,\\
\psi^8&=&\psi^5+s_3\theta^1+s_4\theta^2,\qquad\qquad\qquad\qquad (2.13)
\end{eqnarray*}
where $s_{\alpha}$ are real-valued functions.

A lengthy but elementary calculation now shows that the 1-forms
$\phi^{\alpha}_{\beta}$, $\mu^{\alpha}_{\beta}$, $\psi^{\alpha}$
satisfying (2.7), (2.9), (2.12) are defined up to transformations of
the form
\begin{eqnarray*}
\phi_1^{1*}&=&\phi_1^1+h\theta^1+g\theta^2,\\
\phi_1^{2*}&=&\phi_1^2+\delta h'\theta^1+\delta g'\theta^2,\\
\phi_2^{1*}&=&\phi_2^1+h'\theta^1+g'\theta^2,\\
\phi_2^{2*}&=&\phi_2^2+h\theta^1+g\theta^2,\\
\mu_1^{1*}&=&\mu_1^1+h\omega^1+h'\omega^2+p_1\theta^1+q_1\theta^2,\\
\mu_1^{2*}&=&\mu_1^2+\delta
h'\omega^1+h\omega^2+p_2\theta^1+q_2\theta^2,\\
\mu_2^{1*}&=&\mu_2^1+g\omega^1+g'\omega^2+q_1\theta^1+q_3\theta^2,\\
\mu_2^{2*}&=&\mu_2^2+\delta
g'\omega^1+g\omega^2+q_2\theta^1+q_4\theta^2,\\
\psi^{1*}&=&\psi^1+2\hbox{Re}\,\left(i\overline{p_1}\omega^1+i\delta\overline{p_2}\omega^2\right)+\sigma_1\theta^1+\sigma_2\theta^2,\\
\psi^{2*}&=&\psi^2+2\hbox{Re}\,\left(i\overline{q_1}\omega^1+i\delta\overline{q_2}\omega^2\right)+\sigma_2\theta^1+\sigma_3\theta^2,\\
\psi^{3*}&=&\psi^3+2\hbox{Re}\,\left(i\delta\overline{q_1}\omega^1+i\overline{q_2}\omega^2\right)+\sigma_4\theta^1+\sigma_5\theta^2,\\
\psi^{4*}&=&\psi^4+2\hbox{Re}\,\left(i\delta\overline{q_3}\omega^1+i\overline{q_4}\omega^2\right)+\sigma_5\theta^1+\sigma_6\theta^2,\\
\psi^{5*}&=&\psi^5+2\hbox{Re}\,\left(i\overline{q_2}\omega^1+i\overline{q_1}\omega^2\right)+\sigma_7\theta^1+\sigma_8\theta^2,\\
\psi^{6*}&=&\psi^6+2\hbox{Re}\,\left(i\overline{q_4}\omega^1+i\overline{q_3}\omega^2\right)+\\
&{}&\bigl(\sigma_8-2\hbox{Re}\,(r_1q_2+r_2q_1)\bigr)\theta^1+\sigma_9\theta^2,\\
\psi^{7*}&=&\psi^7+2\hbox{Re}\,\left(i\overline{p_2}\omega^1+i\overline{p_1}\omega^2\right)+\sigma_{10}\theta^1+\sigma_{11}\theta^2,\\
\psi^{8*}&=&\psi^8+2\hbox{Re}\,\left(i\overline{q_2}\omega^1+i\overline{q_1}\omega^2\right)+\\
&{}&\bigl(\sigma_{11}+2\hbox{Re}\,(\delta
r_1q_1+r_2q_2)\bigr)\theta^1+\sigma_{12}\theta^2,\qquad\qquad\qquad\qquad
(2.14)
\end{eqnarray*}
where $g,g',h,h'$ are imaginary-valued, $\sigma_{\alpha}$
are real-valued, $p_{\alpha}$,
$q_{\alpha}$ are complex-valued functions. The same calculation shows that
$\rho_{\alpha}$ are chosen uniquely and therefore are globally
defined on $P^2$.

We now need to introduce extra conditions that  would fix the parameters in
(2.14) uniquely and therefore, taken together with (2.7), (2.9) and (2.12), would fix the forms $\phi^{\alpha}_{\beta}$,
$\mu^{\alpha}_{\beta}$, $\psi^{\alpha}$ uniquely. The first set of
conditions comes from comparing two pairs of equations: (2.12.a),
(2.12.c) and (2.12.b), (2.12.d). From this comparison we get:
\begin{eqnarray*}
\mu_2^2&=&\mu_1^1+a_1\omega^1+b_1\omega^2+c_1\overline{\omega^1}+d_1\overline{\omega^2}+h_1\theta^1+g_1\theta^2+\\
&{}&2i\overline{r_1}\phi^2-id\overline{r_2}+i\overline{r_2}\overline{\phi_1^1}+i\overline{r_1}\overline{\phi_1^2},\\
\mu_2^1&=&\delta\mu_1^2-\overline{b_1}\omega^1+b_2\omega^2+d_1\overline{\omega^1}+d_2\overline{\omega^2}+h_2\theta^1+g_2\theta^2+\\
&{}&2i\delta\overline{r_2}\phi^2-id\overline{r_1}+i\overline{r_2}\overline{\phi_2^1}+i\overline{r_1}\overline{\phi_2^2},\\
\psi^5&=&\psi^1+2\hbox{Re}\,\left(i\overline{h_1}\omega^1+i\overline{h_2}\omega^2\right)+s_5\theta^1+s_6\theta^2+\\
&{}&2\hbox{Re}\,\left(r_2\mu_1^1+r_1\mu_1^2\right),\\
\psi^6&=&\psi^2+2\hbox{Re}\,\left(i\overline{g_1}\omega^1+i\overline{g_2}\omega^2\right)+s_6\theta^1+s_7\theta^2,\\
\mu_2^2&=&\mu_1^1+
a_3\omega^1+b_3\omega^2+c_3\overline{\omega^1}+d_3\overline{\omega^2}+
h_3\theta^1+g_3\theta^2-\\
&{}&2i\overline{r_1}\phi^2+id\overline{r_2}-i\delta\overline{r_1}\overline{\phi_2^1}-i\overline{r_2}\overline{\phi_2^2},\\
\mu_2^1&=&\delta\mu_1^2+a_4\omega^1+(-\delta\overline{a_3}+i(|r_1|^2-\delta|r_2|^2)\omega^2+c_4\overline{\omega^1}+\\
&{}&(\delta
c_3+i(\overline{r_1^2}+\delta\overline{r_2^2}))\overline{\omega^2}+
\delta
h_4\theta^1+\delta g_4\theta^2-
2i\delta\overline{r_2}\phi^2+id\overline{r_1}-i\overline{r_1}\overline{\phi_1^1}-i\delta\overline{r_2}\overline{\phi_1^2},\\
\psi^7&=&\psi^3-2\hbox{Re}\,\left(i\overline{h_4}\omega^1+i\overline{h_3}\omega^2\right)+s_8\theta^1+s_9\theta^2,\\
\psi^8&=&\psi^4-2\hbox{Re}\,\left(i\overline{g_4}\omega^1+i\overline{g_3}\omega^2\right)+s_9\theta^1+s_{10}\theta^2+\\
&{}&2\hbox{Re}\,\left(\delta r_1\mu_2^1+r_2\mu_2^2\right),\qquad\qquad\qquad\qquad\qquad\qquad\qquad(2.15)
\end{eqnarray*}
where $a_1,b_2$ are imaginary-valued,
$s_{\alpha}$ are real-valued, the
rest of the functions are complex-valued, and
$\hbox{Re}\,b_3=\hbox{Im}\,(r_2\overline{r_1})$,
$\hbox{Re}\,a_4=\delta\hbox{Im}\,(r_1\overline{r_2})$.

We now choose $g,g',h,'h'$ in (2.14) so that
$$
a_1^{*}=0,\qquad \hbox{Im}\,a_4^{*}=0,\eqno{(2.16)}
$$
where the functions with asterisks correspond to the forms with asterisks
from (2.14).
This can be achieved by setting
\begin{eqnarray*}
h-\delta g'&=&a_1,\\
h'-g&=&i\hbox{Im}\,a_4.\qquad\qquad\qquad\qquad (2.17)
\end{eqnarray*}
Choice (2.17) uniquely fixes the functions $a_1^{*}$, $a_4^{*}$ and therefore
all the functions $a_{\alpha}^{*}$, $b_{\alpha}^{*}$, $c_{\alpha}^{*}$,
$d_{\alpha}^{*}$.

We also choose $\sigma_{\alpha}$ in (2.14) so that in (2.13), (2.15)
one has 
$$
s_{\alpha}^{*}=0, \qquad \alpha=1,\dots,10,\eqno{(2.18)}
$$
by setting
\begin{eqnarray*}
\delta\sigma_4-\sigma_2&=&s_1,\\
\delta\sigma_5-\sigma_3&=&s_2,\\
\sigma_7-\sigma_{11}-2\hbox{Re}\,(\delta r_1q_1+r_2q_2)&=&s_3,\\
\sigma_8-\sigma_{12}&=&s_4,\\
\sigma_1-\sigma_7+2\hbox{Re}\,(r_2p_1+r_1p_2)&=&s_5,\\
\sigma_2-\sigma_8+2\hbox{Re}\,(r_2q_1+r_1q_2)&=&s_6,\\
\sigma_3-\sigma_9&=&s_7,\\
\sigma_4-\sigma_{10}&=&s_8,\\
\sigma_5-\sigma_{11}&=&s_9,\\
\sigma_6-\sigma_{12}+2\hbox{Re}\,(\delta
r_1q_3+r_2q_4)&=&s_{10}.\qquad\qquad\qquad\qquad(2.19)
\end{eqnarray*}

To introduce further restrictions on the parameters
$g,g',h,h',p_{\alpha},q_{\alpha},\sigma_{\alpha}$ we need to
differentiate equations (2.7). By doing this and using (2.6), (2.7) we get
\begin{eqnarray*}
&{}&\omega^1\wedge\Biggl(-d\phi_1^1+\phi_1^2\wedge\phi_2^1-i\mu_1^1\wedge\overline{\omega^1}-i\mu_2^1\wedge\overline{\omega^2}\Biggr)+\\
&{}&\omega^2\wedge\Biggl(-d\phi_2^1+\phi_2^1\wedge\phi_1^1+\phi_2^2\wedge\phi_2^1-i\mu_2^1\wedge\overline{\omega^1}-i\delta\mu_1^1\wedge\overline{\omega^2}\Biggr)+\\
&{}&\theta^1\wedge\Biggl(-d\mu_1^1+\mu_1^1\wedge\phi_1^1+\mu_1^2\wedge\phi_2^1-\mu_1^1\wedge\phi^1-\\
&{}&\mu_2^1\wedge\biggl(\phi^2+2\hbox{Re}\,(\delta
r_1\omega^1+r_2\omega^2)\biggr)\Biggr)+\\
&{}&\theta^2\wedge\Biggl(-d\mu_2^1+\mu_2^1\wedge\phi_1^1+\mu_2^2\wedge\phi_2^1-\delta\mu_1^1\wedge\phi^2-\\
&{}&\mu_2^1\wedge\biggl(\phi^1+2\hbox{Re}\,(r_2\omega^1+r_1\omega^2)\biggr)\Biggr)=0,\\
&{}&\omega^1\wedge\Biggl(-d\phi_1^2+\phi_1^1\wedge\phi_1^2+\phi_1^2\wedge\phi_2^2-i\mu_1^2\wedge\overline{\omega^1}-i\mu_2^2\wedge\overline{\omega^2}\Biggr)+\\
&{}&\omega^2\wedge\Biggl(-d\phi_2^2+\phi_2^1\wedge\phi_1^2-i\mu_2^2\wedge\overline{\omega^1}-i\delta\mu_1^2\wedge\overline{\omega^2}\Biggr)+\\
&{}&\theta^1\wedge\Biggl(-d\mu_1^2+\mu_1^1\wedge\phi_1^2+\mu_1^2\wedge\phi_2^2-\mu_1^2\wedge\phi^1-\\
&{}&\mu_2^2\wedge\biggl(\phi^2+2\hbox{Re}\,(\delta
r_1\omega^1+r_2\omega^2)\biggr)\Biggr)+\\
&{}&\theta^2\wedge\Biggl(-d\mu_2^2+\mu_2^1\wedge\phi_1^2+\mu_2^2\wedge\phi_2^2-\delta\mu_1^2\wedge\phi^2-\\
&{}&\mu_2^2\wedge\biggl(\phi^1+2\hbox{Re}\,(r_2\omega^1+r_1\omega^2)\biggr)\Biggr)=0.\qquad\qquad\qquad
(2.20)
\end{eqnarray*}

Let $\nu^1:=\hbox{Im}\,\phi_1^1$, $\nu^2:=\hbox{Im}\,\phi_1^2$. Then
it is easy to see that for $x\in P^2$,\linebreak
$(\theta^{\alpha}(x),\omega^{\alpha}(x),\overline{\omega^{\alpha}}(x),
\phi^{\alpha}(x), \nu^{\alpha}(x),
\mu_1^1(x),\mu_1^2(x),\overline{\mu_1^1}(x),\overline{\mu_1^2}(x), \psi^1(x),\psi^3(x))$ is
a coframe at $x$. From now on we will use the independent 1-forms
 $\theta^{\alpha},\omega^{\alpha},\overline{\omega^{\alpha}},
\phi^{\alpha}, \nu^{\alpha},
\mu_1^1,\mu_1^2,\overline{\mu_1^1},\overline{\mu_1^2},\psi^1,\psi^3$ as the standard
basis in which we will be writing the expansions of all differential forms that
we will need in the future. Equations (2.7), (2.9), (2.15), (2.20) imply
\begin{eqnarray*}
\omega^1\wedge d\phi_1^1+\delta\omega^2\wedge
d\phi_1^2+A&\equiv&0\qquad (\hbox{mod}\,\theta^{\alpha}),\\
\omega^1\wedge d\phi_1^2+\omega^2\wedge
d\phi_1^1+B&\equiv&0\qquad
(\hbox{mod}\,\theta^{\alpha}),\qquad\qquad\qquad (2.21)
\end{eqnarray*}
where
\begin{eqnarray*}
A&:=&\omega^1\wedge\Biggl(\frac{1}{2}\overline{\rho_2}\phi^2\wedge\overline{\omega^1}-\frac{1}{2}\overline{r_1}\phi^1\wedge\overline{\omega^2}+\delta\biggl(\frac{1}{2}\overline{\rho_1}-\frac{5}{2}\overline{r_2}\biggr)\phi^2\wedge\overline{\omega^2}+i\overline{\rho_2}\nu^2\wedge\overline{\omega^1}+\\
&{}&i\overline{r_1}\nu^1\wedge\overline{\omega^2}+i\delta\biggl(\overline{\rho_1}+\overline{r_2}\biggr)\nu^2\wedge
\overline{\omega^2}+i\mu_1^1\wedge\overline{\omega^1}+i\biggl(\delta\mu_1^2+(b_2-i\delta\overline{r_2}\rho_1-i\overline{r_1}\rho_2)\omega^2+\\
&{}&d_1\overline{\omega^1}-id\overline{r_1}\biggr)\wedge\overline{\omega^2}\Biggr)+
\omega^2\wedge\Biggl(\biggl(\frac{1}{2}\overline{\rho_2}-\frac{1}{2}\overline{r_1}\biggr)\phi^1\wedge\overline{\omega^1}-\delta\frac{5}{2}\overline{r_2}\phi^2\wedge\overline{\omega^1}+\delta\frac{1}{2}\overline{\rho_1}\phi^1\wedge\overline{\omega^2}+\\
&{}&i\biggl(\overline{r_1}-\overline{\rho_2}\biggr)\nu^1\wedge\overline{\omega^1}+i\delta\biggl(\overline{r_2}-2\overline{\rho_1}\biggr)\nu^2\wedge\overline{\omega^1}-i\delta\overline{\rho_1}\nu^1\wedge\overline{\omega^2}-2i\delta\overline{\rho_2}\nu^2\wedge\overline{\omega^2}+\\
&{}&i\biggl(\delta\mu_1^2-(\overline{b_1}+i\overline{r_1}\rho_1+i\overline{r_2}\rho_2)\omega^1+(d_2-i\delta\overline{\rho_1^2}+i\overline{\rho_2^2})\overline{\omega^2}-id\overline{r_1}+id\overline{\rho_2}\biggr)\wedge\overline{\omega^1}+\\
&{}&i\biggl(\delta\mu_1^1+i(\delta|\rho_1|^2+|\rho_2|^2)\omega^1+
i\delta d\overline{\rho_1}
\biggr)\wedge\overline{\omega^2}\Biggr),\\
B&:=&\omega^1\wedge\Biggl(\frac{1}{2}\overline{\rho_1}\phi^2\wedge\overline{\omega^1}-\frac{1}{2}\overline{r_2}\phi^1\wedge\overline{\omega^2}+\biggl(\frac{1}{2}\overline{\rho_2}-\frac{5}{2}\overline{r_1}\biggr)\phi^2\wedge\overline{\omega^2}+i\overline{\rho_1}\nu^2\wedge\overline{\omega^1}+\\
&{}&i\overline{r_2}\nu^1\wedge\overline{\omega^2}+i\biggl(\overline{\rho_2}+\overline{r_1}\biggr)\nu^2\wedge
\overline{\omega^2}+i\mu_1^2\wedge\overline{\omega^1}+i\biggl(\mu_1^1+b_1\omega^2+\\
&{}&c_1\overline{\omega^1}-id\overline{r_2}\biggr)\wedge\overline{\omega^2}\Biggr)+
\omega^2\wedge\Biggl(\biggl(\frac{1}{2}\overline{\rho_1}-\frac{1}{2}\overline{r_2}\biggr)\phi^1\wedge\overline{\omega^1}-\frac{5}{2}\overline{r_1}\phi^2\wedge\overline{\omega^1}+\frac{1}{2}\overline{\rho_2}\phi^1\wedge\overline{\omega^2}+\\
&{}&i\biggl(\overline{r_2}-\overline{\rho_1}\biggr)\nu^1\wedge\overline{\omega^1}+i\biggl(\overline{r_1}-2\overline{\rho_2}\biggr)\nu^2\wedge\overline{\omega^1}-i\overline{\rho_2}\nu^1\wedge\overline{\omega^2}-2i\delta\overline{\rho_1}\nu^2\wedge\overline{\omega^2}+\\
&{}&i\biggl(\mu_1^1+a_1\omega^1+d_1\overline{\omega^2}-id\overline{r_2}+id\overline{\rho_1}\biggr)\wedge\overline{\omega^1}+\\
&{}&i\biggl(\delta\mu_1^2+2i\hbox{Re}\,(\rho_1\overline{\rho_2})\omega^1+id\overline{\rho_2}
\biggr)\wedge\overline{\omega^2}\Biggr).\qquad\qquad
(2.22)
\end{eqnarray*}

It follows from (2.21) that
\begin{eqnarray*}
\omega^1\wedge A&\equiv&-\delta\omega^2\wedge B \qquad (\hbox{mod}\,\theta^{\alpha}),\\
\omega^2\wedge A&\equiv&-\omega^1\wedge B \qquad (\hbox{mod}\,\theta^{\alpha}),
\end{eqnarray*}
which together with (2.22) gives the following expressions for
$dr_{\alpha}$
\begin{eqnarray*}
dr_1&=&\frac{1}{2}r_1\phi^1+\delta\frac{5}{2}r_2\phi^2+ir_1\nu^1+i\delta
r_2\nu^2+\\
&{}&t_1\omega^1+t_2\omega^2+t_3\overline{\omega^1}+t_4\overline{\omega^2}+v_1\theta^1+v_2\theta^2,\\
dr_2&=&\frac{1}{2}r_2\phi^1+\frac{5}{2}r_1\phi^2+ir_2\nu^1+i
r_1\nu^2+\\
&{}&t_5\omega^1+t_6\omega^2+t_7\overline{\omega^1}+t_8\overline{\omega^2}+v_3\theta^1+v_4\theta^2,\qquad
(2.23)
\end{eqnarray*}
and for $d\rho_{\alpha}$
\begin{eqnarray*}
d\rho_1&=&\frac{1}{2}\rho_1\phi^1-\frac{1}{2}\rho_2\phi^2+i\rho_1\nu^1+3i
\rho_2\nu^2+\\
&{}&u_1\omega^1+u_2\omega^2+u_3\overline{\omega^1}+u_4\overline{\omega^2}+w_1\theta^1+w_2\theta^2,\\
d\rho_2&=&\frac{1}{2}\rho_2\phi^1-\delta\frac{1}{2}r_1\phi^2+i\rho_2\nu^1+3i\delta
\rho_1\nu^2+\\
&{}&u_5\omega^1+u_6\omega^2+u_7\overline{\omega^1}+u_8\overline{\omega^2}+w_3\theta^1+w_4\theta^2,\qquad
(2.24)
\end{eqnarray*}
where $t_{\alpha},v_{\alpha},u_{\alpha},w_{\alpha}$ are complex-valued
functions, and $t_{\alpha},u_{\alpha}$ are fixed uniquely.
It now follows from (2.9), (2.15), (2.23) that
\begin{eqnarray*}
\mu_2^2&=&\mu_1^1+(a_1-i\overline{t_7})\omega^1+(b_1-i\overline{t_8})\omega^2+(c_1-i\overline{t_5})\overline{\omega^1}+(d_1-i\overline{t_6})\overline{\omega^2}+\\&{}&(h_1-i\overline{v_3})\theta^1+(g_1-i\overline{v_4})\theta^2=\\
&{}&\mu_1^1+(a_3+i\delta\overline{r_1}\rho_2+i\overline{r_2}\rho_1+i\overline{t_7})\omega^1+(b_3+i\overline{r_1}\rho_1+i\overline{r_2}\rho_2+i\overline{t_8})\omega^2+\\
&{}&(c_3+i\overline{t_5})\overline{\omega^1}+(d_3+i\overline{t_6})\overline{\omega^2}+(h_3+i\overline{v_3})\theta^1+(g_3+i\overline{v_4})\theta^2,\\
\mu_2^1&=&\delta\mu_1^2+(a_4+i\overline{t_3})\omega^1+(-\delta\overline{a_3}+i|r_1|^2-i\delta|r_2|^2+i\overline{t_4})\omega^2+(c_4+i\overline{t_1})\overline{\omega^1}+\\
&{}&(\delta
c_3+i\overline{r_1^2}+i\delta\overline{r_2^2}+i\overline{t_2})\overline{\omega^2}+(\delta
h_4+i\overline{v_1})\theta^1+(\delta
g_4+i\overline{v_2})\theta^2=\\
&{}&\delta\mu_1^2-(\overline{b_1}+i\overline{r_1}\rho_1+i\overline{r_2}\rho_2+i\overline{t_3})\omega^1+(b_2-i\overline{r_1}\rho_2-i\delta\overline{r_2}\rho_1-i\overline{t_4})\omega^2+\\
&{}&(d_1-i\overline{t_1})\overline{\omega^1}+(d_2-i\overline{t_2})\overline{\omega^2}+(h_2-i\overline{v_1})\theta^1+(g_2-i\overline{v_2})\theta^2.
\qquad\qquad
(2.25)
\end{eqnarray*}
To obtain
\begin{eqnarray*}
\mu_2^{2*}&\equiv&\mu_1^{1*}\qquad
(\hbox{mod}\,\omega^{\alpha},\overline{\omega^{\alpha}}),\\
\mu_2^{1*}&\equiv&\delta\mu_1^{2*}\qquad
(\hbox{mod}\,\omega^{\alpha},\overline{\omega^{\alpha}}),\qquad\qquad
(2.26)
\end{eqnarray*}
we set
\begin{eqnarray*}
p_1-q_2&=&h_1-i\overline{v_3}=h_3+i\overline{v_3},\\
q_1-q_4&=&g_1-i\overline{v_4}=g_3+i\overline{v_4},\\
\delta p_2-q_1&=&\delta h_4+i\overline{v_1}=h_2-i\overline{v_1},\\
\delta q_2-q_3&=&\delta g_4+i\overline{v_2}=g_2-i\overline{v_2}.\qquad\qquad\qquad\qquad
(2.27)
\end{eqnarray*}

From now on we assume that the 1-forms $\phi^{\alpha}_{\beta}$,
$\mu^{\alpha}_{\beta}$, $\psi^{\alpha}$ are chosen so that (2.7),
(2.9), (2.12), (2.16), (2.18), (2.26) are satisfied. It follows from
(2.13), (2.15), (2.25) that this set of
conditions is equivalent to (2.7), (2.9), (2.12) and
\begin{eqnarray*}
\mu_2^2&=&\mu_1^1+\tilde a_1\omega^1+\tilde b_1\omega^2+\tilde
c_1\overline{\omega^1}+\tilde d_1\overline{\omega^2},\\
\mu_2^1&=&\delta\mu_1^2+\tilde a_2\omega^1+\tilde b_2\omega^2+\tilde c_2\overline{\omega^1}+
\tilde d_2\overline{\omega^2},\\
\psi^2&=&\delta\psi^3,\\
\psi^8&=&\psi^5,\\
\psi^5&=&\psi^1+2\hbox{Re}\,\left(v_3\omega^1+v_1\omega^2\right)+2\hbox{Re}\,\left(r_2\mu_1^1+r_1\mu_1^2\right),\\
\psi^6&=&\psi^2+2\hbox{Re}\,\left(v_4\omega^1+v_2\omega^2\right),\\
\psi^7&=&\psi^3+2\hbox{Re}\,\left(\delta v_1\omega^1+v_3\omega^2\right),\\
\psi^8&=&\psi^4+2\hbox{Re}\,\left(\delta v_2\omega^1+v_4\omega^2\right)+2\hbox{Re}\,\left(\delta r_1\mu_2^1+r_2\mu_2^2\right),
\qquad\qquad
(2.28)
\end{eqnarray*}
where
\begin{eqnarray*}
\tilde a_1&:=&-i\overline{t_7}=a_3+i\delta\overline{r_1}\rho_2+i\overline{r_2}\rho_2+i\overline{t_7},\\
\tilde b_1&:=&b_1-i\overline{t_8}=b_3+i\overline{r_1}\rho_1+i\overline{r_2}\rho_2+i\overline{t_8},\\
\tilde c_1&:=&c_1-i\overline{t_5}=c_3+i\overline{t_5},\\
\tilde d_1&:=&d_1-i\overline{t_6}=d_3+i\overline{t_6},\\
\tilde a_2&:=&\delta\hbox{Im}\,(r_1\overline{r_2})+i\overline{t_3}=
-(\overline{b_1}+i\overline{r_1}\rho_1+i\overline{r_2}\rho_2+i\overline{t_3}),\\
\tilde
b_2&:=&-\delta\overline{a_3}+i|r_1|^2-i\delta|r_2|^2+i\overline{t_4}=
b_2-i\overline{r_1}\rho_2-i\delta\overline{r_2}\rho_1-i\overline{t_4},\\
\tilde c_2&:=&c_4+i\overline{t_1}=d_1-i\overline{t_1},\\
\tilde d_2&:=&\delta
c_3+i\overline{r_1^2}+i\delta\overline{r_2^2}+i\overline{t_2}=d_2-i\overline{t_2},\qquad\qquad
(2.29)
\end{eqnarray*}
are fixed, and therefore globally defined on $P^2$, functions.

Now (2.14), (2.16)-(2.19), (2.26), (2.27) give that $\phi^{\alpha}_{\beta}$,
$\mu^{\alpha}_{\beta}$, $\psi^{\alpha}$ are fixed up to
transformations of the form
\begin{eqnarray*}
\phi_1^{1*}&=&\phi_1^1+h\theta^1+g\theta^2,\\
\phi_1^{2*}&=&\phi_1^2+\delta g\theta^1+h\theta^2,\\
\phi_2^{1*}&=&\phi_2^1+g\theta^1+\delta h\theta^2,\\
\phi_2^{2*}&=&\phi_2^2+h\theta^1+g\theta^2,\\
\mu_1^{1*}&=&\mu_1^1+h\omega^1+g\omega^2+p\theta^1+q\theta^2,\\
\mu_1^{2*}&=&\mu_1^2+\delta
g\omega^1+h\omega^2+\delta q\theta^1+p\theta^2,\\
\mu_2^{1*}&=&\mu_2^1+g\omega^1+\delta h\omega^2+q\theta^1+\delta p\theta^2,\\
\mu_2^{2*}&=&\mu_2^2+h\omega^1+g\omega^2+p\theta^1+q\theta^2,\\
\psi^{1*}&=&\psi^1+2\hbox{Re}\,\left(i\overline{p}\omega^1+i\overline{q}\omega^2\right)+\sigma\theta^1+s\theta^2,\\
\psi^{2*}&=&\psi^2+2\hbox{Re}\,\left(i\overline{q}\omega^1+i\delta\overline{p}\omega^2\right)+s\theta^1+\delta\sigma\theta^2,\\
\psi^{3*}&=&\psi^3+2\hbox{Re}\,\left(i\delta\overline{q}\omega^1+i\overline{p}\omega^2\right)+\delta
s\theta^1+\sigma\theta^2,\\
\psi^{4*}&=&\psi^4+2\hbox{Re}\,\left(i\overline{p}\omega^1+i\overline{q}\omega^2\right)+\sigma\theta^1+s\theta^2,\\
\psi^{5*}&=&\psi^5+2\hbox{Re}\,\left(i\overline{p}\omega^1+i\overline{q}\omega^2\right)+\\
&{}&(\sigma+2\hbox{Re}\,(\delta
r_1q+r_2p))\theta^1+(s+2\hbox{Re}\,(r_1p+r_2q))\theta^2,\\
\psi^{6*}&=&\psi^6+2\hbox{Re}\,\left(i\overline{q}\omega^1+i\delta\overline{p}\omega^2\right)+s\theta^1+\delta\sigma\theta^2,\\
\psi^{7*}&=&\psi^7+2\hbox{Re}\,\left(i\delta\overline{q}\omega^1+i\overline{p}\omega^2\right)+\delta
s\theta^1+\sigma\theta^2,\\
\psi^{8*}&=&\psi^8+2\hbox{Re}\,\left(i\overline{p}\omega^1+i\overline{q}\omega^2\right)+\\
&{}&\bigl(\sigma+2\hbox{Re}\,(\delta
r_1q+r_2p)\bigr)\theta^1+(s+2\hbox{Re}\,(r_1p+r_2q))\theta^2,\qquad\qquad
(2.30)
\end{eqnarray*}
where $h,g$ are imaginary-valued, $\sigma, s$ are real-valued, $p,q$ are
complex-valued functions.

Consider the following matrix-valued 1-form
$$
\omega:=
\left(
\begin{array}{ccc}
\displaystyle -\frac{1}{3}\underline{\underline{\phi}}-\frac{1}{3}\underline{\phi}&\underline{\omega}&2\underline{\theta}\\
\displaystyle -i\overline{\underline{\mu}}&\frac{2}{3}\underline{\underline{\phi}}-\frac{1}{3}\underline{\phi}&2i\overline{\underline{\omega}}\\
\displaystyle -\frac{1}{4}\underline{\psi}&\frac{1}{2}\underline{\mu}&\frac{1}{3}\overline{\underline{\underline{\phi}}}+\frac{1}{3}\underline{\phi}\\
\end{array}
\right),\eqno{(2.31)}
$$
where
\begin{eqnarray*}
\underline{\theta}&:=&\left(
\begin{array}{ccc}
\displaystyle \theta^1&\delta\theta^2\\
\displaystyle \theta^2&\theta^1
\end{array}
\right),\qquad
\underline{\omega}:=
\left(
\begin{array}{ccc}
\displaystyle \omega^1&\delta\omega^2\\
\displaystyle \omega^2&\omega^1
\end{array}
\right),\\
\underline{\phi}&:=&
\left(
\begin{array}{ccc}
\displaystyle \phi^1&\delta\phi^2\\
\displaystyle \phi^2&\phi^1
\end{array}
\right),\qquad
\underline{\underline{\phi}}:=
\left(
\begin{array}{ccc}
\displaystyle \phi_1^1&\delta\phi_1^2\\
\displaystyle \phi_1^2&\phi_1^1
\end{array}
\right),\\
\underline{\mu}&:=&
\left(
\begin{array}{ccc}
\displaystyle \mu_1^1&\delta\mu_1^2\\
\displaystyle \mu_1^2&\mu_1^1
\end{array}
\right),\qquad
\underline{\psi}:=
\left(
\begin{array}{ccc}
\displaystyle \psi^1&\delta\psi^3\\
\displaystyle \psi^3&\psi^1
\end{array}
\right).
\end{eqnarray*}
It is clear from (1.3), (2.9) that the form $\omega$ defined in (2.31) takes values in
${\frak {su}}(2,1)$. Also, for any point $x\in P^2$, $\omega(x)$ is an
isomorphism between $T_x(P^2)$ and ${\frak {su}}(2,1)$. However, the
form $\omega$ is defined only locally. We now need to fix the free
parameters $h,g,p,q,\sigma,s$ from (2.30) to make the choice of
the corresponding forms unique. This will turn $\omega$ into a
globally defined ${\frak {su}}(2,1)$-valued form on $P^2$, and it will
be the parallelism that we are looking for.

To fix the free parameters from (2.30) we consider the curvature
form $\Omega$ of $\omega$
$$
\Omega:=d\omega-\frac{1}{2}[\omega,\omega]=d\omega-\omega\wedge\omega,
$$
which is a ${\frak {su}}(2,1)$-valued 2-form. In more detail, $\Omega$
is given by $\Omega=(\Omega_i^j)_{0\le i,j\le 2}$ with
$\Omega_i^j\in{\frak A}$, and
\begin{eqnarray*}
\Omega_0^0&=&-\frac{1}{3}d\underline{\underline{\phi}}-\frac{1}{3}d\underline{\phi}+i\underline{\omega}\wedge\overline{\underline{\mu}}+\frac{1}{2}\underline{\theta}\wedge\underline{\psi},\\
\Omega_1^1&=&-2i\hbox{Im}\,\Omega_0^0=\frac{2}{3}d\underline{\underline{\phi}}-\frac{1}{3}d\underline{\phi}+i\overline{\underline{\mu}}\wedge\underline{\omega}+i\underline{\mu}\wedge\overline{\underline{\omega}},\\
\Omega_2^2&=&-\overline{\Omega_0^0},\\
\Omega_0^1&=&d\underline{\omega}-\underline{\omega}\wedge\underline{\underline{\phi}}-\underline{\theta}\wedge\underline{\mu},\\
\Omega_1^2&=&2i\overline{\Omega_0^1},\\
\Omega_0^2&=&2\left(d\underline{\theta}-i\underline{\omega}\wedge\overline{\underline{\omega}}-\underline{\theta}\wedge\underline{\phi}\right),\\
\Omega_2^1&=&\frac{1}{2}\left(d\underline{\mu}+\underline{\mu}\wedge\overline{\underline{\underline{\phi}}}+\frac{1}{2}\underline{\psi}\wedge\underline{\omega}\right),\\
\Omega_1^0&=&-2i\overline{\Omega_2^1},\\
\Omega_2^0&=&-\frac{1}{4}\left(d\underline{\psi}-2i\underline{\mu}\wedge\overline{\underline{\mu}}+\underline{\psi}\wedge\underline{\phi}\right).\qquad\qquad
(2.32)
\end{eqnarray*}
Now, to fix the parameters from (2.30) we concentrate on the
components $\Omega_1^1$, $\Omega_2^1$, $\Omega_2^0$ and impose certain
conditions on their expansions. 

We start with $\Omega_1^1$. It follows from (2.32) that
$$
\Omega_1^1=
\left(
\begin{array}{ccc}
\displaystyle \Phi_1^1&\delta\Phi_1^2\\
\displaystyle \Phi_1^2&\Phi_1^1\\
\end{array}
\right),
$$
where
\begin{eqnarray*}
\Phi_1^1&:=&\frac{2}{3}d\phi_1^1-\frac{1}{3}d\phi^1+i\overline{\mu_1^1}\wedge\omega^1+i\delta\overline{\mu_1^2}\wedge\omega^2+i\mu_1^1\wedge\overline{\omega^1}+i\delta\mu_1^2\wedge\overline{\omega^2},\qquad\qquad(2.33.a)\\
\Phi_1^2&:=&\frac{2}{3}d\phi_1^2-\frac{1}{3}d\phi^2+i\overline{\mu_1^2}\wedge\omega^1+i\overline{\mu_1^1}\wedge\omega^2+i\mu_1^2\wedge\overline{\omega^1}+i\mu_1^1\wedge\overline{\omega^2}.\qquad\qquad
(2.33.b)
\end{eqnarray*}
Then (2.12), (2.21)--(2.24), (2.33) imply
\begin{eqnarray*}
\Phi_1^1&\equiv&\frac{2}{3}\Biggl(\frac{1}{2}(\rho_2-r_1)\phi^2\wedge\omega^1+\delta\frac{1}{2}(\rho_1-r_2)\phi^2\wedge\omega^2+\frac{1}{2}(\overline{r_1}-\overline{\rho_2})\phi^2\wedge\overline{\omega^1}+\\
&{}&\delta\frac{1}{2}(\overline{r_2}-\overline{\rho_1})\phi^2\wedge\overline{\omega^2}-i\rho_2\nu^2\wedge\omega^1-i\delta\rho_1\nu^2\wedge\omega^2-\\
&{}&i\overline{\rho_2}\nu^2\wedge\overline{\omega^1}-i\delta\overline{\rho_1}\nu^2\wedge\overline{\omega^2}\Biggr)\qquad
(\hbox{mod}\,\theta^{\alpha},\,\,\hbox{terms quadratic in}\,\,\omega^{\alpha},\overline{\omega^{\alpha}}),\\
\Phi_1^2&\equiv&\frac{2}{3}\Biggl(\frac{1}{2}(\rho_1-r_2)\phi^2\wedge\omega^1+\frac{1}{2}(\rho_2-r_1)\phi^2\wedge\omega^2+\frac{1}{2}(\overline{r_2}-\overline{\rho_1})\phi^2\wedge\overline{\omega^1}+\\
&{}&\frac{1}{2}(\overline{r_1}-\overline{\rho_2})\phi^2\wedge\overline{\omega^2}-i\rho_1\nu^2\wedge\omega^1-i\rho_2\nu^2\wedge\omega^2-\\
&{}&i\overline{\rho_1}\nu^2\wedge\overline{\omega^1}-i\overline{\rho_2}\nu^2\wedge\overline{\omega^2}\Biggr)\qquad
(\hbox{mod}\,\theta^{\alpha},\,\,\hbox{terms quadratic in}\,\,\omega^{\alpha},\overline{\omega^{\alpha}}).\qquad
(2.34)
\end{eqnarray*}
Let us consider the parts of the expansions of $\Phi^{\alpha}_{\beta}$
that are quadratic in $\omega^{\gamma}$, $\overline{\omega^{\gamma}}$: 
$$
\Phi^{\alpha}_{\beta}=S^{\alpha}_{\beta\gamma\tau}\omega^{\gamma}\wedge\omega^{\tau}+S^{\alpha}_{\beta\gamma\overline{\tau}}\omega^{\gamma}\wedge\overline{\omega^{\tau}}+S^{\alpha}_{\beta\overline{\gamma}\overline{\tau}}\overline{\omega^{\gamma}}\wedge\overline{\omega^{\tau}}+\dots
\eqno{(2.35)}
$$
Since $\Phi^{\alpha}_{\beta}$ are imaginary-valued, the coefficients
$S^1_{11\overline{1}},S^2_{12\overline{2}}$ are real-valued
functions. It now follows from (2.6), (2.33), (2.34) that under transformation (2.30) they change as:
\begin{eqnarray*}
S^{1*}_{11\overline{1}}&=&S^1_{11\overline{1}}+i\frac{8}{3}h,\\
S^{2*}_{12\overline{2}}&=&S^2_{12\overline{2}}+i\frac{8}{3}g.
\end{eqnarray*}
We now fix the parameters $h,g$ by the conditions
\begin{eqnarray*}
S^{1*}_{11\overline{1}}&=&0,\\
S^{2*}_{12\overline{2}}&=&0,\qquad\qquad\qquad\qquad (2.36)
\end{eqnarray*}
i.e. we set
\begin{eqnarray*}
h&=&i\frac{3}{8}S^1_{11\overline{1}},\\
g&=&i\frac{3}{8}S^2_{12\overline{2}}.
\end{eqnarray*}
Thus, the forms $\phi^{\alpha}_{\beta}$ are fixed and therefore are
defined globally on $P^2$. It then follows that the
functions $v_{\alpha}$, $w_{\alpha}$ from (2.23), (2.24), (2.28) are
also fixed and defined globally on $P^2$. 

From now on we assume that (2.7), (2.9), (2.12),
(2.28), (2.29), (2.36) are satisfied. It then follows from (2.30) that
$\mu^{\alpha}_{\beta}$, $\psi^{\alpha}$ are fixed up to
transformations of the form
\begin{eqnarray*}
\mu_1^{1*}&=&\mu_1^1+p\theta^1+q\theta^2,\\
\mu_1^{2*}&=&\mu_1^2+\delta q\theta^1+p\theta^2,\\
\mu_2^{1*}&=&\mu_2^1+q\theta^1+\delta p\theta^2,\\
\mu_2^{2*}&=&\mu_2^2+p\theta^1+q\theta^2,\\
\psi^{1*}&=&\psi^1+2\hbox{Re}\,\left(i\overline{p}\omega^1+i\overline{q}\omega^2\right)+\sigma\theta^1+s\theta^2,\\
\psi^{2*}&=&\psi^2+2\hbox{Re}\,\left(i\overline{q}\omega^1+i\delta\overline{p}\omega^2\right)+s\theta^1+\delta\sigma\theta^2,\\
\psi^{3*}&=&\psi^3+2\hbox{Re}\,\left(i\delta\overline{q}\omega^1+i\overline{p}\omega^2\right)+\delta
s\theta^1+\sigma\theta^2,\\
\psi^{4*}&=&\psi^4+2\hbox{Re}\,\left(i\overline{p}\omega^1+i\overline{q}\omega^2\right)+\sigma\theta^1+s\theta^2,\\
\psi^{5*}&=&\psi^5+2\hbox{Re}\,\left(i\overline{p}\omega^1+i\overline{q}\omega^2\right)+\\
&{}&(\sigma+2\hbox{Re}\,(\delta
r_1q+r_2p))\theta^1+(s+2\hbox{Re}\,(r_1p+r_2q))\theta^2,\\
\psi^{6*}&=&\psi^6+2\hbox{Re}\,\left(i\overline{q}\omega^1+i\delta\overline{p}\omega^2\right)+s\theta^1+\delta\sigma\theta^2,\\
\psi^{7*}&=&\psi^7+2\hbox{Re}\,\left(i\delta\overline{q}\omega^1+i\overline{p}\omega^2\right)+\delta
s\theta^1+\sigma\theta^2,\\
\psi^{8*}&=&\psi^8+2\hbox{Re}\,\left(i\overline{p}\omega^1+i\overline{q}\omega^2\right)+\\
&{}&\bigl(\sigma+2\hbox{Re}\,(\delta
r_1q+r_2p)\bigr)\theta^1+(s+2\hbox{Re}\,(r_1p+r_2q))\theta^2,\qquad\qquad
(2.37)
\end{eqnarray*}
where $\sigma, s$ are real-valued, $p,q$ are
complex-valued functions.

We now consider $\Omega_2^1$. It follows from (2.32) that
$$
\Omega_2^1=
\left(
\begin{array}{ccc}
\displaystyle \Phi^1&\delta\Phi^2\\
\displaystyle \Phi^2&\Phi^1\\
\end{array}
\right),
$$
where
\begin{eqnarray*}
\Phi^1&:=&\frac{1}{2}\Biggl(d\mu_1^1+\frac{1}{2}\mu_1^1\wedge\phi^1+\delta\frac{1}{2}\mu_1^2\wedge\phi^2+i\nu^1\wedge\mu_1^1+i\delta\nu^2\wedge\mu_1^2+\\
&{}&\frac{1}{2}\psi^1\wedge\omega^1
+\delta\frac{1}{2}\psi^3\wedge\omega^2\Biggr),\\
\Phi^2&:=&\frac{1}{2}\Biggl(d\mu_1^2+\frac{1}{2}\mu_1^2\wedge\phi^1+\frac{1}{2}\mu_1^1\wedge\phi^2+i\nu^1\wedge\mu_1^2+i\nu^2\wedge\mu_1^1+\\
&{}&\frac{1}{2}\psi^3\wedge\omega^1
+\frac{1}{2}\psi^1\wedge\omega^2\Biggr).\qquad\qquad
(2.38)
\end{eqnarray*}
To get information about the expansions of $\Phi^{\alpha}$ we return
to (2.20) and consider terms containing $\theta^{\alpha}$ there. Let
such terms in the expansions of $\phi_1^1$, $\phi_1^2$ be
\begin{eqnarray*}
d\phi_1^1&=&\lambda^1\wedge\theta^1+\lambda^2\wedge\theta^2+\dots,\\
d\phi_1^2&=&\lambda^3\wedge\theta^1+\lambda^4\wedge\theta^2+\dots,\qquad\qquad\qquad\qquad
(2.39)
\end{eqnarray*}
for some 1-forms $\lambda^{\alpha}$.
Then (2.7), (2.9), (2.20), (2.28) imply
\begin{eqnarray*}
\theta^1\wedge\Sigma^1+\delta\theta^2\wedge\Sigma^2+U&\equiv&0\qquad
(\hbox{mod}\,\,\theta^1\wedge\theta^2,\hbox{terms quardratic
in}\,\,\omega^{\alpha},\overline{\omega^{\alpha}}),\\
\theta^1\wedge\Sigma^2+\theta^2\wedge\Sigma^1+V&\equiv&0\qquad
(\hbox{mod}\,\,\theta^1\wedge\theta^2,\hbox{terms quardratic
in}\,\,\omega^{\alpha},\overline{\omega^{\alpha}}),\qquad(2.40)
\end{eqnarray*}
where
\begin{eqnarray*}
\Sigma^1&:=&d\mu_1^1+\frac{1}{2}\mu_1^1\wedge\phi^1+\delta\frac{1}{2}\mu_1^2\wedge\phi^2+i\nu^1\wedge\mu_1^1+i\delta\nu^2\wedge\mu_1^2,\\
\Sigma^2&:=&d\mu_1^2+\frac{1}{2}\mu_1^2\wedge\phi^1+\frac{1}{2}\mu_1^1\wedge\phi^2+i\nu^1\wedge\mu_1^2+i\nu^2\wedge\mu_1^1,\qquad\qquad
(2.41)
\end{eqnarray*}
and
\begin{eqnarray*}
U&:=&\theta^1\wedge\Biggl(\mu_1^2\wedge\biggl(2\hbox{Re}\,(r_1\omega^1+\delta
r_2\omega^2)+\overline{\rho_2}\overline{\omega^1}+\delta\overline{\rho_1}\overline{\omega^2}\biggr)+\omega^2\wedge\biggl(\overline{\rho_2}\overline{\mu_1^1}+\delta\overline{\rho_1}\overline{\mu_1^2}\biggr)+\\
&{}&\biggl(\tilde{a_2}\omega^1+\tilde{b_2}\omega^2+\tilde{c_2}\overline{\omega^1}+\tilde{d_2}\overline{\omega^2}\biggr)\wedge\phi^2+\omega^1\wedge\lambda^1+\delta\omega^2\wedge\lambda^3\Biggr)+\\
&{}&\theta^2\wedge\Biggl(\mu_1^1\wedge\biggl(\overline{\rho_2}\overline{\omega^1}+\delta\overline{\rho_1}\overline{\omega^2}\biggr)+2\delta\mu_1^2\wedge\hbox{Re}\,\biggl(r_2\omega^1+r_1\omega^2\biggr)+\delta\omega^2\wedge\biggl(\overline{\rho_1}\overline{\mu_1^1}+\overline{\rho_2}\overline{\mu_1^2}\biggr)+\\
&{}&\biggl(\tilde{a_2}\omega^1+\tilde{b_2}\omega^2+\tilde{c_2}\overline{\omega^1}+\tilde{d_2}\overline{\omega^2}\biggr)\wedge\phi^1+
\frac{1}{2}\biggl((\tilde{b_2}-\delta\tilde{a_1})\omega^1+\delta(\tilde{a_2}-\tilde{b_1})\omega^2+\\
&{}&(\tilde{d_2}-\delta\tilde{c_1})\overline{\omega^1}+\delta(\tilde{c_2}-\tilde{d_1})\overline{\omega^2}\biggr)\wedge\phi^2+
2i\nu^1\wedge\biggl(\tilde{c_2}\overline{\omega^1}+\tilde{d_2}\overline{\omega^2}\biggr)+\\
&{}&i\nu^2\wedge\biggl((\delta\tilde{a_1}-\tilde{b_2})\omega^1+\delta(\tilde{b_1}-\tilde{a_2})\omega^2+(\delta\tilde{c_1}+\tilde{d_2})\overline{\omega^1}+\delta(\tilde{d_1}+\tilde{c_2})\overline{\omega^2}\biggr)+\\
&{}&d\tilde{a_2}\wedge\omega^1+d\tilde{b_2}\wedge\omega^2+d\tilde{c_2}\wedge\overline{\omega^1}+d\tilde{d_2}\wedge\overline{\omega^2}+\omega^1\wedge\lambda^2+\delta\omega^2\wedge\lambda^4\Biggr),\\
V&:=&\theta^1\wedge\Biggl(2\mu_1^1\wedge\hbox{Re}\,\biggl(\delta r_1\omega^1+
r_2\omega^2\biggr)+\mu_1^2\wedge\biggl(\overline{\rho_1}\overline{\omega^1}+\overline{\rho_2}\overline{\omega^2}\biggr)+\omega^2\wedge\biggl(\overline{\rho_1}\overline{\mu_1^1}+\overline{\rho_2}\overline{\mu_1^2}\biggr)+\\
&{}&\biggl(\tilde{a_1}\omega^1+\tilde{b_1}\omega^2+\tilde{c_1}\overline{\omega^1}+\tilde{d_1}\overline{\omega^2}\biggr)\wedge\phi^2+\omega^1\wedge\lambda^3+\omega^2\wedge\lambda^1\Biggr)+\\
&{}&\theta^2\wedge\Biggl(\mu_1^1\wedge\biggl(2\hbox{Re}\,(r_2\omega^1+r_1\omega^2)+\overline{\rho_1}\overline{\omega^1}+\overline{\rho_2}\overline{\omega^2}\biggr)+\omega^2\wedge\biggl(\delta\overline{\rho_1}\overline{\mu_1^2}+\overline{\rho_2}\overline{\mu_1^1}\biggr)+\\
&{}&\biggl(\tilde{a_1}\omega^1+\tilde{b_1}\omega^2+\tilde{c_1}\overline{\omega^1}+\tilde{d_1}\overline{\omega^2}\biggr)\wedge\phi^1+
\frac{1}{2}\biggl((\tilde{b_1}-\tilde{a_2})\omega^1+(\delta\tilde{a_1}-\tilde{b_2})\omega^2+\\
&{}&(\tilde{d_1}-\tilde{c_2})\overline{\omega^1}+(\delta\tilde{c_1}-\tilde{d_2})\overline{\omega^2}\biggr)\wedge\phi^2+
2i\nu^1\wedge\biggl(\tilde{c_1}\overline{\omega^1}+\tilde{d_1}\overline{\omega^2}\biggr)+\\
&{}&i\nu^2\wedge\biggl((\tilde{a_2}-\tilde{b_1})\omega^1+(\tilde{b_2}-\delta\tilde{a_1})\omega^2+(\tilde{c_2}+\tilde{d_1})\overline{\omega^1}+(\tilde{d_2}+\delta\tilde{c_1})\overline{\omega^2}\biggr)+\\
&{}&d\tilde{a_1}\wedge\omega^1+d\tilde{b_1}\wedge\omega^2+d\tilde{c_1}\wedge\overline{\omega^1}+d\tilde{d_1}\wedge\overline{\omega^2}+
\omega^1\wedge\lambda^4+\omega^2\wedge\lambda^2\Biggr).\qquad\qquad
(2.42)
\end{eqnarray*}
Equations (2.40) imply
\begin{eqnarray*}
\theta^1\wedge U&\equiv&-\delta\theta^2\wedge V \qquad
(\hbox{mod}\,\,\theta^1\wedge\theta^2,\hbox{terms quadratic in}\,\, \omega^{\alpha},\overline{\omega^{\alpha}}),\\
\theta^2\wedge U&\equiv&\theta^1\wedge V \qquad
(\hbox{mod}\,\,\theta^1\wedge\theta^2,\hbox{terms quadratic in}\,\, \omega^{\alpha},\overline{\omega^{\alpha}}),
\end{eqnarray*}
which together with (2.42) gives
\begin{eqnarray*}
d\tilde{a_1}&\equiv&\lambda^4-\lambda^1-r_2\mu_1^1+r_1\mu_1^2+\tilde{a_1}\phi^1+\frac{1}{2}\biggl(\tilde{b_1}-3\tilde{a_2}\biggr)\phi^2+\\
&{}&i\biggl(\tilde{b_1}-\tilde{a_2}\biggr)\nu^2\qquad(\hbox{mod}\,\,\theta^{\alpha},\omega^{\alpha},
\overline{\omega^{\alpha}}),\\
d\tilde{b_1}&\equiv&\lambda^2-\delta\lambda^3-r_1\mu_1^1+\delta
r_2\mu_1^2+\tilde{b_1}\phi^1+\frac{1}{2}\biggl(\delta\tilde{a_1}-3\tilde{b_2}\biggr)\phi^2+\\
&{}&i\biggl(\delta\tilde{a_1}-\tilde{b_2}\biggr)\nu^2\qquad(\hbox{mod}\,\,\theta^{\alpha},\omega^{\alpha},
\overline{\omega^{\alpha}}),\\
d\tilde{c_1}&\equiv&-\biggl(\overline{r_2}+\overline{\rho_1}\biggr)\mu_1^1+\biggl(\overline{r_1}+\overline{\rho_2}\biggr)\mu_1^2+\tilde{c_1}\phi^1+\frac{1}{2}\biggl(\tilde{d_1}-3\tilde{c_2}\biggr)\phi^2-\\
&{}&2i\tilde{c_1}\nu^1-i\biggl(\tilde{d_1}+\tilde{c_2}\biggr)\nu^2\qquad(\hbox{mod}\,\,\theta^{\alpha},\omega^{\alpha},
\overline{\omega^{\alpha}}),\\
d\tilde{d_1}&\equiv&-\biggl(\overline{r_1}+\overline{\rho_2}\biggr)\mu_1^1+\delta\biggl(\overline{r_2}+\overline{\rho_1}\biggr)\mu_1^2+\tilde{d_1}\phi^1+\frac{1}{2}\biggl(\delta\tilde{c_1}-3\tilde{d_2}\biggr)\phi^2-\\
&{}&2i\tilde{d_1}\nu^1-i\biggl(\delta\tilde{c_1}+\tilde{d_2}\biggr)\nu^2\qquad(\hbox{mod}\,\,\theta^{\alpha},\omega^{\alpha},
\overline{\omega^{\alpha}}),\\
d\tilde{a_2}&\equiv&\lambda^2-\delta\lambda^3+r_1\mu_1^1-\delta
r_2\mu_1^2+\tilde{a_2}\phi^1+\frac{1}{2}\biggl(\tilde{b_2}-3\delta\tilde{a_1}\biggr)\phi^2+\\
&{}&i\biggl(\tilde{b_2}-\delta\tilde{a_1}\biggr)\nu^2\qquad(\hbox{mod}\,\,\theta^{\alpha},\omega^{\alpha},
\overline{\omega^{\alpha}}),\\
d\tilde{b_2}&\equiv&\delta\biggl(\lambda^4-\lambda^1)+\delta
r_2\mu_1^1-\delta
r_1\mu_1^2+\tilde{b_2}\phi^1+\delta\frac{1}{2}\biggl(\tilde{a_2}-3\tilde{b_1}\biggr)\phi^2+\\
&{}&i\delta\biggl(\tilde{a_2}-\tilde{b_1}\biggr)\nu^2\qquad(\hbox{mod}\,\,\theta^{\alpha},\omega^{\alpha},
\overline{\omega^{\alpha}}),\\
d\tilde{c_2}&\equiv&\biggl(\overline{r_1}-\overline{\rho_2}\biggr)\mu_1^1+\delta\biggl(\overline{\rho_1}-\overline{r_2}\biggr)\mu_1^2+\tilde{c_2}\phi^1+\frac{1}{2}\biggl(\tilde{d_2}-3\delta\tilde{c_1}\biggr)\phi^2-\\
&{}&2i\tilde{c_2}\nu^1-i\biggl(\tilde{d_2}+\delta\tilde{c_1}\biggr)\nu^2\qquad(\hbox{mod}\,\,\theta^{\alpha},\omega^{\alpha},
\overline{\omega^{\alpha}}),\\
d\tilde{d_2}&\equiv&\delta\biggl(\overline{r_2}-\overline{\rho_1}\biggr)\mu_1^1+\delta\biggl(\overline{\rho_2}-\overline{r_1}\biggr)\mu_1^2+\tilde{d_2}\phi^1+\delta\frac{1}{2}\biggl(\tilde{c_2}-3\tilde{d_1}\biggr)\phi^2-\\
&{}&2i\tilde{d_2}\nu^1-i\delta\biggl(\tilde{d_1}+\tilde{c_2}\biggr)\nu^2\qquad(\hbox{mod}\,\,\theta^{\alpha},\omega^{\alpha},
\overline{\omega^{\alpha}}).\qquad\qquad (2.43)
\end{eqnarray*}
It now follows from (2.40), (2.42), (2.43) that
\begin{eqnarray*}
\Sigma^1&\equiv&\biggl(\lambda^1-r_1\mu_1^2+\tilde{a_2}\phi^2\biggr)\wedge\omega^1+\biggl(\delta\lambda^3-\delta
r_2\mu_1^2+\overline{\rho_2}\overline{\mu_1^1}+\delta\overline{\rho_1}\overline{\mu_1^2}+\tilde{b_2}\phi^2\biggr)\wedge\omega^2+\\
&{}&\biggl(-(\overline{r_1}+\overline{\rho_2})\mu_1^2+\tilde{c_2}\phi^2\biggr)\wedge\overline{\omega^1}+\biggl(-\delta(\overline{r_2}+\overline{\rho_1})\mu_1^2+\tilde{d_2}\phi^2\biggr)\wedge\overline{\omega^2}\\
&{}&(\hbox{mod}\,\,\theta^{\alpha},\hbox{terms
quadratic in}\,\,\omega^{\alpha},\overline{\omega^{\alpha}}),\\
\Sigma^2&\equiv&\biggl(\lambda^3-\delta
r_1\mu_1^1+\tilde{a_1}\phi^2\biggr)\wedge\omega^1+\biggl(\lambda^1-\delta
r_2\mu_1^1+\delta\overline{\rho_1}\overline{\mu_1^1}+\delta\overline{\rho_2}\overline{\mu_1^2}+\tilde{b_1}\phi^2\biggr)\wedge\omega^2+\\
&{}&\biggl(-\delta\overline{r_1}\mu_1^1-\overline{\rho_1}\mu_1^2+\tilde{c_1}\phi^2\biggr)\wedge\overline{\omega^1}+\biggl(-\overline{r_2}\mu_1^1-\overline{\rho_2}\mu_1^2+\tilde{d_1}\phi^2\biggr)\wedge\overline{\omega^2}\\
&{}&(\hbox{mod}\,\,\theta^{\alpha},\hbox{terms
quadratic
in}\,\,\omega^{\alpha},\overline{\omega^{\alpha}}).\qquad\qquad (2.44)
\end{eqnarray*}
Formulas (2.38), (2.41), (2.44) give
\begin{eqnarray*}
\Phi^1&\equiv&\frac{1}{2}\Biggl(
\biggl(\lambda^1-r_1\mu_1^2+\tilde{a_2}\phi^2+\frac{1}{2}\psi^1\biggr)\wedge\omega^1+\biggl(\delta\lambda^3-\delta
r_2\mu_1^2+\overline{\rho_2}\overline{\mu_1^1}+\delta\overline{\rho_1}\overline{\mu_1^2}+\\
&{}&\tilde{b_2}\phi^2+\delta\frac{1}{2}\psi^3\biggr)\wedge\omega^2+
\biggl(-(\overline{r_1}+\overline{\rho_2})\mu_1^2+\tilde{c_2}\phi^2\biggr)\wedge\overline{\omega^1}+\\
&{}&\biggl(-\delta(\overline{r_2}+\overline{\rho_1})\mu_1^2+\tilde{d_2}\phi^2\biggr)\wedge\overline{\omega^2}\Biggr)\\
&{}&(\hbox{mod}\,\,\theta^{\alpha},\hbox{terms
quadratic in}\,\,\omega^{\alpha},\overline{\omega^{\alpha}}),\\
\Phi^2&\equiv&\frac{1}{2}\Biggl(\biggl(\lambda^3-\delta
r_1\mu_1^1+\tilde{a_1}\phi^2+\frac{1}{2}\psi^3\biggr)\wedge\omega^1+\biggl(\lambda^1-\delta
r_2\mu_1^1+\delta\overline{\rho_1}\overline{\mu_1^1}+\delta\overline{\rho_2}\overline{\mu_1^2}+\\
&{}&\tilde{b_1}\phi^2+\frac{1}{2}\psi^1\biggr)\wedge\omega^2+\biggl(-\delta\overline{r_1}\mu_1^1-\overline{\rho_1}\mu_1^2+\tilde{c_1}\phi^2\biggr)\wedge\overline{\omega^1}+\\
&{}&\biggl(-\overline{r_2}\mu_1^1-\overline{\rho_2}\mu_1^2+\tilde{d_1}\phi^2\biggr)\wedge\overline{\omega^2}\Biggr)\\
&{}&(\hbox{mod}\,\,\theta^{\alpha},\hbox{terms
quadratic
in}\,\,\omega^{\alpha},\overline{\omega^{\alpha}}).\qquad (2.45)
\end{eqnarray*}

We will now show that the expansions of $\Phi^{\alpha}$ do not contain
terms involving $\omega^{\beta}\wedge\psi^{\gamma}$. It is clear from
(2.45) that for this we only need to prove that the expansions of
$\lambda^1$, $\lambda^3$ have the form
\begin{eqnarray*}
\lambda^1&=&-\frac{1}{2}\psi^1+\hbox{terms not containing
$\psi^{\alpha}$},\qquad\qquad\qquad\qquad (2.46.a)\\
\lambda^3&=&-\frac{1}{2}\psi^3+\hbox{terms not containing
$\psi^{\alpha}$}.
\qquad\qquad\qquad\qquad (2.46.b)
\end{eqnarray*}
Let
\begin{eqnarray*}
\lambda^1&=&\chi_1\psi^1+\chi_2\psi^3+\hbox{terms not containing
$\psi^{\alpha}$},\\
\lambda^3&=&\chi_3\psi^1+\chi_4\psi^3+\hbox{terms not containing
$\psi^{\alpha}$}.\qquad\qquad\qquad\qquad (2.47)
\end{eqnarray*}
Identities (2.41), (2.44) give
\begin{eqnarray*}
d\mu_1^1&=&\biggl(\chi_1\psi^1+\chi_2\psi^3\biggr)\wedge\omega^1+\delta\biggl(\chi_3\psi^1+\chi_4\psi^3\biggr)\wedge\omega^2+\\
&{}&\hbox{terms
not containing $\psi^{\alpha}$},\\
d\mu_1^2&=&\biggl(\chi_3\psi^1+\chi_4\psi^3\biggr)\wedge\omega^1+\biggl(\chi_1\psi^1+\chi_2\psi^3\biggr)\wedge\omega^2+\\
&{}&\hbox{terms
not containing $\psi^{\alpha}$}.\qquad(2.48)
\end{eqnarray*}

We now differentiate (2.33.a) and, using (2.7), (2.9), (2.28), (2.48),
in the right-hand side of the resulting equation collect terms
containing $\omega^1\wedge\overline{\omega^1}\wedge\psi^{\alpha}$:
$$
-2\omega^1\wedge\overline{\omega^1}\wedge\biggl(\hbox{Im}\,\chi_1\psi^1+\hbox{Im}\,\chi_2\psi^3\biggr).\eqno{(2.49)}
$$
On the other hand, it follows from (2.12), (2.28), (2.33)--(2.36),
(2.39), (2.47) that such terms in the left-hand side are:
$$
-i\omega^1\wedge\overline{\omega^1}\wedge\left(\left(\frac{2}{3}\chi_1+\frac{1}{3}\right)\psi^1+\frac{2}{3}\chi_2\psi^3\right).\eqno{(2.50)}
$$
Comparing (2.49), (2.50) we get
$$
\chi_1=-\frac{1}{2},\qquad\chi_2=0,
$$
and (2.46.a) is proved.

Similarly, we differentiate (2.33.b) and in the right-hand side of the
resulting equation collect terms containing
$\omega^2\wedge\overline{\omega^2}\wedge\psi^{\alpha}$:
$$
-2\delta\omega^2\wedge\overline{\omega^2}\wedge\biggl(\hbox{Im}\,\chi_3\psi^1+\hbox{Im}\,\chi_4\psi^3\biggr).\eqno{(2.51)}
$$
Such terms in the left-hand side are:
$$
-i\delta\omega^2\wedge\overline{\omega^2}\wedge\left(\frac{2}{3}\chi_3\psi^1+
\left(\frac{2}{3}\chi_4+\frac{1}{3}\right)\psi^3\right).\eqno{(2.52)}
$$
Comparing (2.51), (2.52) we obtain
$$
\chi_3=0,\qquad\chi_4=-\frac{1}{2},
$$
and (2.46.b) is proved.

We are now ready to fix the parameters $p,q$ in (2.37). Let us
consider the parts of the expansions of $\Phi^{\alpha}$ that are quadratic in
$\omega^{\beta},\overline{\omega^{\beta}}$:
$$
\Phi^{\alpha}=S^{\alpha}_{\beta\gamma}\omega^{\beta}\wedge\omega^{\gamma}+S^{\alpha}_{\beta\overline{\gamma}}\omega^{\beta}\wedge\overline{\omega^{\gamma}}+S^{\alpha}_{\overline{\beta}\overline{\gamma}}\overline{\omega^{\beta}}\wedge
\overline{\omega^{\gamma}}+\dots
\eqno{(2.53)}
$$
It follows from (2.6), (2.38), (2.45), (2.46) that under transformation
(2.37) the coefficients $S^1_{1\overline{1}}, S^2_{2\overline{2}}$
change as:
\begin{eqnarray*}
S^{1*}_{1\overline{1}}&=&S^1_{1\overline{1}}+i\frac{3}{4}p,\\
S^{2*}_{2\overline{2}}&=&S^2_{2\overline{2}}+i\frac{3}{4}q.
\end{eqnarray*}
We now fix $p,q$ by the conditions:
\begin{eqnarray*}
S^{1*}_{1\overline{1}}&=&0,\\
S^{2*}_{2\overline{2}}&=&0,
\qquad\qquad\qquad\qquad (2.54)
\end{eqnarray*}
i.e. we set
\begin{eqnarray*}
p&=&i\frac{4}{3}S^1_{1\overline{1}},\\
q&=&i\frac{4}{3}S^2_{2\overline{2}}.
\end{eqnarray*}
The forms $\mu^{\alpha}_{\beta}$ are now fixed and thus are globally
defined on $P^2$.

From now on we assume that (2.7), (2.9), (2.12),
(2.28), (2.29), (2.36), (2.54) are satisfied. It then follows from (2.37) that
$\psi^{\alpha}$ are fixed up to transformations of the form
\begin{eqnarray*}
\psi^{1*}&=&\psi^1+\sigma\theta^1+s\theta^2,\\
\psi^{2*}&=&\psi^2+s\theta^1+\delta\sigma\theta^2,\\
\psi^{3*}&=&\psi^3+\delta
s\theta^1+\sigma\theta^2,\\
\psi^{4*}&=&\psi^4+\sigma\theta^1+s\theta^2,\\
\psi^{5*}&=&\psi^5+\sigma\theta^1+s\theta^2,\\
\psi^{6*}&=&\psi^6+s\theta^1+\delta\sigma\theta^2,\\
\psi^{7*}&=&\psi^7+\delta
s\theta^1+\sigma\theta^2,\\
\psi^{8*}&=&\psi^8+\sigma\theta^1+s\theta^2,\qquad\qquad
(2.55)
\end{eqnarray*}
where $\sigma, s$ are real-valued functions.

To fix the parameters $\sigma,s$ in (2.55) we consider
$\Omega^0_2$. It follows from (2.32) that
$$
\Omega^0_2=
\left(
\begin{array}{ccc}
\displaystyle \Psi^1&\delta\Psi^2\\
\displaystyle \Psi^2&\Psi^1\\
\end{array}
\right),
$$
where
\begin{eqnarray*}
\Psi^1&:=&-\frac{1}{4}\Biggl(d\psi^1-2i\biggl(\mu_1^1\wedge\overline{\mu_1^1}+\delta\mu_1^2\wedge\overline{\mu_1^2}\Biggr)+\psi^1\wedge\phi^1+\delta\psi^3\wedge\phi^2\Biggr),\\
\Psi^2&:=&-\frac{1}{4}\Biggl(d\psi^3-2i\biggl(\mu_1^2\wedge\overline{\mu_1^1}+\mu_1^1\wedge\overline{\mu_1^2}\Biggr)+\psi^3\wedge\phi^1+\psi^1\wedge\phi^2\Biggr).\qquad\qquad
(2.56)
\end{eqnarray*}
Let us consider the parts of the expansions of $\Psi^{\alpha}$ that are
quadratic in $\omega^{\beta},\overline{\omega^{\beta}}$:
$$
\Psi^{\alpha}=T^{\alpha}_{\beta\gamma}\omega^{\beta}\wedge\omega^{\gamma}+T^{\alpha}_{\beta\overline{\gamma}}\omega^{\beta}\wedge\overline{\omega^{\gamma}}+T^{\alpha}_{\overline{\beta}\overline{\gamma}}\overline{\omega^{\beta}}\wedge
\overline{\omega^{\gamma}}+\dots
\eqno{(2.57)}
$$
Since $\Psi^{\alpha}$ are real-valued, the coefficients
$T^1_{1\overline{1}},T^2_{2\overline{2}}$ are imaginary-valued. It
follows from (2.6), (2.56) that under
transformation (2.55) they change as:
\begin{eqnarray*}
T^{1*}_{1\overline{1}}&=&T^1_{1\overline{1}}-i\frac{1}{4}\sigma,\\
T^{2*}_{2\overline{2}}&=&T^2_{2\overline{2}}-i\frac{1}{4}s.
\end{eqnarray*}
We now fix $\sigma,s$ by the conditions:
\begin{eqnarray*}
T^{1*}_{1\overline{1}}&=&0,\\
T^{2*}_{2\overline{2}}&=&0,\qquad\qquad\qquad\qquad (2.58)
\end{eqnarray*}
i.e. we set
\begin{eqnarray*}
\sigma&=&-4iT^1_{1\overline{1}},\\
s&=&-4iT^2_{2\overline{2}}.
\end{eqnarray*}

We have proved that the forms
$\phi^{\alpha}_{\beta},\mu^{\alpha}_{\beta},\psi^{\alpha}$ are
uniquely fixed by conditions (2.7), (2.9), (2.12),
(2.28), (2.29), (2.36), (2.54), (2.58), and therefore we can assume that they
are now defined globally on $P^2$. The form $\omega$ defined via these
forms as in (2.31) is the parallelism that we needed to construct. We
also note that the functions $r_{\alpha}$ from (2.6), $\rho_{\alpha}$
from (2.9), $t_{\alpha},v_{\alpha}$ from (2.23),
$u_{\alpha},w_{\alpha}$ from (2.24),
$\tilde{a_{\alpha}},\tilde{b_{\alpha}},\tilde{c_{\alpha}},\tilde{d_{\alpha}}$
from (2.28) are $CR$-invariant functions, i.e. scalar invariants.

In the remainder of this section we will find a transformation formula
for $\omega$ under the action of $G^2$ on $P^2$. Let $\eta\in G^2$ be
given by matrix (1.2), where $C\overline{C}=E$ and
$$
C=
\left(
\begin{array}{ccc}
\displaystyle C_1&\delta C_2\\
\displaystyle C_2&C_1
\end{array}
\right),\qquad
A=
\left(
\begin{array}{ccc}
\displaystyle A_1&\delta A_2\\
\displaystyle A_2&A_1
\end{array}
\right),\qquad
R=
\left(
\begin{array}{ccc}
\displaystyle R_1&\delta R_2\\
\displaystyle R_2&R_1
\end{array}
\right),
$$
with $C_{\alpha},A_{\alpha}\in\CC$, $R_{\alpha}\in\RR$. Let $L_{\eta}$
denote the mapping of $P^2$ induced by $\chi(\eta)$ (see (1.5)).

It turns out that one can find the
transformation rule for $\omega$ under $L_{\eta}$ in the form
$$
L_{\eta}^{*}\omega=\hbox{Ad}(\eta)\omega+
\left(
\begin{array}{ccc}
   \displaystyle -\frac{1}{3}\Bigl(\Pi\overline{\underline{\omega}}-\overline{\Pi}\underline{\omega}+\Gamma\Bigr)&0&0\\
   \displaystyle -i(\overline{M}+\overline{P}+\overline{\Delta})&\frac{2}{3}\Bigl(\Pi\overline{\underline{\omega}}-\overline{\Pi}\underline{\omega}+\Gamma\Bigr)&0\\
\displaystyle -\frac{1}{4}\Bigl(N+\overline{N}+\Theta+\Lambda\Bigr)&\frac{1}{2}\Bigl(M+P+\Delta\Bigr)&\frac{1}{3}\Bigl(\overline{\Pi}\underline{\omega}-\Pi\overline{\underline{\omega}}-\Gamma\Bigr),\\
\end{array}
\right), \eqno{(2.59)}
$$
where $\hbox{Ad}$ is the adjoint representation of $\hbox{Aut}_{e}(Q_H)$
on ${\frak {su}}(2,1)$ and 
\begin{eqnarray*}
\Pi&:=&\left(
\begin{array}{ccc}
\displaystyle \Pi_1&\delta\Pi_2\\
\displaystyle \Pi_2&\Pi_1\\
\end{array}
\right),\qquad
\Gamma:=\left(
\begin{array}{ccc}
\displaystyle
\Gamma_1\theta^1+\Gamma_2\theta^2&\delta(\Gamma_3\theta^1+\Gamma_4\theta^2)\\
\displaystyle
\Gamma_3\theta^1+\Gamma_4\theta^2&\Gamma_1\theta^1+\Gamma_2\theta^2\\
\end{array}
\right),\\
M&:=&\left(
\begin{array}{ccc}
\displaystyle
M_1\omega^1+M_2\omega^2&\delta(M_3\omega^1+M_4\omega^2)\\
\displaystyle M_3\omega^1+M_4\omega^2&M_1\omega^1+M_2\omega^2\\
\end{array}
\right),\qquad
P:=\left(
\begin{array}{ccc}
\displaystyle
P_1\overline{\omega^1}+P_2\overline{\omega^2}&\delta(P_3\overline{\omega^1}+P_4\overline{\omega^2})\\
\displaystyle P_3\overline{\omega^1}+P_4\overline{\omega^2}&P_1\overline{\omega^1}+P_2\overline{\omega^2}\\
\end{array}
\right),\\
\Delta&:=&\left(
\begin{array}{ccc}
\displaystyle
\Delta_1\theta^1+\Delta_2\theta^2&\delta(\Delta_3\theta^1+\Delta_4\theta^2)\\
\displaystyle
\Delta_3\theta^1+\Delta_4\theta^2&\Delta_1\theta^1+\Delta_2\theta^2\\
\end{array}
\right),\qquad
N:=\left(
\begin{array}{ccc}
\displaystyle
N_1\omega^1+N_2\omega^2&\delta(N_3\omega^1+N_4\omega^2)\\
\displaystyle N_3\omega^1+N_4\omega^2&N_1\omega^1+N_2\omega^2\\
\end{array}
\right),\\
\Theta&:=&\left(
\begin{array}{ccc}
\displaystyle
\Theta_1\theta^1+\Theta_2\theta^2&\delta(\Theta_3\theta^1+\Theta_4\theta^2)\\
\displaystyle
\Theta_3\theta^1+\Theta_4\theta^2&\Theta_1\theta^1+\Theta_2\theta^2\\
\end{array}
\right),\qquad
\Lambda:=\left(
\begin{array}{ccc}
\displaystyle
\Lambda_1\phi^2&\delta\Lambda_2\phi^2\\
\displaystyle
\Lambda_2\phi^2&\Lambda_1\phi^2\\
\end{array}
\right),\qquad\qquad\qquad (2.60)
\end{eqnarray*}
with $\Pi_{\alpha},M_{\alpha},P_{\alpha},\Delta_{\alpha},N_{\alpha}$
complex-valued, $\Gamma_{\alpha}$ imaginary-valued, and
$\Theta_{\alpha},\Lambda_{\alpha}$ real-valued functions on $P^2$.

To determine the
parameters in (2.60) we plug the right-hand side of (2.59) in (2.6), (2.7),
(2.9), (2.12), (2.16), (2.18), (2.26), (2.36), (2.54), (2.58). The
computations are elementary, but lengthy, and we only list
the final results here.

Plugging in (2.6) we get
$$
\left(
\begin{array}{ccc}
\displaystyle r_1^{*}\\
\displaystyle r_2^{*}\\
\end{array}
\right):=
L_{\eta}^{*}
\left(
\begin{array}{ccc}
\displaystyle r_1\\
\displaystyle r_2\\
\end{array}
\right)=
C\left(
\begin{array}{ccc}
\displaystyle r_1\\
\displaystyle r_2\\
\end{array}
\right).
$$

Plugging in (2.7) gives
\begin{eqnarray*}
\Pi_1&=&\hbox{det}\,C\overline{C_2\left(C_1\rho_2+\delta
C_2\rho_1\right)},\\
\Pi_2&=&\hbox{det}\,C\overline{C_2\left(C_1\rho_1+
C_2\rho_2\right)},\\
\rho_1^{*}&:=&L_{\eta}^{*}\rho_1=\left(C_1|C_1|^2-2\delta C_1|C_2|^2+\delta\overline{C_1}C_2^2\right)\rho_1+
\left(2C_2|C_1|^2-\delta
C_2|C_2|^2-\overline{C_2}C_1^2\right)\rho_2,\\
\rho_2^{*}&:=&L_{\eta}^{*}\rho_2=\left(2\delta C_2|C_1|^2-C_2|C_2|^2-\delta\overline{C_2}C_1^2\right)\rho_1+
\left(C_1|C_1|^2-2\delta
C_1|C_2|^2+\delta\overline{C_1}C_2^2\right)\rho_2,\\
P_1&=&-2\Biggl(\delta
A_2C_2\overline{\rho_1^{*}}+A_2C_1\overline{\rho_2^{*}}+A_2\overline{r_1}-A_1\Pi_1-\delta
A_2\Pi_2\Biggr),\\
P_2&=&-2\delta\Biggl(
A_2C_1\overline{\rho_1^{*}}+A_2C_2\overline{\rho_2^{*}}+A_2\overline{r_2}-A_1\Pi_2-
A_2\Pi_1\Biggr),\\
P_3&=&-2\Biggl(
A_2C_1\overline{\rho_1^{*}}+A_2C_2\overline{\rho_2^{*}}+\delta A_1\overline{r_1}-A_1\Pi_2-
A_2\Pi_1\Biggr),\\
P_4&=&-2\Biggl(\delta
A_2C_2\overline{\rho_1^{*}}+A_2C_1\overline{\rho_2^{*}}+A_1\overline{r_2}-A_1\Pi_1-\delta
A_2\Pi_2\Biggr),\\
M_1&=&2\Biggl(\delta
\overline{A_2}\overline{C_2}\overline{\rho_1^{*}}+\overline{A_1}\overline{C_2}\overline{\rho_2^{*}}-A_2r_1-A_1\overline{\Pi_1}-\delta
A_2\overline{\Pi_2}\Biggr)+\Gamma_1\overline{C_1}+\delta\Gamma_3\overline{C_2},\\
M_2&=&2\Biggl(\delta
\overline{A_2}\overline{C_1}\overline{\rho_1^{*}}+\overline{A_1}\overline{C_1}\overline{\rho_2^{*}}-\delta
A_2r_2-\delta A_1\overline{\Pi_2}-\delta
A_2\overline{\Pi_1}\Biggr)+\delta\Gamma_1\overline{C_2}+\delta\Gamma_3\overline{C_1},\\
M_3&=&2\Biggl(
\overline{A_1}\overline{C_2}\overline{\rho_1^{*}}+\overline{A_2}\overline{C_2}\overline{\rho_2^{*}}-\delta
A_1r_1-A_2\overline{\Pi_1}-
A_1\overline{\Pi_2}\Biggr)+\Gamma_1\overline{C_2}+\Gamma_3\overline{C_1},\\
M_4&=&2\Biggl(
\overline{A_1}\overline{C_1}\overline{\rho_1^{*}}+\overline{A_2}\overline{C_1}\overline{\rho_2^{*}}-A_1r_2-A_1\overline{\Pi_1}-\delta
A_2\overline{\Pi_2}\Biggr)+\Gamma_1\overline{C_1}+\delta\Gamma_3\overline{C_2},\\
\tilde{a_1^{*}}&:=&L_{\eta}^{*}\tilde{a_1}=\hbox{det}\,C\Biggl(\overline{C_1^2}\tilde{a_1}-\overline{C_1}\overline{C_2}\tilde{b_1}+\overline{C_1}\overline{C_2}\tilde{a_2}-\overline{C_2^2}\tilde{b_2}+2\Bigl((A_1\overline{C_2}+A_2\overline{C_1})r_1-\\
&{}&(A_1\overline{C_1}+\delta
A_2\overline{C_2})r_2\Bigr)\Biggr)-\Gamma_1+\Gamma_4,\\
\tilde{b_1^{*}}&:=&L_{\eta}^{*}\tilde{b_1}=\hbox{det}\,C\Biggl(-\delta\overline{C_1}\overline{C_2}\tilde{a_1}+\overline{C_1^2}\tilde{b_1}-\delta\overline{C_2^2}\tilde{a_2}+\overline{C_1}\overline{C_2}\tilde{b_2}+2\Bigl(-(A_1\overline{C_1}+\delta
A_2\overline{C_2})r_1+\\
&{}&\delta(A_1\overline{C_2}+
A_2\overline{C_1})r_2\Bigr)\Biggr)+\Gamma_2-\delta\Gamma_3,\\
\tilde{c_1^{*}}&:=&L_{\eta}^{*}\tilde{c_1}=\hbox{det}\,\overline{C}\Biggl(|C_1|^2\tilde{c_1}-\overline{C_1}C_2\tilde{d_1}+\overline{C_2}C_1\tilde{c_2}-|C_2|^2\tilde{d_2}+2\Bigl((A_1C_2+A_2C_1)\overline{r_1}-\\
&{}&(A_1C_1+\delta
A_2C_2)\overline{r_2}\Bigr)\Biggr)-2\Biggl(A_1\overline{\rho_1^{*}}-A_2\overline{\rho_2^{*}}\Biggr),\\
\tilde{d_1^{*}}&:=&L_{\eta}^{*}\tilde{d_1}=\hbox{det}\,\overline{C}\Biggl(-\delta\overline{C_1}C_2\tilde{c_1}+|C_1|^2\tilde{d_1}-\delta|C_2|^2\tilde{c_2}+\overline{C_2}C_1\tilde{d_2}+2\Bigl(-(A_1C_1+\delta
A_2C_2)\overline{r_1}+\\
&{}&\delta(A_1C_2+
A_2C_1)\overline{r_2}\Bigr)\Biggr)+2\Biggl(\delta A_2\overline{\rho_1^{*}}-A_1\overline{\rho_2^{*}}\Biggr),\\
\tilde{a_2^{*}}&:=&L_{\eta}^{*}\tilde{a_2}=\hbox{det}\,C\Biggl(\delta\overline{C_1}\overline{C_2}\tilde{a_1}-\delta\overline{C_2^2}\tilde{b_1}+\overline{C_1^2}\tilde{a_2}-\overline{C_1}\overline{C_2}\tilde{b_2}+2\Bigl((A_1\overline{C_1}+\delta
A_2\overline{C_2})r_1-\\
&{}&\delta(A_1\overline{C_2}+
A_2\overline{C_1})r_2\Bigr)\Biggr)+\Gamma_2-\delta \Gamma_3,\\
\tilde{b_2^{*}}&:=&L_{\eta}^{*}\tilde{b_2}=\hbox{det}\,C\Biggl(-\overline{C_2^2}\tilde{a_1}+\delta\overline{C_1}\overline{C_2}\tilde{b_1}-\delta\overline{C_1}\overline{C_2}\tilde{a_2}+\overline{C_1^2}\tilde{b_2}+2\Bigl(-\delta(A_1
\overline{C_2}+A_2\overline{C_1})r_1+\\
&{}&(\delta A_1\overline{C_1}+
A_2\overline{C_2})r_2\Bigr)\Biggr)-\delta\Gamma_1+\delta\Gamma_4,\\
\tilde{c_2^{*}}&:=&L_{\eta}^{*}\tilde{c_2}=\hbox{det}\,\overline{C}\Biggl(\delta\overline{C_2}C_1\tilde{c_1}-\delta|C_2|^2\tilde{d_1}+|C_1|^2\tilde{c_2}-\overline{C_1}C_2\tilde{d_2}+2\Bigl((A_1C_1+\delta
A_2C_2)\overline{r_1}-\\
&{}&\delta(A_1C_2+
A_2C_1)\overline{r_2}\Bigr)\Biggr)+2\Biggl(\delta
A_2\overline{\rho_1^{*}}-A_1\overline{\rho_2^{*}}\Biggr),\\
\tilde{d_2^{*}}&:=&L_{\eta}^{*}\tilde{d_2}=\hbox{det}\,\overline{C}\Biggl(-|C_2|^2\tilde{c_1}+\delta\overline{C_2}C_1\tilde{d_1}-\delta\overline{C_1}C_2\tilde{c_2}+|C_1|^2\tilde{d_2}+2\Bigl(-\delta(A_1C_2+A_2C_1)\overline{r_1}+\\
&{}&(\delta
A_1C_1+
A_2C_2)\overline{r_2}\Bigr)\Biggr)-2\delta\Biggl(A_1\overline{\rho_1^{*}}-A_2\overline{\rho_2^{*}}\Biggr),\\
\delta\Delta_3-\Delta_2&=&4\delta\Bigl(A_1\overline{A_2}-A_2\overline{A_1}\Bigr)\overline{\rho_1^{*}}+4\Bigl(|A_1|^2-\delta|A_2|^2\Bigr)\overline{\rho_2^{*}}+2\Bigl(A_1\tilde{a_2^{*}}+A_2\tilde{b_2^{*}}+\overline{A_1}\tilde{c_2^{*}}+\\
&{}&\overline{A_2}\tilde{d_2^{*}}-A_1\Gamma_2+\delta
A_2\Gamma_1+\delta
A_1\Gamma_3-\delta A_2\Gamma_4\Bigr),\\
\Delta_1-\Delta_4&=&4\Bigl(|A_1|^2-\delta|A_2|^2\Bigr)\overline{\rho_1^{*}}+4\Bigl(A_1\overline{A_2}-A_2\overline{A_1}\Bigr)\overline{\rho_2^{*}}+2\Bigl(A_1\tilde{a_1^{*}}+A_2\tilde{b_1^{*}}+\overline{A_1}\tilde{c_1^{*}}+\\
&{}&\overline{A_2}\tilde{d_1^{*}}+A_1\Gamma_1-
A_2\Gamma_2-
A_1\Gamma_4+\delta A_2\Gamma_3\Bigr).\qquad\qquad\qquad (2.61)
\end{eqnarray*}
It follows from (2.61) that $\delta\Delta_3-\Delta_2$ and
$\Delta_1-\Delta_4$ do not in fact depend on $\Gamma_{\alpha}$.

Equations (2.9) are satisfied automatically, and we now plug the
right-hand side of (2.59) in (2.12). This gives
\begin{eqnarray*}
\Lambda_1&=&4\hbox{Re}\,\Bigl(A_1r_1^{*}+\delta A_2r_2^{*}\Bigr),\\
\Lambda_2&=&4\hbox{Re}\,\Bigl(A_1r_2^{*}+A_2r_1^{*}\Bigr),\\
N_1&=&-8\Bigl(R_2r_1+i(\overline{C_2}\overline{A_1}+\overline{C_1}\overline{A_2})\hbox{Re}\,(A_1r_1^{*}+\delta
A_2r_2^{*})\Bigr)+i\Bigl(\overline{C_1}\overline{\Delta_1}+\delta\overline{C_2}\overline{\Delta_3}+\\
&{}&2(A_1\overline{P_1}+\delta
A_2\overline{P_3}-\overline{A_1}M_1-\delta\overline{A_2}M_3)\Bigr),\\
N_2&=&-8\Bigl(\delta R_2r_2+i(\overline{C_1}\overline{A_1}+\delta\overline{C_2}\overline{A_2})\hbox{Re}\,(A_1r_1^{*}+\delta
A_2r_2^{*})\Bigr)+i\Bigl(\delta\overline{C_1}\overline{\Delta_3}+\delta\overline{C_2}\overline{\Delta_1}+\\
&{}&2(A_1\overline{P_2}+\delta
A_2\overline{P_4}-\overline{A_1}M_2-\delta\overline{A_2}M_4)\Bigr),\\
N_3&=&-4\Bigl(\delta R_1r_1+R_2r_2+2i(\overline{C_2}\overline{A_1}+\overline{C_1}\overline{A_2})\hbox{Re}\,(A_1r_2^{*}+
A_2r_1^{*})\Bigr)+i\Biggl(\delta\overline{C_1}\overline{\Delta_2}+\overline{C_2}\overline{\Delta_4}+\\
&{}&2\Bigl(A_1\overline{P_3}+
A_2\overline{P_1}-\overline{A_1}M_3-\overline{A_2}M_1
-\tilde{a_1}(\overline{C_1}\overline{A_2}+\overline{C_2}\overline{A_1})-\delta\tilde{a_2}(\overline{C_1}\overline{A_1}+\delta\overline{C_2}\overline{A_2})+\\
&{}&\overline{\tilde{c_1}}(C_1A_2+C_2A_1)+\delta\overline{\tilde{c_2}}(C_1A_1+\delta
C_2A_2)\Bigr)\Biggr),\\
N_4&=&-4\Bigl(R_1r_2+R_2r_1+2i(\overline{C_1}\overline{A_1}+\delta\overline{C_2}\overline{A_2})\hbox{Re}\,(A_1r_2^{*}+
A_2r_1^{*})\Bigr)+i\Biggl(\overline{C_1}\overline{\Delta_4}+\overline{C_2}\overline{\Delta_2}+\\
&{}&2\Bigl(A_1\overline{P_4}+
A_2\overline{P_2}-\overline{A_1}M_4-\overline{A_2}M_2
-\tilde{b_1}(\overline{C_1}\overline{A_2}+\overline{C_2}\overline{A_1})-\delta\tilde{b_2}(\overline{C_1}\overline{A_1}+\delta\overline{C_2}\overline{A_2})+\\
&{}&\overline{\tilde{d_1}}(C_1A_2+C_2A_1)+\delta\overline{\tilde{d_2}}(C_1A_1+\delta
C_2A_2)\Bigr)\Biggr),\\
\Theta_2-\delta\Theta_3&=&-16\hbox{Re}\,\Bigl(R_1(A_1r_1^{*}+\delta
A_2r_2^{*})-\delta
R_2(A_1r_2^{*}+A_2r_1^{*})\Bigr)-\\
&{}&4\hbox{Im}\,\Bigl(A_1(\overline{\Delta_2}-\delta\overline{\Delta_3})+
\delta A_2(\overline{\Delta_4}-\overline{\Delta_1})\Bigr),\\
\Theta_1-\Theta_4&=&16\hbox{Re}\,\Bigl(R_1(A_1r_2^{*}+
A_2r_1^{*})-
R_2(A_1r_1^{*}+\delta A_2r_2^{*})\Bigr)-\\
&{}&4\hbox{Im}\,\Bigl(A_1(\overline{\Delta_1}-\overline{\Delta_4})+
A_2(\delta\overline{\Delta_3}-\overline{\Delta_2})\Bigr)-\\
&{}&8\hbox{Im}\,\Biggl(\overline{A_1}\Bigl(A_1(\tilde{a_1^{*}}-\delta\tilde{b_2^{*}})+A_2(\tilde{b_1^{*}}-\tilde{a_2^{*}})+\overline{A_1}(\tilde{c_1^{*}}-\delta\tilde{d_2^{*}})\Bigr)+\\
&{}&\overline{A_2}\Bigl(A_1(\tilde{a_2^{*}}-\tilde{b_1^{*}})+A_2(\tilde{b_2^{*}}-\delta\tilde{a_1^{*}})+\overline{A_2}(\tilde{d_2^{*}}-\delta\tilde{c_1^{*}})\Bigr)\Biggr)
+4\hbox{Re}\,\Bigl(A_1G_1+A_2G_2\Bigr),\qquad(2.62)
\end{eqnarray*}
where $G_1,G_2$ are found from the following relations
\begin{eqnarray*}
\overline{C_1}G_1+\overline{C_2}G_2&=&v_3-\delta v_2-\delta
r_1\tilde{a_2}-r_2\tilde{a_1}-\delta\overline{r_1}\overline{\tilde{c_2}}-\overline{r_2}\overline{\tilde{c_1}}-8(R_1r_2-R_2r_1)+\\
&{}&2i\Bigl(-\tilde{a_1}(\overline{C_1}\overline{A_1}+\delta\overline{C_2}\overline{A_2})-\tilde{a_2}(\overline{C_1}\overline{A_2}+\overline{C_2}\overline{A_1})+\overline{\tilde{c_1}}(C_1A_1+\delta
C_2A_2)+\\
&{}&\overline{\tilde{c_2}}(C_1A_2+C_2A_1)
-\overline{A_1}(\overline{C_1}\tilde{a_1^{*}}+\overline{C_2}\tilde{b_1^{*}})-
\overline{A_2}(\overline{C_1}\tilde{a_2^{*}}+\overline{C_2}\tilde{b_2^{*}})-\\
&{}&\overline{C_2}(\delta\overline{A_2}\overline{\tilde{a_1^{*}}}+\overline{A_1}\overline{\tilde{b_1^{*}}}+\delta
A_2\overline{\tilde{c_1^{*}}}-A_2\overline{\tilde{d_2^{*}}})-
\overline{C_1}(\overline{A_2}\overline{\tilde{a_2^{*}}}+\delta\overline{A_1}\overline{\tilde{b_2^{*}}}-
A_1\overline{\tilde{c_1^{*}}}+\delta A_1\overline{\tilde{d_2^{*}}})
\Bigr)+\\
&{}&i\overline{C_1}(\overline{\Delta_4}-\overline{\Delta_1})+i\overline{C_2}(\overline{\Delta_2}-\delta\overline{\Delta_3}),\\
\overline{C_1}G_2+\delta\overline{C_2}G_1&=&v_1-v_4-\delta
r_1\tilde{b_2}-r_2\tilde{b_1}-\delta\overline{r_1}\overline{\tilde{d_2}}-\overline{r_2}\overline{\tilde{d_1}}-8(R_1r_1-\delta
R_2r_2)+\\
&{}&2i\Bigl(-\tilde{b_1}(\overline{C_1}\overline{A_1}+\delta\overline{C_2}\overline{A_2})-\tilde{b_2}(\overline{C_1}\overline{A_2}+\overline{C_2}\overline{A_1})
+\overline{\tilde{d_1}}(C_1A_1+\delta
C_2A_2)+\\
&{}&\overline{\tilde{d_2}}(C_1A_2+C_2A_1)
-\overline{A_1}(\delta\overline{C_2}\tilde{a_1^{*}}+\overline{C_1}\tilde{b_1^{*}})-
\delta\overline{A_2}(\overline{C_2}\tilde{a_2^{*}}+\delta\overline{C_1}\tilde{b_2^{*}})-\\
&{}&\overline{C_1}(\delta\overline{A_2}\overline{\tilde{a_1^{*}}}+\overline{A_1}\overline{\tilde{b_1^{*}}}+\delta
A_2\overline{\tilde{c_1^{*}}}-A_2\overline{\tilde{d_2^{*}}})-
\delta\overline{C_2}(\overline{A_2}\overline{\tilde{a_2^{*}}}+\delta\overline{A_1}\overline{\tilde{b_2^{*}}}-
A_1\overline{\tilde{c_1^{*}}}+\delta A_1\overline{\tilde{d_2^{*}}})
\Bigr)+\\
&{}&i\overline{C_1}(\overline{\Delta_2}-\delta\overline{\Delta_3})+i\delta\overline{C_2}(\overline{\Delta_4}-\overline{\Delta_1}).\qquad\qquad\qquad
(2.63)
\end{eqnarray*}

Further, plugging in (2.16) gives
\begin{eqnarray*}
\Gamma_1-\Gamma_4&=&2i\hbox{Re}\,\Bigl((C_1\overline{C_2}\overline{r_1}
+\delta|C_2|^2\overline{r_2})\rho_1+(C_1\overline{C_2}\overline{r_2}+|C_2|^2\overline{r_1})\rho_2\Bigr)+\\
&{}&4i\hbox{Im}\,\Bigl(A_2r_1^{*}-A_1r_2^{*}+C_1\overline{C_2}\tilde{a_2}\Bigr)+2C_1\overline{C_2}\hbox{Im}\,\Bigl(r_1\overline{r_2}\Bigr)+2|C_2|^2\Bigl(\tilde{b_2}+\delta\overline{\tilde{a_1}}-i|r_1|^2+i\delta|r_2|^2\Bigr),\\
\Gamma_2-\delta\Gamma_3&=&2i\hbox{Re}\,\Bigl(\delta(C_1\overline{C_2}\overline{r_2}
+|C_2|^2\overline{r_1})\rho_1+(C_1\overline{C_2}\overline{r_1}+\delta|C_2|^2\overline{r_2})\rho_2\Bigr)+\\
&{}&4i\hbox{Im}\,\Bigl(\delta
A_2r_2^{*}-A_1r_1^{*}+\delta C_1\overline{C_2}\tilde{a_1}\Bigr)+2\delta|C_2|^2\Bigl(\tilde{b_1}+\overline{\tilde{a_2}}\Bigr).\qquad\qquad\qquad
(2.64)
\end{eqnarray*}

Equations (2.18), (2.26) are automatically satisfied, so
we now have to use
(2.36), (2.54) and (2.58). The remaining part of the proof goes as
follows: first, in addition to (2.64), we
will obtain two more relations for $\Gamma_{\alpha}$ from (2.36) and
thus determine them; further, in addition to the relations for
$\Delta_{\alpha}$ in (2.61), conditions (2.54) will give two more
relations and thus fix $\Delta_{\alpha}$; finally, in
addition to the relations for $\Theta_{\alpha}$ in (2.62), (2.63), we
obtain two more relations from (2.58) which will determine
$\Theta_{\alpha}$. It is clear from (2.61), (2.62) that the choice of
$\Gamma_{\alpha}, \Delta_{\alpha}, \Theta_{\alpha}$ determines the rest of the
parameters as well. To get these extra relations (as well as for future applications) we need to find the expansions of
$\Phi^{\alpha}_{\beta}, \Phi^{\alpha},\Psi^{\alpha}$ completely
(cf. (2.34), (2.45)). The
calculations turn out to be so enormous that even writing down the
final formulas would be extremely lengthy; therefore we only give here
an outline of the procedure that allowed us to find these expansions. 

The coefficients in (2.35) can be found from (2.21) if one considers
terms cubic in $\omega^{\alpha},\overline{\omega^{\alpha}}$
there. Further, to determine terms containing $\theta^{\gamma}$ in the
expansions of $\Phi^{\alpha}_{\beta}$, we need to find the expressions
of $\lambda^{\alpha}\wedge\theta^{\beta}$ in (2.39). This is done analogously to the
proof of (2.46): we use (2.41), (2.44) to get information about the
expansions of $d\mu_1^1,d\mu_1^2$ in terms of $\lambda^{\alpha}$,
differentiate (2.33) and compare appropriate terms in both sides of
the resulting equations. In addition to this, to find terms containing
$\omega^{\gamma}\wedge\theta^{\delta}$ or $\theta^1\wedge\theta^2$ in the expansions of
$\Phi_{\beta}^{\alpha}$, we use (2.20). Thus, we determine the expansions
of $\Phi_{\beta}^{\alpha}$ completely, and at the same time the
expansions of $\Phi^{\alpha}$ (mod $\theta^{\beta}$), in particular, the
coefficients in (2.53); we also determine the imaginary parts of the
coefficients at $\omega^{\beta}\wedge\theta^{\gamma}$ and obtain
certain symmetries for the real parts of these coefficients. 
Further, (2.20) and the last unused terms in
the differentials of (2.33) give all terms in the expansions of
$\Phi^{\alpha}$ except for the real parts of the coefficients at
$\omega^{\beta}\wedge\theta^{\gamma}$ and for the terms containing
$\overline{\omega^{\beta}}\wedge\theta^{\gamma}$.  To find more terms in the expansions of $\Phi^{\alpha}$ as
well as some terms in the expansions of $\Psi^{\alpha}$ (in particular, terms in (2.57)) we use an analogue of
the above procedure where the differentiation of (2.7) (i.e. equations
(2.20)) is replaced by the differentiation of (2.12) and the
differentiation of (2.33) is replaced by the differentiation of
(2.38).

Eventually we find all the terms in the expansions of
$\Phi^{\alpha}_{\beta},\Phi^{\alpha},\Psi^{\alpha}$ except for the
coefficients at $\overline{\omega^{\gamma}}\wedge\theta^{\delta}$ in
$\Phi^{\alpha}$ and coefficients at
$\omega^{\gamma}\wedge\theta^{\beta},\overline{\omega^{\gamma}}\wedge\theta^{\beta}$
in $\Psi^{\alpha}$. All the coefficients that we have found are expressed in
terms of
$r_{\alpha},\rho_{\alpha},\tilde{a_{\alpha}},\tilde{b_{\alpha}},\tilde{c_{\alpha}},\tilde{d_{\alpha}}$.
An observation that is going to be important for future references is: if $r_1\equiv 0,r_2\equiv
0, \rho_1\equiv 0,\rho_2\equiv 0$, then $\tilde{a_{\alpha}}\equiv
0,\tilde{b_{\alpha}}\equiv 0, \tilde{c_{\alpha}}\equiv
0,\tilde{d_{\alpha}}\equiv 0$ and
\begin{eqnarray*}
\Omega_1^1&\equiv& 0,\\
\Omega_2^1&=&Q\overline{\underline{\omega}}\wedge\underline{\theta},\\
\Omega^0_2&=&\Bigl(P\underline{\omega}+\overline{P}\overline{\underline{\omega}}\Bigr)\wedge\underline{\theta},
\end{eqnarray*}
where $Q$ and $P$ are ${\frak A}$-valued functions on $P^2$.

Now that we know the expansions for $\Phi_{\beta}^{\alpha}$,
$\Phi^{\alpha}$, $\Psi^{\alpha}$, we can determine all the parameters
in (2.60). Namely, conditions (2.36) determine $\Gamma_1,\Gamma_3$,
conditions (2.54) determine $\Delta_1,\Delta_3$, conditions (2.58)
determine $\Theta_1,\Theta_3$. This determines the right-hand side in
(2.59) completely in terms of $r_{\alpha}$, $\rho_{\alpha}$,
$a_{\alpha}$, $b_{\alpha}$, $c_{\alpha}$, $d_{\alpha}$.

As one can see, most of the work in the above proof came from dealing
with $\phi_{\beta}^{\alpha},\mu_{\beta}^{\alpha},\psi^{\alpha}$. Since
the difference $L_{\eta}^{*}\omega-\hbox{Ad}(\eta)\omega$ does not
contain these forms, we may say that the form
$\omega$ on $P^2$ is reasonably close to being a Cartan connection.

The theorem is proved.\hfill $\Box$

\section{Corollaries and Applications}

In this section we derive some corollaries from Theorem 1.1 and the
construction of parallelism in Section 2.

First, we discuss the question whether the sequence of bundles
$P^2\rightarrow P^1\rightarrow M$ actually reduces to a single
principle bundle over $M$ with structure group $\hbox{Aut}_{0,e}(Q_H)$.

\begin{proposition} If $r_1\equiv 0,r_2\equiv 0$, then $P^2$ is a
principle $\hbox{Aut}_{0,e}(Q_H)$-bundle over $M$. If, in addition,
$\rho_1\equiv 0$, $\rho_2\equiv 0$, then the form $\omega$ is a Cartan
connection on $P^2$.
\end{proposition}

{\bf Proof.} Let $x\in P^2$ be
$x=(\underline{\tilde{\theta}}(y),\underline{\tilde{\omega}},\underline{\overline{\tilde{\omega}}},\underline{\tilde{\phi}})$
for $y=\underline{\theta}(p)\in P^1$, $p\in M$, where
$\underline{\tilde{\theta}}:=(\tilde{\theta^1},\tilde{\theta^2})$,
$\underline{\tilde{\omega}}:=(\tilde{\omega^1},\tilde{\omega^2})$,
$\underline{\tilde{\phi}}:=(\tilde{\phi^1},\tilde{\phi^2})$,
$\underline{\theta}(p):=(\theta^1(p),\theta^2(p))$,
$\tilde{\underline{\theta}}=\pi^{1*}\underline{\theta}$,
for
$\underline{\theta}\in{\frak M}_p$,
$\tilde{\underline{\omega}}=\pi^{1*}\underline{\omega}(y)$ for some complex
covectors $\underline{\omega}:=(\omega^1,\omega^2)$ at $p$, $\tilde{\phi^{\alpha}}$ are real
covectors at $y$, $\pi^1(y)=p$, $\pi^2(x)=y$. Let an element
$\eta\in\hbox{Aut}_{0,e}(Q_H)$ be represented by matrix (1.2).
We then define $F_{\eta}(x)$ as follows
$$
F_{\eta}(x):=\Bigl(\underline{\tilde{\theta}}(y'),C^{-1}\underline{\tilde{\omega^{*}}}-2A\underline{\tilde{\theta}}(y'),\overline{C^{-1}\underline{\tilde{\omega^{*}}}-2A\underline{\tilde{\theta}}(y')},\underline{\tilde{\phi^{*}}}-2iC^{-1}\overline{A}\underline{\tilde{\omega^{*}}}+2i\overline{C^{-1}}A\overline{\underline{\tilde{\omega^{*}}}}-4R\underline{\theta}(y')\Bigr),\eqno{(3.1)}
$$
where $y'=C^{-1}\overline{C^{-1}}\underline{\theta}(p)$ and
$\underline{\tilde{\omega^{*}}}, \underline{\tilde{\phi^{*}}}$ are the
pull-backs of the covectors
$\underline{\tilde{\omega}},\underline{\tilde{\phi}}$ respectively
under the diffeomorphism $\Phi_C$ of $P^1$ locally over a
neighbourhood of $p$ given by
$$
\Phi_C(D\underline{\theta}(q))=C\overline{C}D\underline{\theta}(q),
$$
for $D\in G^1$.

It is now easy to check that (3.1) indeed defines an action of
$\hbox{Aut}_{0,e}(Q_H)$ on $P^2$, provided $r_1\equiv 0,
r_2\equiv 0$.

Further, one can derive an analogue of transformation law (2.59) for
the form $\omega$ under the action of $\hbox{Aut}_{0,e}(Q_H)$ on
$P^2$ by the procedure described in Section 2. The transformation formula has
the same form (2.59), but the error terms in the right-hand side turn
out to be zero due to the identical vanishing of $r_{\alpha}$ and
$\rho_{\alpha}$, i.e. we get
$$
L_{\eta}^{*}\omega=\hbox{Ad}\,(\eta)\omega.
$$

The proposition is proved.\hfill $\Box$

\begin{remark} \rm If $r_1,r_2$ do not necessarily vanish, one can still
define for any $\eta\in \hbox{Aut}_{0,e}(Q_H)$ the mapping
$\tilde{F_{\eta}}$ as
\begin{eqnarray*}
\tilde{F_{\eta}}(x)&:=&\Biggl(\underline{\tilde{\theta}}(y'),C^{-1}\underline{\tilde{\omega^{*}}}-2A\underline{\tilde{\theta}}(y'),\overline{C^{-1}\underline{\tilde{\omega^{*}}}-2A\underline{\tilde{\theta}}(y')},\underline{\tilde{\phi^{*}}}-2iC^{-1}\overline{A}\underline{\tilde{\omega^{*}}}+2i\overline{C^{-1}}A\overline{\underline{\tilde{\omega^{*}}}}-4R\underline{\theta}(y')+\\
&{}&tC\overline{C}\left(
\begin{array}{ccc}
\displaystyle r_1(y')&\delta r_2(y')\\
\displaystyle r_2(y')&r_1(y')
\end{array}
\right)
\tilde{\underline{\omega^{*}}}
\Biggr),
\end{eqnarray*}
where $t$ is determined by
$$
C^{-1}\overline{C^{-1}}=
\left(
\begin{array}{ccc}
\displaystyle s&\delta t\\
\displaystyle t&s
\end{array}
\right).
$$
However, $\tilde{F_{\eta}}$ does not give a group action
unless $r_{\alpha}$ identically vanish.
\end{remark}

Next, we characterize $CR$-flat manifolds, i.e. manifolds
for which the form $\Omega$ in (2.32) vanishes.

\begin{proposition} The form $\Omega$ identically vanishes on $P^2$ if
and only if $M$ is locally $CR$-equivalent to $Q_H$. 
\end{proposition}

{\bf Proof.} First, we will explicitly calculate the bundles $P^1$,
$P^2$ and the forms $\omega$, $\Omega$ for $Q_H$. To do this,
we identify $(z_1,z_2)$ and $(w_1,w_2)$ with the matrices
$$
Z:=\left(
\begin{array}{ccc}
\displaystyle z_1&\delta z_2\\
\displaystyle z_2&z_1
\end{array}
\right)\qquad\hbox{and}\qquad
W:=\left(
\begin{array}{ccc}
\displaystyle w_1&\delta w_2\\
\displaystyle w_2&w_1
\end{array}
\right)
$$
respectively and write the equation of $Q_H$ in the form
$$
V=Z\overline{Z},\eqno{(3.2)}
$$
with $W=U+iV$. Let the orientation of $Q_H$ be given by the form
$\underline{\theta}:=\frac{1}{2}(dU-i\overline{Z}dZ+iZd\overline{Z})$.
For
$\underline{\omega}:=dZ$ we have
$$
d\underline{\theta}=i\underline{\omega}\wedge\overline{\underline{\omega}}.
$$
The bundle $P^1$ then consists of $D\underline{\theta}$ with
$D=C\overline{C}$ for $C\in{\frak A}^{*}$. The tautological form
$\tilde{\underline{\theta}}$ is given by
$$
\tilde{\underline{\theta}}=\frac{1}{2}D\left(dU-i\overline{Z}dZ+iZd\overline{Z}\right),
$$
and therefore $P^2$ consists of coframes of the form
$$
x:=\Bigl(\tilde{\underline{\theta}},T\tilde{\underline{\theta}}+C\tilde{\underline{\omega}},\overline{T\tilde{\underline{\theta}}+C\tilde{\underline{\omega}}},S\tilde{\underline{\theta}}+iC\overline{T}\tilde{\underline{\omega}}-i\overline{C}T\overline{\tilde{\underline{\omega}}}-D^{-1}dD\Bigr),
$$
with $T,C,S$ as in (1.4), $C\overline{C}=E$, where
$\tilde{\underline{\omega}}:=\sqrt{D}dZ$ (it is an easy exercise to prove
the existence of a locally smooth operation of taking real square root on
$G^1$).

Now one can check that
\begin{eqnarray*}
\hat{\underline{\theta}}&=&\frac{1}{2}D\left(dU-i\overline{Z}dZ+iZd\overline{Z}\right),\\
\hat{\underline{\omega}}&=&T\hat{\underline{\theta}}+C\sqrt{D}dZ,\\
\hat{\underline{\phi}}&=&S\hat{\underline{\theta}}+iC\overline{T}\sqrt{D}dZ-i\overline{C}T\sqrt{D}d\overline{Z}-D^{-1}dD,\\
\underline{\underline{\phi}}&=&\frac{1}{2}(S-3iT\overline{T})\hat{\underline{\theta}}+2i\overline{T}\hat{\underline{\omega}}+iT\overline{\hat{\underline{\omega}}}-\overline{C}dC-\frac{1}{2}D^{-1}dD,\\
\underline{\mu}&=&iT^2\overline{T}\hat{\underline{\theta}}+\frac{1}{2}(S-iT\overline{T})\hat{\underline{\omega}}-iT^2\overline{\hat{\underline{\omega}}}-dT+\overline{C}TdC-\frac{1}{2}TD^{-1}dD,\\
\underline{\psi}&=&\frac{1}{2}(S^2-3T^2\overline{T^2})\hat{\underline{\theta}}+(iS\overline{T}+T\overline{T^2})\hat{\underline{\omega}}+(-iST+T^2\overline{T})\overline{\hat{\underline{\omega}}}-dS-\\
&{}&i\overline{T}dT+iTd\overline{T}+
2i\overline{C}T\overline{T}dC-SD^{-1}dD,\qquad\qquad\qquad\qquad
(3.3)
\end{eqnarray*}
and all the functions $r_{\alpha},\rho_{\alpha},
\tilde{a_{\alpha}},\tilde{b_{\alpha}},\tilde{c_{\alpha}},\tilde{d_{\alpha}}$
identically vanish. It now follows from (2.32) and (3.3) that for $\omega$
defined as in (2.31) its curvature form $\Omega\equiv 0$.  

Now, if $M$ is a $CR$-manifold with $\Omega\equiv 0$, then $P^2$ locally can be
mapped by a diffeomorphism onto a neighbourhood of identity in $\hbox{Aut}_{e}(Q_H)$ in such a way that $\omega$ transforms
into the Maurer-Cartan form on $\hbox{Aut}_{e}(Q_H)$ (see \cite
{St}). Therefore, there exists a local diffeomorphism of $P^2$ to the
corresponding bundle over $Q_H$ that preserves
parallelism. By
Theorem 1.1 this diffeomorphism is the lift of a local
$CR$-diffeomorphism between $M$ and $Q_H$, and thus $M$ is
locally $CR$-equivalent to $Q_H$.

The proposition is proved. \hfill $\Box$

We will now try to understand what
proper analogues of chains in the case of hyperbolic and
elliptic $CR$-manifolds should be. We define the {\it chain
distribution} on $P^2$ by
$$
\underline{\omega}=0,\qquad
\underline{\mu}=0.\eqno{(3.4)}
$$
This distribution is analogous to the one of Chern
\cite{CM} that, in the case of $CR$-codimension 1, was used to define chains on $M$ as the projections of
the integral manifolds of the distribution. However, in contrast with
\cite{CM}, distribution (3.4) may not be
integrable. It follows from the expansions of $\Phi_{\beta}^{\alpha}$,
$\Phi^{\alpha}$, $\Psi^{\alpha}$ that can be found as outlined in Section
2, that it is integrable if and only if the following conditions are satisfied
(cf. Proposition 3.1):
\begin{eqnarray*}
r_1&\equiv& r_2\equiv 0,\\
\rho_1&\equiv& \rho_2\equiv 0,\qquad\qquad\qquad\qquad (3.5)
\end{eqnarray*}
(note that the functions $r_{\alpha},\rho_{\alpha}$ are scalar
$CR$-invariants). For manifolds satisfying conditions (3.5)
one can project the integral manifolds of chain distribution (3.4) to $M$. The resulting two-dimensional submanifolds of $M$ we
call {\it $G$-chains}.

It is clear from (2.23), (2.24) that conditions (3.5) are equivalent to
$$
r\equiv 0,\qquad \rho\equiv 0,\eqno{(3.6)}
$$
where $r:=r_1r_2$, $\rho:=\rho_1\rho_2$. Elliptic or hyperbolic
$CR$-manifolds satisfying (3.6) we call {\it semi-flat}. It follows
from (2.6), (2.7), (2.32) that semi-flat manifolds are characterized
by the curvature conditions
$$
\Omega_0^1\equiv 0,\qquad \Omega_0^2\equiv 0.\eqno{(3.7)}
$$
Note that conditions (3.7) are always satisfied for the parallelism
constructed in \cite{CM}. For semi-flat manifolds the main formulas in
Section 2 reduce to matrix forms of Chern's formulas in the case of
$CR$-dimension and $CR$-codimension 1. Namely, we
have
\begin{eqnarray*}
d\underline{\theta}&=&i\underline{\omega}\wedge\overline{\underline{\omega}}+\underline{\theta}\wedge\underline{\phi},\\
d\underline{\omega}&=&\underline{\omega}\wedge\underline{\underline{\phi}}+\underline{\theta}\wedge\underline{\mu},\\
\hbox{Re}\,\underline{\underline{\phi}}&=&\frac{1}{2}\underline{\phi},\\
d\underline{\phi}&=&2\hbox{Re}\,\Bigl(i\underline{\mu}\wedge\overline{\underline{\omega}}\Bigr)+\underline{\theta}\wedge\underline{\psi},\\
\Omega_1^1&\equiv& 0,\\
\Omega_2^1&=&Q\overline{\underline{\omega}}\wedge\underline{\theta},\\
\Omega_2^0&=&\Bigl(P\underline{\omega}+\overline{P}\overline{\underline{\omega}}\Bigr)\wedge\underline{\theta},
\end{eqnarray*}
where $Q$ and $P$ are ${\frak A}$-valued functions on $P^2$, and
\begin{eqnarray*}
\mu_2^2&=&\mu_1^1,\\
\mu_2^1&=&\delta\mu_1^2,\\
\psi^2&=&\delta\psi^3,\\
\psi^4&=&\psi^1,\\
\psi^5&=&\psi^1,\\
\psi^6&=&\delta\psi^3,\\
\psi^7&=&\psi^3,\\
\psi^8&=&\psi^1.
\end{eqnarray*}

It follows from the proof of Proposition 3.3 that the quadric $Q_H$ is a semi-flat
manifold. In the next proposition we describe $G$-chains on $Q_H$.

\begin{proposition} Any $G$-chain on $Q_H$ passing through the origin
is the intersection of $Q_H$ with
$Z=AW$ for some $A\in{\frak A}$.
\end{proposition}

{\bf Proof.} Formulas (3.3) imply that distribution (3.4) in the case
of $Q_H$ is given by
\begin{eqnarray*}
T\hat{\underline{\theta}}+C\sqrt{D}dZ&=&0,\qquad\qquad\qquad\qquad (3.8.a)\\
\frac{1}{2}\Bigl(TS+iT^2\overline{T}\Bigr)\hat{\underline{\theta}}+dT&=&0.\qquad\qquad\qquad\qquad
(3.8.b)
\end{eqnarray*}

First we show that along the integral manifolds of (3.8)
passing through points of the fibre of $P^2$ over the origin,
$G:=C-iT\sqrt{D}\overline{Z}$ is non-degenerate. To do this, we
differentiate
$G\overline{G}$ and
plug in the resulting expression $dZ$ and $dT$ found from (3.8). It
then easily follows that
$d\Bigl(G\overline{G}\Bigr)\equiv 0$ and thus
$\hbox{det}\,|G|^2\equiv\hbox{const}$. Since $\hbox{det}\,G\ne 0$ for
$Z=0$, $G$ is non-degenerate everywhere.

Next, it follows from (3.8.a) that
$$
dZ=-\frac{1}{2}TDG^{-1}\Bigl(dU+i\overline{Z}dZ+iZd\overline{Z}\Bigr).
$$
Therefore, to show that the projections of the integral manifolds of
distribution (3.8) passing through points of the fibre of $P^2$ over
the origin have the desired form, we need only prove that
$$
TDG^{-1}=\hbox{const}\eqno{(3.9)}
$$
along the integral manifolds of (3.8). To prove (3.9) we
differentiate $TDG^{-1}$ and by using (3.8) conclude that
$d(TDG^{-1})\equiv 0$.

The proposition is proved. \hfill $\Box$

As we have seen, semi-flat manifolds possess some of the nice properties that
were observed in \cite{CM} for manifolds of $CR$-codimension 1. The quadric $Q_H$ is
an example of such manifolds. Many more examples come from considering
{\it matrix surfaces} in $\CC^4$, i.e. real-analytic surfaces locally
near the origin given in the form
$$
V=Z\overline{Z}+\sum_{k+l+m\ge 2}A_{k,\overline{l},m}Z^k\overline{Z^l}U^m,\eqno{(3.10)}
$$
where $A_{k,\overline{l},m}\in{\frak A}$, $A_{k,\overline{l},m}=\overline{A_{l,\overline{k},m}}$,
$A_{1,\overline{1},0}=0$, and the power series in the right-hand side converges
in a neighbourhood of the origin. The quadric $Q_H$ written in the
form (3.2) is a matrix surface.

Hyperbolic ($\delta=1$) matrix surfaces are easily described. Indeed,
the mapping
\begin{eqnarray*}
z_1^{*}&=&z_1+z_2,\\
z_2^{*}&=&z_1-z_2,\\
w_1^{*}&=&w_1+w_2,\\
w_2^{*}&=&w_1-w_2,\qquad\qquad\qquad (3.11)
\end{eqnarray*}
maps a hyperbolic matrix surface into a direct product of real
hypersurfaces in $\CC^2$. In the next proposition we show that
semi-flat hyperbolic manifolds can be characterized in a similar way.

\begin{proposition} A semi-flat hyperbolic manifold is locally $CR$-equivalent
to a product of 3-dimensional Levi non-degenerate $CR$-manifolds of codimension 1.
\end{proposition}

{\bf Proof.} It is not hard to show that a manifold $M$ is semi-flat
if and only if near every point $p\in M$ there exist
complex 1-forms $\omega^1,\omega^2$ that at every point $q$ are complex
linear on $T_q^c(M)$, real 1-forms
$\theta^1,\theta^2$ whose common annihilator at every point is $T_q^c(M)$, real 1-forms $\phi^1,\phi^2$ and complex 1-forms $\lambda^1,\lambda^2,\mu^1,\mu^2$
such that near $p$ the following holds
\begin{eqnarray*}
d\theta^1&=&i\left(\omega^1\wedge\overline{\omega^1}+\delta\omega^2\wedge\overline{\omega^2}\right)+\theta^1\wedge\phi^1+\delta\theta^2\wedge\phi^2,\\
d\theta^1&=&i\left(\omega^1\wedge\overline{\omega^2}+\omega^2\wedge\overline{\omega^1}\right)+\theta^1\wedge\phi^2+\theta^2\wedge\phi^1,\\
d\omega^1&=&\omega^1\wedge\lambda^1+\delta\omega^2\wedge\lambda^2+\theta^1\wedge\mu^1+\delta\theta^2\wedge\mu^2,\\
d\omega^2&=&\omega^1\wedge\lambda^2+\omega^2\wedge\lambda^1+\theta^1\wedge\mu^2+\theta^2\wedge\mu^1,
\end{eqnarray*}
and
$(\theta^{\alpha},\hbox{Re}\,\omega^{\alpha},\hbox{Im}\,\omega^{\alpha})$
at every point form a coframe. Then, in the case of hyperbolic manifolds, the forms
\begin{eqnarray*}
\theta^{1'}&:=&\theta^1+\theta^2,\\
\theta^{2'}&:=&\theta^1-\theta^2,\\
\omega^{1'}&:=&\omega^1+\omega^2,\\
\omega^{2'}&:=&\omega^1-\omega^2,\\
\phi^{1'}&:=&\phi^1+\phi^2,\\
\phi^{2'}&:=&\phi^1-\phi^2,\\
\lambda^{1'}&:=&\lambda^1+\lambda^2,\\
\lambda^{2'}&:=&\lambda^1-\lambda^2,\\
\mu^{1'}&:=&\mu^1+\mu^2,\\
\mu^{2'}&:=&\mu^1-\mu^2
\end{eqnarray*}
satisfy
\begin{eqnarray*}
d\theta^{1'}&=&i\omega^{1'}\wedge\overline{\omega^{1'}}+\theta^{1'}\wedge\phi^{1'},\\
d\theta^{2'}&=&i\omega^{2'}\wedge\overline{\omega^{2'}}+\theta^{2'}\wedge\phi^{2'},\\
d\omega^{1'}&=&\omega^{1'}\wedge\lambda^{1'}+\theta^{1'}\wedge\mu^{1'},\\
d\omega^{2'}&=&\omega^{2'}\wedge\lambda^{2'}+\theta^{2'}\wedge\mu^{2'}.\qquad\qquad\qquad
(3.12)
\end{eqnarray*}
Formulas (3.12) now imply that the distribution
$$
\theta^{\alpha'}=0,\qquad\omega^{\alpha'}=0,
$$
for each $\alpha=1,2$ is integrable and thus gives a foliation of $M$
near $p$ by 3-dimensional Levi non-degenerate $CR$-manifolds of
$CR$-dimension 1. Let $M^1,M^2$ be the leaves of the first and the
second foliation respectively that pass through $p$. Then $M$ is
clearly $CR$-equivalent near $p$ to the product $M^1\times M^2$ (see also Proposition 5.8 in \cite{Ch}).

The proposition is proved. \hfill $\Box$

\begin{remark} \rm If in the above proposition a semi-flat hyperbolic manifold $M$ is in
addition real-analytic, then $M^1,M^2$ are also real-analytic and
therefore admit real-analytic $CR$-embeddings in $\CC^3$ as
hypersurfaces \cite{AH}. Mapping the image of the point $p$ into the origin and
applying the transformation inverse to (3.11) we see that $M$ near $p$
is $CR$-equivalent to a surface of the form (3.10), and therefore real-analytic
semi-flat hyperbolic manifolds are characterized locally as matrix
surfaces.
\end{remark}

\section {Normal Forms}

In this section we consider real-analytic hyperbolic and elliptic
$CR$-manifolds that by \cite{AH} can be assumed to be locally embedded in $\CC^4$
near the origin. We denote coordinates in $\CC^4$ by 
$z:=(z_1,z_2)$,
$w:=(w_1,w_2)$, $u:=(u_1,u_2):=(\hbox{Re}\,w_1,\hbox{Re}\,w_2)$,
$v:=(v_1,v_2):=(\hbox{Im}\,w_1,\hbox{Im}\,w_2)$. Let $M$ be such a manifold
and suppose that the coordinates are chosen
so that $T_0(M)$ is spanned by $z,u$ and $T_0^c(M)$ by $z$. Then $M$ is
given by an equation of the form
$$
v=H^{\delta}(z,z)+F(z,\overline{z},u),
$$
where $F$ is an $\RR^2$-valued real-analytic function such that
$F(0)=0$, $dF(0)=0$, $\frac{\partial^2 F}{\partial z_i\partial\overline{z_j}}(0)=0$. 

Important
examples of hyperbolic and elliptic manifolds are matrix surfaces
(3.10). As we noted in Section 3, transformation (3.11) maps a
hyperbolic matrix surface into a direct product:
\begin{eqnarray*}
v_1&=&|z_1|^2+F^1(z_1,\overline{z_1},u_1),\\
v_2&=&|z_2|^2+F^2(z_2,\overline{z_2},u_2),\qquad\qquad\qquad\qquad
(4.1)
\end{eqnarray*}
where $F^k$ are real-analytic functions, $F^k(0)=0$, $dF^k(0)=0$,
$\frac
{\partial^2 F^k}{\partial z_k\partial\overline{z_k}}(0)=0$.  In the elliptic 
case, the transformation
\begin{eqnarray*}
z_1^{*}&=&z_1+iz_2,\\
z_2^{*}&=&z_1-iz_2,\\
w_1^{*}&=&w_1,\\
w_2^{*}&=&w_2\qquad\qquad\qquad\qquad (4.2)
\end{eqnarray*}
maps the surface into a surface of the form
$$
{\cal V}=z_1\overline{z_2}+F(z_1,\overline{z_2},{\cal U}),\eqno{(4.3)}
$$
where ${\cal U}:=u_1+iu_2$, ${\cal V}:=v_1+iv_2$ and $F$ is a
$\CC$-valued analytic function, $F(0)=0$, $dF(0)=0$, $\frac{\partial^2 F}{\partial z_1\partial\overline{z_2}}(0)=0$. For convenience, we will use the forms
(4.1), (4.3) for matrix surfaces instead of (3.10).

Equations of hyperbolic and elliptic manifolds can be written in
normal forms \cite{Lo1}, \cite{ES2} that may be viewed as
generalizations of the Chern-Moser normal forms for real-analytic Levi
non-degenerate hypersurfaces in $\CC^n$. Here we write them in a
modified way as follows (the center of normalization is assumed to be
at the origin).
\medskip\\

\noindent The hyperbolic normal form:
\begin{eqnarray*}
v_1-|z_1|^2&=&N^1(z,\overline{z},u):=N_1^1(z_1,\overline{z_1},u_1)+N_2^1(z,\overline{z},u),\\
v_2-|z_2|^2&=&N^2(z,\overline{z},u):=N_2^2(z_2,\overline{z_2},u_2)+N_1^2(z,\overline{z},u),\qquad\qquad\qquad\qquad
(4.4)
\end{eqnarray*}
where $N_1^1$ and $N_2^2$ are in the Chern-Moser
normal form, i.e.
$$
N_j^j=2\hbox{Re}\,(h^j_{4,\overline{2}}(u_j)z_j^4\overline{z_j}^2)+\sum_{
\begin{array}{lll}
k,l\ge 2,\\
k+l\ge 7
\end{array}
}h^j_{k,\overline{l}}(u_j)z_j^k\overline{z_j}^l,
$$ 
each
monomial in $N_2^1$ contains at least one of the variables
$z_2,\overline{z_2},u_2$ and satisfies the following conditions
\begin{eqnarray*}
N^1_{2\,\,\,k,\overline{0}}&\equiv& 0,\qquad k\ge 1,\\
N^1_{2\,\,\,1,\overline{1}}&\equiv& 0,\\
\frac{\partial N^1_{2\,\,\,k,\overline{1}}}{\partial
\overline{z_1}}&\equiv& 0,\qquad k\ge 2,\\
\frac{\partial^2 N^1_{2\,\,\,2,\overline{1}}}{\partial
z_1\partial z_2}&\equiv& 0,\\
\frac{\partial^4 N^1_{2\,\,\,2,\overline{2}}}{\partial
z_1\partial z_2\partial\overline{z_1}\partial\overline{z_2}}&\equiv&
0,
\end{eqnarray*}
and $N_1^2$ has the same properties as $N_2^1$ above with interchanged
indices 1,2.
\smallskip\\

\noindent The elliptic normal form:
$$
{\cal V}-z_1\overline{z_2}=N(z,\overline{z},u):=N_1(z_1,\overline{z_2},{\cal
U})+N_2(z,\overline{z},{\cal U},\overline{{\cal U}}),\eqno{(4.5)}
$$
where the real and imaginary parts of
$N_1(\zeta,\overline{\zeta},\tau)$, $\zeta\in\CC$, $\tau\in\RR$, are
in the Chern-Moser normal form, each monomial in $N_2$
contains at least one of the variables
$\overline{z_1},z_2,\overline{{\cal U}}$ and satisfies the 
conditions
\begin{eqnarray*}
N_{2\,\,\,k,\overline{0}}&\equiv& 0,\qquad k\ge 1,\\
N_{2\,\,\,1,\overline{1}}&\equiv& 0,\\
\frac{\partial N_{2\,\,\,k,\overline{1}}}{\partial
\overline{z_2}}&\equiv& 0,\qquad k\ge 2,\\
\frac{\partial^2 N_{2\,\,\,2,\overline{1}}}{\partial z_1\partial
z_2}&\equiv& 0,\\
\frac{\partial^4 N_{2\,\,\,2,\overline{2}}}{\partial z_1\partial
z_2\partial\overline{z_1}\partial\overline{z_2}}&\equiv& 0.
\end{eqnarray*}
\medskip\\

We emphasize that in representations (4.4), (4.5) all other conditions of
the normal forms from \cite{Lo1}, \cite{ES2} are satisfied automatically.

We call 
$
\left(
\begin{array}{ccc}
\displaystyle N_1^1\\
\displaystyle N_2^2
\end{array}
\right)
$ 
and $N_1$ in (4.4), (4.5) respectively the {\it
matrix non-quadratic part} of the normal form;
$
\left(
\begin{array}{ccc}
\displaystyle N_2^1\\
\displaystyle N_1^2
\end{array}
\right)
$
and $N_2$ 
the {\it non-matrix part}. We say that a normal form of a hyperbolic
(resp. elliptic) surface $M$ is a {\it matrix normal form}
if $N_2^1\equiv N_1^2\equiv 0$
(resp. $N_2\equiv 0$).

\begin{proposition} Let $M$ be a matrix surface given by a power
series of the form (3.10) that converges in a neighbourhood $\Omega$ of the
origin. Then any normal form of $M$ with center at any point of $\Omega$ is
a matrix normal form.
\end{proposition}

{\bf Proof.} We recall \cite{CM} that for a real-analytic Levi non-degenerate
hypersurface in $\CC^2$ given by
$$
v=|z|^2+G(z,\overline{z},u),
$$
with $G(0)=dG(0)=\frac{\partial^2 G}{\partial
z\partial\overline{z}}(0)=0$, a normalizing mapping with initial data
$(c,a,r)$ (see \cite{V} for the definition)
can be obtained in the
form
$$
z\mapsto \frac{c(z+aw+\phi(z,w))}{1-2i\overline{a}z-(r+i|a|^2)w},\qquad
w\mapsto \frac{|c|^2(w+\psi(z,w))}{1-2i\overline{a}{z}-(r+i|a|^2)w},\qquad d\phi(0)=d\psi(0)=0,
$$
where the holomorphic functions $\phi,\psi$ are uniquely determined by
$c,a,r$ and $G$. Analogously, for any set of initial data $(C,A,R)$
\cite{Lo1}, \cite{ES2} we can find a mapping of the form
\begin{eqnarray*}
Z&\mapsto&
C\Biggl(Z+AW+\Phi(Z,W)\Biggr)\Biggl(E-2i\overline{A}Z-(R+iA\overline{A})W\Biggr)^{-1}, \qquad
d\Phi(0)=0,\\
W&\mapsto&
C\overline{C}\Biggl(W+\Psi(Z,W)\Biggr)\Biggl(E-2i\overline{A}Z-(R+iA\overline{A})W\Biggr)^{-1},\qquad
d\Psi(0)=0,
\end{eqnarray*}
that transforms $M$ into a surface given in the form (3.10) with
$A_{k,\overline{0},m}=0$, $A_{k,\overline{1},m}=0$,
$A_{2,\overline{2},m}=0$, $A_{2,\overline{3},m}=0$,
$A_{3,\overline{3},m}=0$ for all $k,m$. By transformations
(3.11), (4.2) such equations are mapped into equations in the matrix normal
form.  

It now follows from the uniqueness of normalization with prescribed
initial data (see \cite{Lo2},
\cite{ES2}) that 
any normal form of $M$ with center at the origin is a matrix normal
form. Since a transformation that maps a fixed point of $\Omega$ into the
origin can be chosen in the matrix form,  
any normal form of $M$ with
center at any point of $\Omega$ is also matrix.

The proposition is proved. \hfill $\Box$

Suppose $M$ is given in normal form (4.4) or (4.5). We write the
equation in a standard way as a sum of weighted homogeneous
polynomials \cite{Lo1}, \cite{ES2}. Let $\nu$ be the smallest weight
for which there exists a nonvanishing polynomial in the non-matrix
part of the equation. We define $\kappa(0):=\frac{1}{\nu}$. Clearly,
$\kappa(0)=0$ if and only if $M$ is in a matrix normal form. It
follows from Proposition 4.1 that, for a matrix surface, $\kappa(0)=0$
for any normalization.

\begin{proposition} The number $\kappa(0)$ does not depend on the choice of
normalization for any $M$.
\end{proposition}

{\bf Proof.} Without loss of generality we can assume that the
differential of the renormalizing transformation at the origin is
identical on $T_0^c(M)$. We follow the scheme of normalization in \cite{Lo1},
\cite{ES2} and  split the transformation into a matrix
and non-matrix parts. The non-matrix part begins with a term of weight
$\ge\frac{1}{\kappa(0)}$, and the lowest
weight of the new non-matrix part of the equation in normal form
remains $\nu$.   

The proposition is proved. \hfill $\Box$

Moving the center of normalization, we define a function $\kappa$ on
all of $M$. By Proposition 4.2, $\kappa$ is a holomorphic
invariant. This invariant measures  the ``non-matricity'' of the surface and
is analogous to the $CR$-invariant
functions $r,\rho$ on $P^2$ from (3.6) that measure the
``non-semi-flatness'' of the manifold.

We recall that a 2-dimensional real-analytic submanifold
$\Gamma\subset M$ through the origin is called a {\it chain}, if in
some normal coordinates it is defined by $\{z=0,v=0\}$ \cite{Lo1},
\cite{ES2}. For example, such chains on quadrics (3.2) are the
intersections of the quadrics with complex planes $Z=AW$ for
$A\in{\frak A}$ and thus coincide with $G$-chains by Proposition
3.4. Chains form a holomorphically invariant family and are defined by
a mixed system of ODE's and PDE's. The main disadvantage of
2-dimensional chains as opposed the one-dimensional Chern-Moser chains,
is that in general translations along 2-dimensional chains spoil the normal form
conditions. However, matrix manifolds are a class of manifolds with
proper chains: chains on matrix manifolds are given by a matrix
generalization of the Chern-Moser chains. More precisely, any chain on
matrix surface (3.10) has the parametric form
$Z=P(U),W=Q(U)$ with $P,Q$ satisfying certain equations
\begin{eqnarray*}
D^2P&=&\Phi'(DP,P,Q,U),\\
D^3Q&=&\Psi'(D^2Q,DQ,Q,U),
\end{eqnarray*}
where the operator
$$
D:=
\left(
\begin{array}{ccc}
\displaystyle{\frac{\partial}{\partial u_1}}&\displaystyle{\delta\frac{\partial}{\partial u_2}}\\
\displaystyle{\frac{\partial}{\partial u_2}}&\displaystyle{\frac{\partial}{\partial
u_1}}\\
\end{array}
\right)
$$
is a matrix analogue of $\frac{\partial}{\partial u}$.

In Section 3 we defined $G$-chains on any semi-flat manifold, in
particular, on any matrix surface. It follows from a matrix analogue
of calculations on pp. 265--268 in \cite{CM} that, in the case of
matrix surfaces, $G$-chains coincide with chains (note that Proposition 3.4 gives an independent proof of this fact
for the quadrics).  

\section {Final Remarks and Questions}

We conclude the paper with a list of questions related to this work
that we believe are important for further study of the subject.
\medskip\\

\noindent {\bf 1.} Let $H$ be a non-degenerate $\RR^k$-valued Hermitian form
on $\CC^n$ and $Q_H$ the quadric in $\CC^{n+k}$ associated with
$H$. It is
well-known that the algebra  ${\frak g}_H$ of infinitesimal
automorphisms of $Q_H$ is a graded Lie algebra: ${\frak
g}_H=\oplus_{k=-2}^{2}{\frak g}_H^k$ \cite{Sa}, \cite{B1}. Is it true
that ${\frak g}_H$ is isomorphic to the maximal prolongation
$\tilde{{\frak g}}_H$ of $\oplus_{k=-2}^{0}{\frak g}_H^k$ in the sense
of Tanaka (see \cite{Ta1})? So far, we have not been able to find any
references on this matter except for the case $k=1$
\cite{Ta1} and the situation considered in
\cite{La}, and we have produced proofs that give positive
answers to this question for each
of $H^1$, $H^{-1}$, $H^0$ (here $k=2,n=2$). Note that ${\frak g}_H$ is
always isomorphic to a subalgebra of $\tilde{{\frak g}}_H$ by a
mapping that preserves grading.\footnote{After this work had been
completed, the first two authors proved that ${\frak g}_H$ and
$\tilde{\frak g}_H$ are always isomorphic thus giving a positive
answer to question 1. The proof will appear
elsewhere. In the meantime see \cite{EI}.}
\medskip\\

\noindent {\bf 2.} As we noted in the Introduction, every known result on the
equivalence problem for $CR$-manifolds falls in one of the two types: strongly uniform
manifolds or
weakly uniform manifolds with certain generic Levi forms. 
It is a reasonable question whether these two
groups of results treat in fact non-intersecting collections of manifolds. 
Namely, let $M$ be a weakly uniform connected $CR$-manifold, $p\in M$,
and the Levi form at $p$ is {\it not} in general position as in
\cite{ES6}; is then $M$ strongly uniform?

One can ask a stronger question as follows. Let $Q_{H^1}$,
$Q_{H^2}$ be two irreducible (i.e. not equivalent to direct products) quadrics and $\hbox{Aut}_{lin,e}(Q_{H^1})$ is
isomorphic to $\hbox{Aut}_{lin,e}(Q_{H^2})$ in such a way that the
isomorphism extends to an isomorphism between
$\hbox{Aut}_{0,e}(Q_{H^1})$ and $\hbox{Aut}_{0,e}(Q_{H^2})$. Suppose
that at least one of $Q_{H^1}$, $Q_{H^2}$ is not in general position
as in \cite{ES6}. Is it then true that $H^1$ is equivalent to $H^2$? 
\medskip\\

{\bf 3.} In Proposition 3.5 we characterized semi-flat hyperbolic
manifolds. In particular, in the real-analytic case they turned out to
be locally $CR$-equivalent to matrix surfaces. It would be
interesting to obtain analogues of these facts for elliptic
manifolds. In particular, is it true that a semi-flat elliptic
real-analytic manifold is locally equivalent to a matrix surface?

So far, we have been able to obtain only the following
result that we mention below without a proof (cf. the proof of Proposition 3.5).  

\begin{proposition} Let $M$ be a semi-flat elliptic manifold. Then
there are two foliations $\Sigma_1, \Sigma_2$ of $M$ by complex curves
such that for every point $p\in M$, the complex tangent space to $M$
at $p$ is spanned by the tangent spaces to the leaves of
$\Sigma_1,\Sigma_2$ at $p$.
\end{proposition}

{\bf 4.} The parallelism $\omega$ that we constructed in Theorem 1.1 is not in general a
Cartan connection. Is it possible to find a parallelism that is at the
same time a
Cartan connection for hyperbolic and elliptic $CR$-manifolds and for
general strongly uniform manifolds (cf. \cite{Ta1})?

It can be shown that the parallelism $\omega$ turns into a Cartan connection on $P^2$ with
respect to the action of $G^2$ if and only if
the manifold is semi-flat (cf. Proposition 3.1). Is it true that the
existence of a $CR$-invariant
Cartan connection on some other fibre bundle over a hyperbolic or
elliptic $CR$-manifold $M$ implies the semi-flatness of $M$?

One can introduce
matrix manifolds whenever $CR\hbox{dim}(M)$=$CR\hbox{codim}(M)$
by using matrix algebras other than ${\frak A}$ \cite{ES7}; can one
construct Cartan connections for such manifolds? 

In \cite{M}, \cite{GM} $CR$-invariant connections ({\it not} Cartan
connections) were constructed for certain weakly uniform
$CR$-structures. In these cases the groups 
$\hbox{Aut}_{0,e}(Q_{{\cal L}(M)(p)})$ contain only linear
automorphisms given by certain diagonal matrices; thus to establish that
the ${\frak g}_{{\cal L}(M)(p)}^0$-valued forms constructed in
\cite{M}, \cite{GM} are indeed connections, one needs to find a
transformation law only with respect to a very small group. Note that
although $\omega$ from Theorem 1.1 is not a Cartan connection,
it follows from the proof that its transformation law
is  in fact that of a Cartan connection if one acts not by the whole group
$G^2$ on $P^2$, but only by its subgroup containing linear
automorphisms given by diagonal matrices as in \cite{M},
\cite{GM}.\footnote{After this work had been completed, we discovered
the preprint \cite{CS}. The results in \cite{CS} seem to imply the
existence of a $CR$-invariant Cartan connection for hyperbolic and
elliptic $CR$-manifolds. The construction in \cite{CS} is,
however, much less explicit than ours.}
\medskip\\

{\bf 5.} Are chains as defined in Section 4 the
projections of some of submanifolds of $P^2$ that are tangent to chain
distribution (3.4) at some point?
\medskip\\

{\bf 6.} It can be shown that if the chain distribution is
integrable, then the manifold is semi-flat. Therefore, it is natural
to ask the following question: suppose that, for a hyperbolic or
elliptic real-analytic manifold $M$, all  chains (that a priori are
chains only at a single point) turn
out to be chains at each point; is it then true that $M$ is locally
equivalent to a matrix surface?

\bigskip

{\obeylines
Department of Pure Mathematics
The University of Adelaide
Box 498, G.P.O. Adelaide
South Australia 5001
AUSTRALIA
E-mail address: vezhov@spam.maths.adelaide.edu.au
\smallskip

Centre for Mathematics and Its Applications 
The Australian National University 
Canberra, ACT 0200
AUSTRALIA 
E-mail address: Alexander.Isaev@anu.edu.au
\smallskip

Mathematische Institut
Universit\"at Bonn
Beringstra{\ss}e 1
D-53115, Bonn
GERMANY
E-mail address: schmalz@uni-bonn.de
}


\begin{thebibliography}{ABCDE}

\bibitem[AH]{AH} Andreotti, A. and Hill, C. D., Complex characteristic
coordinates and tangential Cauchy-Riemann equations, {\it
Ann. Scuola Norm. Sup. Pisa, Sci. Fis. Mat., Ser III} 26(1972), 299--234.

\bibitem[B1]{B1} Beloshapka, V., A uniqueness theorem for
automorphisms of a nondegenerate surface in a complex space
(translated from Russian), {\it Math. Notes} 47(1990), 239--242.

\bibitem[B2]{B2} Beloshapka, V., On holomorphic transformations of a
quadric (translated from Russian), {\it Math. USSR. Sb.} 72(1992), 189--205.

\bibitem[BS]{BS} Burns, D. and Shnider, S., Real hypersurfaces in
complex manifolds, {\it Several Complex variables (Proc. Symp. Pure
Math. Vol XXX, Part 2, Williams Coll, Williamstown, Mass., 1975)},
141--168, {\it Amer. Math. Soc. 1977}. 

\bibitem[Ca]{Ca} Cartan, \'E., Sur la g\'eometrie pseudo-conforme
des hypersurfaces de l'espace de deux variables complexes: I, {\it
Ann. Math. Pura Appl.} 11(1932), 17--90; II, {\it Ann. Scuola
Norm. Sup. Pisa} 1(1932), 333--354.

\bibitem[Ch]{Ch} Chirka, E. M., Introduction to the geometry of
$CR$-manifolds (translated from Russian), {\it Russian Math. Surveys}
46(1991), 95--197.

\bibitem[\v CS]{CS} \v Cap, A. and Schichl, H., Parabolic geometries
and canonical Cartan connections, ESI Preprint 450, 1997. 
                    
\bibitem[CM]{CM}  Chern, S.\ S. and Moser, J., Real hypersurfaces
in complex manifolds, {\it Acta Math.} 133(1974), 219--271.

\bibitem[EI]{EI} Ezhov, V. V. and Isaev, A. V., Canonical isomorphism
of two Lie algebras arising in $CR$-geometry, Preprint CMA, ANU, 1998,
No. MRR 003-98.

\bibitem[ES1]{ES1} Ezhov, V. and Schmalz, G., Holomorphic
automorphisms of quadrics, {\it Math. Z.} 216(1994), 453--470.

\bibitem[ES2]{ES2} Ezhov, V. and Schmalz G., Normal form and
two-dimensional chains on an elliptic surface of codimension 2, {\it
J. Geom. Analysis}, to appear.

\bibitem[ES3]{ES3} Ezhov, V. and Schmalz G., Normal form and
holomorphic invariants of a parabolic 2-codimensional $CR$-surface, {\it in preparation}.

\bibitem[ES4]{ES4} Ezhov, V. and Schmalz G., Canonical form of a
strictly pseudoconvex 4-codimensional $CR$-surface in $\CC^6$, {\it in
preparation}.

\bibitem[ES5]{ES5} Ezhov, V. and Schmalz, G., Holomorphic
automorphisms of non-degenerate $CR$-quadrics: explicit description,
{\it J. Geom. Analysis}, to appear.

\bibitem[ES6]{ES6} Ezhov, V. and Schmalz, G., X-Starrheit hermitischer
Quadriken in allgemeiner Lage, {\it Math. Nachrichten}, to appear.

\bibitem[ES7]{ES7} Ezhov, V. and Schmalz, G., A matrix Poincar\'e
formula for holomorphic automorphisms of quadrics of higher
codimension. Real associative quadrics, {\it J. Geom. Analysis}, to appear. 

\bibitem[F]{F} Forstneri\v{c}, F., Mappings of quadric Cauchy-Riemann
manifolds, {\it Math. Ann.} 292(1992), 163--180.

\bibitem[GM]{GM} Garrity, T. and Mizner, R., The equivalence problem
for higher-codimensional CR structures, {\it Pacific. J. Math.}
177(1997), 211--235.

\bibitem[J]{J} Jacobowitz, H., Induced connections on hypersurfaces in
$\CC^{n+1}$, {\it Invent. Math.} 43(1977), 109--123.

\bibitem[KT]{KT} Khenkin, G. and Tumanov, A., Local characterization of
holomorphic automorphisms of Siegel domains (translated from Russian),
{\it Funct. Anal. Appl.} 17(1983), 285--294.

\bibitem[K]{K} Kuranishi, M., $CR$ geometry and Cartan geometry, {\it
Forum Math.} 7(1995), 147--205.

\bibitem[La]{La} Lai, H.-F., Real submanifolds of codimension two in
complex manifolds, {\it Trans. Amer. Math. Soc.} 264(1981), 331--352.

\bibitem[Lo1]{Lo1} Loboda, A., Real-analytic generating manifolds of
codimension 2 in $\CC^4$ and their biholomorphic mappings (translated
from Russian), {\it Math. USSR. Izv.} 33(1989), 295--315.

\bibitem[Lo2]{Lo2} Loboda, A., Linearizability of holomorphic mappings
of generating manifolds of codimension 2 in $\CC^4$ (translated from
Russian), {\it Math. USSR. Izv.} 36(1991) 655--667.

\bibitem[M]{M} Mizner, R., $CR$ structures of codimension 2, {\it
J. Diff. Geom.}, 30(1989), 167--190.

\bibitem[Sa]{Sa} Satake, I., {\it Algebraic Structures of Symmetric
Domains}, Kan\^o Memorial Lectures 4, Iwanami Shoten, Tokyo; Princeton
University Press, 1980.
 
\bibitem[St]{St} Sternberg, S., {\it Lectures on Differential Geometry},
Prentice Hall, 1964.

\bibitem[Ta1]{Ta1} Tanaka, N., On generalized graded Lie algebras and
geometric structures I, {\it J. Math. Soc. Japan}, 19(1967), 215--254.

\bibitem[Ta2]{Ta2} Tanaka, N., On non-degenerate real hypersurfaces,
graded Lie algebras and Cartan connections, {\it Japan. J. Math.},
2(1976), 131--190.

\bibitem[Ta3]{Ta3} Tanaka, N., On the equivalence problem associated
with simple graded Lie algebras, {\it Hokkaido Math. J.}, 8(1979), 23--84.

\bibitem[Tu1]{Tu1} Tumanov, A.,
Geometry of $CR$-manifolds (translated from Russian), {\it Encycl. Math. Sci.} 9 -- Several Complex Variables III,
Springer-Verlag, 1989, 201--221.

\bibitem[Tu2]{Tu2} Tumanov, A., Finite-dimensionality of the group of
$CR$ automorphisms of a standard $CR$ manifold, and proper holomorphic
mappings of Siegel domains (translated from Russian), {\it
Math. USSR. Izv.} 32(1989), 655--662.

\bibitem[V]{V} Vitushkin, A. G., Holomorphic mappings and the geometry
of hypersurfaces (translated from Russian), {\it Several Complex
Variables I. Introduction to Complex Analysis, Encycl. Math. Sci.}
7(1990), 159--214, Springer-Verlag.  
  

\end{thebibliography}
\end{document}